\documentstyle[11pt]{article}

\input amssym
\catcode`\@=11\relax
\newwrite\@unused
\def\typeout#1{{\let\protect\string\immediate\write\@unused{#1}}}

\def\psglobal#1{
\immediate\special{ps: plotfile #1 }}
\def\psfiginit{\typeout{psfiginit}
\immediate\psglobal{figtex.pro}%
\special{ps:: /TeXMagnification {\the\mag} def}
}

\def\@nnil{\@nil}
\def\@empty{}
\def\@psdonoop#1\@@#2#3{}
\def\@psdo#1:=#2\do#3{\edef\@psdotmp{#2}\ifx\@psdotmp\@empty \else
    \expandafter\@psdoloop#2,\@nil,\@nil\@@#1{#3}\fi}
\def\@psdoloop#1,#2,#3\@@#4#5{\def#4{#1}\ifx #4\@nnil \else
       #5\def#4{#2}\ifx #4\@nnil \else#5\@ipsdoloop #3\@@#4{#5}\fi\fi}
\def\@ipsdoloop#1,#2\@@#3#4{\def#3{#1}\ifx #3\@nnil
       \let\@nextwhile=\@psdonoop \else
      #4\relax\let\@nextwhile=\@ipsdoloop\fi\@nextwhile#2\@@#3{#4}}
\def\@tpsdo#1:=#2\do#3{\xdef\@psdotmp{#2}\ifx\@psdotmp\@empty \else
    \@tpsdoloop#2\@nil\@nil\@@#1{#3}\fi}
\def\@tpsdoloop#1#2\@@#3#4{\def#3{#1}\ifx #3\@nnil
       \let\@nextwhile=\@psdonoop \else
      #4\relax\let\@nextwhile=\@tpsdoloop\fi\@nextwhile#2\@@#3{#4}}
\def\psdraft{
	\def\@psdraft{0}
	\def\@psdraftspecial{100}
}
\def\psdraftspecial{
	\def\@psdraft{0}
	\def\@psdraftspecial{0}
}
\def\psfull{
	\def\@psdraft{100}
}
\psfull

\newif\if@prologfile
\newif\if@postlogfile
\newif\if@bbllx
\newif\if@bblly
\newif\if@bburx
\newif\if@bbury
\newif\if@height
\newif\if@width
\newif\if@rheight
\newif\if@rwidth
\newif\if@clip
\newif\if@right
\newif\if@left
\newif\if@toplines
\newif\if@box
\newif\if@caption
\newif\if@surround
\newif\if@captionwidth
\newif\if@captionwrite
\newif\if@captionopen
\def\@p@@sclip#1{\@cliptrue}
\def\@p@@sfile#1{
		\def\@p@sfile{#1}
}
\def\@p@@sfigure#1{
		\def\@p@sfile{#1}
}
\def\@p@sfake{\hbox to 0pt{\hss Whatever\hss}}
\def\@p@@sbbllx#1{
		\@bbllxtrue
		\@d@mscratch=#1
		\edef\@p@sbbllx{\number\@d@mscratch}
}
\def\@p@@sbblly#1{
		\@bbllytrue
		\@d@mscratch=#1
		\edef\@p@sbblly{\number\@d@mscratch}
}
\def\@p@@sbburx#1{
		\@bburxtrue
		\@d@mscratch=#1
		\edef\@p@sbburx{\number\@d@mscratch}
}
\def\@p@@sbbury#1{
		\@bburytrue
		\@d@mscratch=#1
		\edef\@p@sbbury{\number\@d@mscratch}
}
\def\@p@@sheight#1{
		\@heighttrue
		\@d@mscratch=#1
   		\edef\@p@sheight{\number\@d@mscratch}
}
\def\@p@@swidth#1{
		\@widthtrue
		\@d@mscratch=#1
		\edef\@p@swidth{\number\@d@mscratch}
}
\def\@p@@srheight#1{
		\@rheighttrue
		\@d@mscratch=#1
		\edef\@p@srheight{\number\@d@mscratch}
}
\def\@p@@srwidth#1{
		\@rwidthtrue
		\@d@mscratch=#1
		\edef\@p@srwidth{\number\@d@mscratch}
}
\def\@p@@sright#1{\@righttrue \@surroundtrue}
\def\@p@@sleft#1{\@lefttrue \@surroundtrue}
\def\@p@@sextraheight#1{\@d@mextraheight=#1}
\def\@p@@sbox#1{\@boxtrue}
\def\@p@@scaption#1{\@captiontrue}
\def\@p@@stoplines#1{
		\@toplinestrue
		\@c@ttoplines=#1
}
\def\@p@@scaptionwidth#1{
		\@captionwidthtrue
	  	\@d@mcaptionwidth=#1
}
\def\@p@@scaptionwrite#1{
		\global\@captionwritetrue
		\global\@w@rname=\expandafter{\jobname_captions.tex}
		\typeout{Captions are written to \the\@w@rname}
}

\def\@p@@sprolog#1{\@prologfiletrue\def\@prologfileval{#1}}
\def\@p@@spostlog#1{\@postlogfiletrue\def\@postlogfileval{#1}}
\def\@cs@name#1{\csname #1\endcsname}
\def\@setparms#1=#2,{\@cs@name{@p@@s#1}{#2}}
\def\ps@init@parms{
		\@bbllxfalse \@bbllyfalse
		\@bburxfalse \@bburyfalse
		\@heightfalse \@widthfalse
		\@rheightfalse \@rwidthfalse
		\def\@p@sbbllx{}\def\@p@sbblly{}
		\def\@p@sbburx{}\def\@p@sbbury{}
		\def\@p@sheight{}\def\@p@swidth{}
		\def\@p@srheight{}\def\@p@srwidth{}
		\def\@p@sfile{}
		\def\@p@scost{10}
		\def\@sc{}
		\@prologfilefalse
		\@postlogfilefalse
		\@clipfalse
		\@rightfalse \@leftfalse
		\@boxfalse \@captionfalse
		\@toplinesfalse \@surroundfalse
		\@d@mextraheight=0pt
 		\@c@ttoplines=0
		\@pshape={} \def\@p@srheight@total{}
		\@captionwidthfalse \@d@mcaptionwidth=0pt
}
\def\parse@ps@parms#1{
	 	\@psdo\@psfiga:=#1\do
		   {\expandafter\@setparms\@psfiga,}}
\newif\ifno@bb
\newif\ifnot@eof
\newread\ps@stream
\newtoks\@linetok
\def\bb@missing{
	\typeout{psfig: searching \@p@sfile \space  for bounding box}
	\openin\ps@stream=\@p@sfile
	\no@bbtrue
	\not@eoftrue
	\catcode`\%=12
	\loop
		\read\ps@stream to \line@in
		\global\@linetok=\expandafter{\line@in}
		\ifeof\ps@stream \not@eoffalse \fi
		\@bbtest{\@linetok}
		\if@bbmatch\not@eoffalse\expandafter\bb@cull\the\@linetok\fi
	\ifnot@eof \repeat
	\catcode`\%=14
}	
\catcode`\%=12
\newif\if@bbmatch
\def\@bbtest#1{\expandafter\@a@\the#1
\long\def\@a@#1
     \ifx\@bbtest#2\@bbmatchfalse\else\@bbmatchtrue\fi}
\long\def\bb@cull#1 #2 #3 #4 #5 {
	\@d@mscratch=#2 bp\edef\@p@sbbllx{\number\@d@mscratch}
	\@d@mscratch=#3 bp\edef\@p@sbblly{\number\@d@mscratch}
	\@d@mscratch=#4 bp\edef\@p@sbburx{\number\@d@mscratch}
	\@d@mscratch=#5 bp\edef\@p@sbbury{\number\@d@mscratch}
	\no@bbfalse
}
\catcode`\%=14
\def\compute@bb{
		\no@bbfalse
		\if@bbllx \else \no@bbtrue \fi
		\if@bblly \else \no@bbtrue \fi
		\if@bburx \else \no@bbtrue \fi
		\if@bbury \else \no@bbtrue \fi
		\ifno@bb \bb@missing \fi
		\ifno@bb \typeout{FATAL ERROR: no bb supplied or found}
			\no-bb-error
		\fi
		\count203=\@p@sbburx
		\count204=\@p@sbbury
		\advance\count203 by -\@p@sbbllx
		\advance\count204 by -\@p@sbblly
		\edef\@bbw{\number\count203}
		\edef\@bbh{\number\count204}
}
\def\in@hundreds#1#2#3{\count240=#2 \count241=#3
		     \count100=\count240	
		     \divide\count100 by \count241
		     \count101=\count100
		     \multiply\count101 by \count241
		     \advance\count240 by -\count101
		     \multiply\count240 by 10
		     \count101=\count240	
		     \divide\count101 by \count241
		     \count102=\count101
		     \multiply\count102 by \count241
		     \advance\count240 by -\count102
		     \multiply\count240 by 10
		     \count102=\count240	
		     \divide\count102 by \count241
		     \count200=#1\count205=0
		     \count201=\count200
			\multiply\count201 by \count100
		     	\advance\count205 by \count201
		     \count201=\count200
			\divide\count201 by 10
		     	\multiply\count201 by \count101
			\advance\count205 by \count201
		     \count201=\count200
			\divide\count201 by 100
			\multiply\count201 by \count102
			\advance\count205 by \count201
		     \edef\@result{\number\count205}
}
\def\compute@wfromh{
		\in@hundreds{\@p@sheight}{\@bbw}{\@bbh}
		\edef\@p@swidth{\@result}
}
\def\compute@hfromw{
		\in@hundreds{\@p@swidth}{\@bbh}{\@bbw}
		\edef\@p@sheight{\@result}
}
\def\compute@handw{
		\if@height
			\if@width
			\else
				\compute@wfromh
			\fi
		\else
			\if@width
				\compute@hfromw
			\else
				\edef\@p@sheight{\@bbh}
				\edef\@p@swidth{\@bbw}
			\fi
		\fi
}
\def\compute@resv{
		\if@rheight \else \edef\@p@srheight{\@p@sheight} \fi
		\if@rwidth \else \edef\@p@srwidth{\@p@swidth} \fi
		\edef\@p@srheight@total{\@p@srheight}
}
\newtoks\@pshape
\def\@c@ttoplines{\count120}
\def\@c@theightcount{\count121}
\def\@c@tshapecount{\count122}
\newdimen\@d@mwidthshape
\newdimen\@d@mextraheight
\newdimen\@d@mscratch
\def\compute@parshape{
	\if@right
		\if@left
	   		\typeout{error: Can't have both left and right set}
			\@leftfalse
		\fi
	\fi
	\@d@mscratch=\@p@swidth truesp
	\divide \@d@mscratch by 19
	\multiply \@d@mscratch by 20
	\edef\@p@swidthdimen{\the\@d@mscratch}
	\@c@tshapecount=\@c@ttoplines
 	\@d@mscratch=\baselineskip
	\multiply \@d@mscratch by \@c@ttoplines
	\advance \@d@mscratch by .4\baselineskip
    	\edef\@p@stopdistance{\the\@d@mscratch }
	\@d@mscratch=\@p@sheight truesp
	\divide \@d@mscratch by 2
	\edef\@p@shalfboxheight{\the\@d@mscratch}
	\if@toplines
		\loop \@pshape=\expandafter{\the\@pshape 0pt \hsize}
		\advance\@c@ttoplines by -1
		\ifnum\@c@ttoplines>0 \repeat
	\fi
   	\@c@theightcount=\@p@srheight@total
	\advance \@c@theightcount by \@d@mextraheight
	\divide  \@c@theightcount by \baselineskip
	\advance \@c@theightcount by 1
    	\advance \@c@tshapecount by \@c@theightcount
	\advance \@c@theightcount by -1
	\@d@mwidthshape=\hsize
     	\advance \@d@mwidthshape by -\@p@swidthdimen
	\if@left
		\def\@moveshape{0pt}
		\ifnum\@c@theightcount>0
		      	\loop
			\@pshape=%
			\expandafter{\the\@pshape %
					\@p@swidthdimen \@d@mwidthshape}
			\advance \@c@theightcount by -1
			\ifnum\@c@theightcount>0 \repeat
		\else
			\advance \@c@tshapecount by 1
		\fi
	\fi
	\if@right
		\@d@mscratch=\hsize
		\advance \@d@mscratch by -\@p@swidth truesp
		\edef\@moveshape{\@d@mscratch}
		\ifnum\@c@theightcount>0
			\loop
			\@pshape=\expandafter{\the\@pshape 0pt \@d@mwidthshape}
			\advance \@c@theightcount by -1
			\ifnum\@c@theightcount>0 \repeat
		\else
			\advance \@c@tshapecount by 1
		\fi
	\fi
	\ifnum \@p@srheight=0
		\@pshape={}
		\@c@tshapecount = 0
	\else
	 	\@pshape=\expandafter{\the\@pshape 0pt \hsize}
	\fi
}
\def\@p@ssetsurroundboxes{
	\global\parshape=\count122 \the\@pshape		
 	\moveright\@moveshape
	\vbox to 0pt\bgroup\hskip0pt\vskip\@p@stopdistance
}
\newtoks\@captiontok
\newbox\@b@xcaption
\newdimen\@d@mcaptionwidth
\newdimen\@d@mcaptionheight
\newwrite\@w@rcaption
\newtoks\@w@rname
\def\setcaption#1{\@captiontok={#1}}
\def\@set@caption{
	\setbox\@b@xcaption\vbox{\hsize\@d@mcaptionwidth
	\tolerance=9000 \baselineskip=11.4pt
	\noindent\relax\the\@captiontok}
	\if@captionwrite
		\if@captionopen
		\else
			\global\@captionopentrue
			\immediate\openout\@w@rcaption=\the\@w@rname
		\fi
		\immediate\write\@w@rcaption{\the\@captiontok}
		\immediate\write\@w@rcaption{}
	\fi
}
\def\compute@caption{
	\if@captionwidth
	\else
		\@d@mcaptionwidth = \@p@swidth truesp
		\divide \@d@mcaptionwidth by 20
		\multiply \@d@mcaptionwidth by 17
	\fi
	\@set@caption
	\@d@mcaptionheight=\ht\@b@xcaption
	\if@rheight
	\else
		\count100=\@p@srheight
	   	\advance \count100 by \@d@mcaptionheight
	   	\advance \count100 by \bigskipamount
	   	\advance \count100 by \medskipamount
	   	\edef\@p@srheight@total{\number\count100}
	\fi
}
\newif\if@alreadyjtem \@alreadyjtemfalse
\def\newpar{\hfil\vadjust{\vskip\parskip}%
	{\count100=\parskip
	\count101=\baselineskip
	\divide\count101 by 10  \multiply\count101 by 3
	\advance \count100 by \count101
	\divide\count100 by \baselineskip
	\advance\count100 by \prevgraf
	\global\prevgraf=\count100}%
	\break\if@alreadyjtem\else\indent\fi%
}
\let\sav@par=\par
\def\jtem#1{%
    	\if@alreadyjtem\else\bgroup\fi
	\def\par{\sav@par\egroup\sav@par}
	\if@alreadyjtem\else\leftskip\parindent\fi
	\@alreadyjtemtrue
	\noindent\hskip0pt
	\llap{#1\ }\ignorespaces
}
\def\compute@sizes{%
	\compute@bb
	\compute@handw
  	\compute@resv
	\if@caption
		\compute@caption
	\fi
	\if@surround
		\compute@parshape
	\fi
}
\def\@p@sdospecials{
	\ifnum\@p@scost<\@psdraft
	       	\typeout{psfig: including \@p@sfile \space }
	\fi
	\special{ps::[begin] 	\@p@swidth \space \@p@sheight \space
			\@p@sbbllx \space \@p@sbblly \space
			\@p@sbburx \space \@p@sbbury \space
			startTexFig \space }
	\ifnum\@p@scost<\@psdraft
		\if@clip
			\typeout{(clip)}
			\special{ps:: \@p@sbbllx \space \@p@sbblly \space
				\@p@sbburx \space \@p@sbbury \space
			    	doclip \space }
		\fi
	\fi
	\if@box
		\typeout{(box)}
  		\special{ps:: \@p@sbbllx \space \@p@sbblly \space
			\@p@sbburx \space \@p@sbbury \space
		    	dobox \space }
	\fi
	\ifnum\@p@scost<\@psdraft
		\if@prologfile
	    		\special{ps: plotfile \@prologfileval \space }
		\fi
		\special{ps: plotfile \@p@sfile \space }
    		\if@postlogfile
			\special{ps: plotfile \@postlogfileval \space }
		\fi
	\fi
	\special{ps::[end] endTexFig \space }
}
\newif\if@putinvbox

\def\psfig#1{{%
	\ifhmode%
		\vbox\bgroup
		\@putinvboxtrue
	\else
		\@putinvboxfalse
	\fi
       	\ps@init@parms
	\parse@ps@parms{#1}
       	\compute@sizes
	\if@surround
		\psfig@for@surround
	\else
		\psfig@for@regular
	\fi
	\if@putinvbox
       		\egroup
	\fi
}}
\def\psfig@for@surround{%
	\@p@ssetsurroundboxes
	\ifnum\@p@scost<\@psdraft
		\@p@sdospecials
		\vbox to \@p@srheight true sp{\vss}
       	\else
		\if@box
			\@p@sdospecials
		\fi
		\vbox to \@p@srheight true sp{
			\vskip\@p@shalfboxheight
			\hbox to \@p@srwidth true sp{
				\hss
				\ifnum\@p@scost<\@psdraftspecial
					\@p@sfile
				\else
					\@p@sfake
				\fi
      				\hss
			}
		\vss
		}
	\fi
	\if@caption
		\medskip
		\hbox to \@p@srwidth true sp{
			\hss
			\box\@b@xcaption
			\hss
		}
 		\medskip
	\fi
	\vss\egroup
	\vskip-\parskip
}

\def\psfig@for@regular{%
	\if@putinvbox
	\else
		\vskip\parskip
	\fi
	\ifnum\@p@scost<\@psdraft
		\@p@sdospecials
		\vbox to \@p@srheight true sp{%
			\hbox to \@p@srwidth true sp{
			\hfil
			}
		\vfil
		}
       	\else
		\if@box
			\@p@sdospecials
		\fi
	    	\vbox to \@p@srheight true sp{
			\vss
			\hbox to \@p@srwidth true sp{
				\hss
				\ifnum\@p@scost<\@psdraftspecial
					\@p@sfile
				\else
					\@p@sfake
				\fi
				\hss
			}
		    	\vss
		}
	\fi
	\if@caption
		\medskip
		\hbox to \@p@srwidth true sp{
			\hss
			\box\@b@xcaption
			\hss
		}
		\bigskip
	\fi
	\if@putinvbox
	\else
		\vskip-\parskip
	\fi
}
\catcode`\@=12\relax

\psfiginit

\font\tinybbfont=msbm6
\font\scriptsizebbfont=msbm7 scaled \magstep 1
 1
\font\footnotesizebbfont=msbm9 scaled \magstep 0
\font\smallbbfont=msbm7 scaled \magstep 2
\font\bbfont=msbm9 scaled \magstep1  
\font\largebbfont=msbm10 scaled \magstep 1
 2
 3
 4

\def\tinyBbb#1{\hbox{\tinybbfont #1}}
\def\scriptsizeBbb#1{\hbox{\scriptsizebbfont #1}}

\def\footnotesizeBbb#1{\hbox{\footnotesizebbfont #1}}
\def\smallBbb#1{\hbox{\smallbbfont #1}}
\def\Bbb#1{\hbox{\bbfont #1}}
\def\largeBbb#1{\hbox{\largebbfont #1}}

\newcommand{\CP}{\mbox{{\Bbb C}{\rm P}}}
\newcommand{\CPlarge}{\mbox{{\largeBbb C}{\rm P}}}
\newcommand{\CPscriptsize}{\mbox{\scriptsize {\scriptsizeBbb C}{\rm P}}}
\newcommand{\CPsmall}{\mbox{\small {\smallBbb C}{\rm P}}}

\newcommand{\Det}{\mbox{\it Det}\,}

\newcommand{\Diag}{\mbox{\it Diag}\,}

\newcommand{\Fl}{\mbox{\it Fl}\,}
\newcommand{\GL}{\mbox{\it GL}\,}
\newcommand{\Gr}{\mbox{\it Gr}}
\newcommand{\HG}{\mbox{\it HG}\,}
\newcommand{\Hom}{\mbox{\it Hom}\,}
\newcommand{\Homsheaf}{\mbox{\it ${\cal H}$om}\,}
\newcommand{\HQ}{\mbox{${\cal H}$${\cal Q}$}}
\newcommand{\HQuot}{\mbox{\it HQuot}\,}

\newcommand{\Ker}{\mbox{\it Ker}\,}

\newcommand{\Mor}{\mbox{\it Mor}}

\newcommand{\Proj}{\mbox{\rm Proj}\,}

\newcommand{\Quot}{\mbox{\it Quot}\,}

\newcommand{\Spec}{\mbox{\it Spec}\,}

\newcommand{\SU}{\mbox{\it SU}\,}
\newcommand{\Sym}{\mbox{\it Sym}}

\newcommand{\degree}{\mbox{\it deg}\;}
\newcommand{\dimm}{\mbox{\it dim}\,}
\newcommand{\ev}{\mbox{\it ev}\,}

\newcommand{\pr}{\mbox{\rm pr}}
\newcommand{\pt}{\mbox{\it pt}}
\newcommand{\rank}{\mbox{\it rank}\,}

\begin{document}

\enlargethispage{23cm}

\begin{titlepage}

$ $

\vspace{-2cm} 

\noindent\hspace{-1cm}
\parbox{6cm}{\small January 2004}\
   \hspace{6.5cm}\
   \parbox{5cm}{math.AG/0401367}

\vspace{1.5cm}

\centerline{\large\bf
 $S^1$-fixed-points in hyper-Quot-schemes
} \vspace{1ex}
\centerline{\large\bf
 and an exact mirror formula for flag manifolds
} \vspace{1ex}
\centerline{\large\bf
  from the extended mirror principle diagram
} 

\vspace{1.5cm}

\centerline{\large
  Chien-Hao Liu, \hspace{1ex}
  Kefeng Liu,    \hspace{1ex} and \hspace{1ex}
  Shing-Tung Yau
 }

\vspace{3em}

\begin{quotation}
\centerline{\bf Abstract}
\vspace{0.3cm}
\baselineskip 12pt  
{\small
 In the series of work [L-L-Y1, III: Sec.\ 5.4] on mirror principle,
  two of the current authors (K.L.\ and S.-T.Y.) with Bong H.\ Lian
  developed a method to compute the integral
  $\int_{X}\tau^{\ast}e^{H\cdot t}\cap {\mathbf 1}_d$
  for a flag manifold $X=\Fl_{r_1,\,\ldots,\,r_I}({\smallBbb C}^n)$
  via an extended mirror principle diagram.
 This integral determines the fundamental hypergeometric series
  $HG[{\mathbf 1}]^X(t)$ and is also related to the computation of
  the Gromov-Witten invariants (string world-sheet instanton numbers)
  on $X$.
 This method turns the required localization computation on the augmented
  moduli stack $\overline{\cal M}_{0,0}(\CPsmall^1\times X)$ of stable
  maps to a localization computation on a hyper-Quot-scheme
  $\mbox{\small\it HQuot}({\cal E}^n)$ of inclusion sequences of
  subsheaves of a trivialized trivial bundle ${\cal E}^n$ of rank $n$
  on $\CPsmall^1$.
 In this article, the detail of this localization computation
  on $\mbox{\small\it HQuot}({\cal E}^n)$ is carried out.
 The necessary ingredients in the computation, notably,
  the $S^1$-fixed-point components and the distinguished ones
   $E_{(A;0)}$ in $\mbox{\small\it HQuot}({\cal E}^n)$,
  the $S^1$-equivariant Euler class of $E_{(A;0)}$ in
   $\mbox{\small\it HQuot}({\cal E}^n)$, and
  a push-forward formula of cohomology classes involved in the problem
   from the total space of a restrictive flag manifold bundle to its
   base manifold are given.
 With these, an exact expression of
  $\int_{X}\tau^{\ast}e^{H\cdot t}\cap {\mathbf 1}_d$ is obtained.
 When $X$ is a Grassmannian manifold, the same route reproduces
  the known exact expression for $HG[{\mathbf 1}]^X(t)$.
 For a general flag manifold $X$, our expression determines
  $\HG[{\mathbf 1}]^X(t)$ implicitly.
 Remarks on what it suggests for general Hori-Vafa formula are given.
 Due to the technical necessity, a discussion on the general construction
  of restrictive flag manifold bundle,
  its natural embedding in a flag manifold bundle, and
  the Thom class of this embedding is also given.
 This work generalizes the result in [L-L-L-Y]. This work gives\
 explicit formulas for mirror principle computations of Calabi-Yau
 manifolds in flag manifolds.
} 
\end{quotation}

\bigskip

\baselineskip 12pt
{\footnotesize
 } 

\medskip

\noindent {\small
MSC number 2000$\,$: 14N35, 14D20, 55N91, 55R91, 81T30.
} 

\bigskip

\baselineskip 11pt
{\footnotesize
\noindent{\bf Acknowledgements.}
 We would like to thank
  Sheldon Katz,
  Jun Li,
  Bong H.\ Lian, and
  Cumrun Vafa
 for motivations, discussions, communications, and comments on this work.
 C.-H.L.\ would like to thank in addition
  Kentaro Hori,
  Shiraz Minwalla,
  Barton Zwiebach, and
  TASI 2003 String Theory Program
 for courses and discussions on string theory that influence the work;
  IPAM-UCLA Symplectic Geometry and String Theory Program, June 2003,
 for an invitation of talk related to the work;
  Hungwen Chang and Ann Willman
 for reminding him of other reference frames in life.
 A large portion of the draft is written in C.-H.L.'s library set up
  by Ling-Miao Chou; no words can express his thanks and moral debts
  to her.
 The work is supported by NSF grants DMS-9803347 and DMS-0074329.
} 

%

\end{titlepage}

\newpage
$ $

\vspace{-4em}  

\centerline{\sc
 Hyper-Quot-Schemes and Exact Mirror Formula}

\vspace{2em}

\baselineskip 14pt  

\begin{flushleft}
{\Large\bf 0. Introduction and outline.}
\end{flushleft}

\begin{flushleft}
{\bf Introduction.}
\end{flushleft}
In the series of work [L-L-Y1, III: Sec.\ 5.4] on Mirror Principle,
 two of the current authors (K.L.\ and S.-T.Y.) with Bong H.\ Lian
 developed a method to compute the integral
 $\int_{X}\tau^{\ast}e^{H\cdot t}\cap {\mathbf 1}_d$
 for a flag manifold $X=\Fl_{r_1,\,\ldots,\,r_I}({\smallBbb C}^n)$
 via an extended mirror principle diagram,
 cf.\ Sec.\ 1 and Sec.\ 3.1.

This integral determines the fundamental hypergeometric series
 $HG[{\mathbf 1}]^X(t)$ and is also related to the computation of
 the Gromov-Witten invariants (string world-sheet instanton numbers)
 on $X$.
This method turns the required localization computation on the augmented
 moduli stack $\overline{\cal M}_{0,0}(\CP^1\times X)$ of stable
 maps to a localization computation on a hyper-Quot-scheme
 $\HQuot({\cal E}^n)$ of inclusion sequences of subsheaves of
 a trivialized trivial bundle ${\cal E}^n$ of rank $n$ on $\CP^1$.
The major purpose of the current work is to carry out the full detail
 of this localization computation on $\HQuot({\cal E}^n)$.

As a first step, the $S^1$-fixed-point components in $\HQuot({\cal E}^n)$
 are described and the distinguished ones $E_{(A;0)}$ are identified,
 cf.\ Sec.\ 2.1, Sec.\ 2.2, and Sec.\ 3.1.
In particular, $E_{(A;0)}$ admits a tower of fibrations with fiber
 restrictive flag manifolds.
Also, by construction, there are canonical morphisms from all these
 $E_{(A;0)}$ to the flag manifold $X$.

For technical necessity, we study a general construction of a restrictive
 flag manifold bundle $W$ over a base manifold $Y$ and its associated
 flag manifold bundle $W^{\prime}$ over the same base,
 cf.\ Sec.\ 3.2.
$W$ naturally embeds in $W^{\prime}$ and we work out the Thom class
 of $W$ in $W^{\prime}$ with respect to this embedding.

The second step involves the computation of the $S^1$-equivariant Euler
 class of the normal bundle $\nu(E_{(A;0)}/HQuot({\cal E}^n))$ of
 $E_{(A;0)}$ in $\HQuot({\cal E}^n)$.
This involves both the understanding of the restriction of the tangent
 bundle $T_{\ast}\HQuot({\cal E}^n)$ of the hyper-Quot-scheme to 
 $E_{(A;0)}$ and the tangent bundle of $E_{(A;0)}$.
Some deformation-theoretical aspects of these spaces and their
 natural decompositions in the $K$-group of $E_{(A;0)}$ are studied.
The computation of the $S^1$-equivariant Euler class of
 $\nu(E_{(A;0)}/\HQuot({\cal E}^n))$ then follows.
Cf.\ Sec.\ 3.3 and Sec.\ 3.5.

The third step involves a push-forward formula of the cohomology classes
 involved in the problem from the total space of a restrictive flag
 manifold bundle to its base manifold.
Consecutive applications of this push-forward formula via the tower
 fibration of $E_{(A;0)}$ give rise to an exact expression of
 the integral $\int_{X}\tau^{\ast}e^{H\cdot t}\cap {\mathbf 1}_d$.
Cf.\ Sec.\ 3.4 and Sec.\ 3.6.

For a general flag manifold, our expression can be interpreted as
 arising from the fundamental hypergeometric series for a product of
 Grassmannian manifolds that contains the flag manifold combined with
 the effect of the Thom class of the induced inclusion of
 $\HQuot({\cal E}^n)$ in a product of Quot-schemes.
When the flag manifold is a Grassmannian manifold, the same route
 reproduces the known expression of $\HG[{\mathbf 1}]^X(t)$ in [B-CF-K]
 and hence the Hori-Vafa formula, conjectured in [H-V] and
 studied also in [B-CF-K] following [L-L-L-Y],
 for the case Grassmannian manifolds.
Cf.\ Sec.\ 4.

The current work generalizes the results in [L-L-L-Y].

\bigskip

\bigskip

\begin{flushleft}
{\bf Outline.}
\end{flushleft}
{\small
\baselineskip 11pt  
\begin{quote}
 1. Essential background and notations for physicists.

 2. The $S^1$-fixed-point components on hyper-Quot-schemes.
    \vspace{-1ex}
    \begin{quote}
     \hspace{-1.3em}
     2.1 \ \parbox[t]{12cm}{
         Inclusion pairs of $S^1$-invariant subsheaves of ${\cal E}^n$.}

     \hspace{-1.3em}
     2.2 \ \parbox[t]{12cm}{
         The $S^1$-fixed-point locus in $HQuot_P({\cal E}^n)$. }
    \end{quote}

 \vspace{-.8ex}
 3. An exact computation of
    $\int_{X}\tau^{\ast}e^{H\cdot t}\cap {\mathbf 1}_d$
    from the mirror
    principle diagram.
    \vspace{-1ex}
    \begin{quote}
     \hspace{-1.3em}
     3.1 \ \parbox[t]{30em}{
         The extended Mirror Principle diagram and the distinguished
         \newline
         $S^1$-fixed-point components in the hyper-Quot-scheme
         ${\cal Q}_d$.}

     \smallskip
     \hspace{-1.3em}
     3.2 \ \parbox[t]{12cm}{
         Bundles with fiber restrictive flag manifolds and
         the class $\Omega({\cal P}_{\bullet})$.}

     \hspace{-1.3em}
     3.3 \ \parbox[t]{12cm}{
         Tautological sheaves on $E_{(A;0)}$ and
         $E_{(A;0)}\times\CPsmall^1$.}

     \hspace{-1.3em}
     3.4 \ \parbox[t]{12cm}{
         The hyperplane-induced classes on $E_{(A;0)}$.}

     \hspace{-1.3em}
     3.5 \ \parbox[t]{12cm}{
         An exact computation of
         $e_{S^1}(\nu(E_{(A;0)}/HQuot_P({\cal E}^n)))$.}

     \hspace{-1.3em}
     3.6 \ \parbox[t]{12cm}{
         An exact computation of
         $\int_{X}\tau^{\ast}e^{H\cdot t}\cap {\mathbf 1}_d$.  }
    \end{quote}

 \vspace{-.8ex}
 4. Remarks on the Hori-Vafa conjecture.
\end{quote}
} 

\bigskip

\baselineskip 14pt  

\bigskip

\section{Essential background and notations for physicists.}
Essential background or its main references used in this article
 and notations for objects involved are collected in this section
 for the convenience of readers.
The list extends that in [L-L-L-Y].

\bigskip

\noindent $\bullet$
{\bf Hyper-Quot-scheme.}
Fix the ample line bundle ${\cal O}_{\CPscriptsize^1}(1)$ over $\CP^1$.
Let
 ${\cal E}^n$ be the trivialized trivial bundle
  ${\cal O}_{\CPscriptsize^1}\otimes{\Bbb C}^n$ of rank $n$ over $\CP^1$,
  $P=(P_1,\,\ldots,\,P_I)$ be a finite sequence of integral polynomials
  $P_i(t)=(n-r_i)t + d_i +(n-r_i)$ with $r_1<\,\ldots\,<r_I$.
Then {\it the hyper-Quot-scheme} $\HQuot_P({\cal E}^n)$ is the fine
 moduli space that parameterizes the set of successive quotients
 $$
  {\cal E}^n \twoheadrightarrow
    {\cal E}^n\!/\mbox{\raisebox{-.4ex}{${\cal V}_1$}}\,
   \twoheadrightarrow\, \cdots\, \twoheadrightarrow
    {\cal E}^n\!/\mbox{\raisebox{-.4ex}{${\cal V}_I$}}\,
 $$
 with Hilbert polynomial
 $P({\cal E}^n\!/\mbox{\raisebox{-.4ex}{${\cal V}_i$}},t)=P_i(t)$.
It is the scheme that represents the {\it hyper}-$\Quot$-functor -
 which generalizes Grothendieck's Quot-functor -
 for ${\cal E}^n$, cf [Gr3].

\bigskip

\noindent $\bullet$
{\bf Hyper-Quot-scheme compactification of
 $\Hom(\CP^1, \Fl_{r_1,\,\ldots,\,r_I}({\Bbb C}^n))$.}
(Cf.\ [CF1], [Kim], [La], and [Str].)
 Let
  $C=\CP^1$ with the very ample line bundle
    ${\cal O}_{\CPscriptsize^1}(1)$,
  ${\cal E}^n={\cal O}_C\otimes{\Bbb C}^n$ be a trivialized trivial
   bundle of rank $n$ over $C$,
  $\Fl_{r_1,\,\ldots,\,r_I}({\Bbb C}^n)$, $r_1<\,\cdots\,<r_I$,
   be the flag manifold that parameterizes inclusion sequences
   $V_{\bullet}: V_1\hookrightarrow \,\cdots\,\hookrightarrow V_I$
   of planes $V_i$ in ${\Bbb C}^n$ of dimension $r_i$, and
  $\Hom(\CP^1, \Fl_{r_1,\,\ldots,\,r_I}({\Bbb C}))$ be the space of
   morphisms from $\CP^1$ to $\Fl_{r_1,\,\ldots,\,r_I}({\Bbb C}^n)$.
 Then an element
  $(f:\CP^1\rightarrow \Fl_{r_1,\,\ldots,\,r_I}({\Bbb C}^n))$
  in $\Hom(\CP^1,\Fl_{r_1,\,\ldots,\,r_I}({\Bbb C}))$ determines
  a unique inclusion sequence  (i.e.\ filtration of ${\cal E}^n$)
  ${\cal V}_{\bullet}:{\cal V}_1\hookrightarrow\,\cdots\,
                                \hookrightarrow {\cal V}_I$
  of subbundles ${\cal V}_i$ of rank $r_i$ in ${\cal E}^n$,
  which corresponds in turn to the element
  ${\cal E}^n
    \twoheadrightarrow {\cal E}^n\!/\mbox{\raisebox{-.4ex}{${\cal V}_1$}}
    \twoheadrightarrow \,\cdots\,
    \twoheadrightarrow {\cal E}^n\!/\mbox{\raisebox{-.4ex}{${\cal V}_I$}}$
  (i.e.\ cofiltration of ${\cal E}^n$) in $\HQuot({\cal E}^n)$.
 This gives a natural embedding of
  $\Hom(\CP^1,\Fl_{r_1,\,\ldots,\,r_I}({\Bbb C}^n))$
  in $\HQuot({\cal E}^n)$.
 The component of $\Hom(\CP^1,\Fl_{r_1,\,\ldots,\,r_I}({\Bbb C}^n))$
  that contains degree $d=(d_1,\,\ldots\,d_I)$ image curves in
  $\Fl_{r_1,\,\ldots\,,r_I}({\Bbb C}^n)$ is embedded in
  $\HQuot_P({\cal E}^n)$ with the Hilbert polynomial
  $P=(P_1(t),\,\ldots\,,P_I(t))$, where $P_i(t)=(n-r_i)t+d_i+(n-r_i)$.
 This gives a compactification of
  $\Hom(\CP^1,\Fl_{r_1,\,\ldots\,,r_I}({\Bbb C}^n))$
  via hyper-Quot-schemes, other than the moduli space
  $\overline{M}_{0,0}(\Fl_{r_1,\,\ldots,\,r_I}({\Bbb C}^n),d)$
  of stable maps from $\CP^1$ to $\Fl_{r_1,\,\ldots\,r_I}({\Bbb C}^n)$.
 Recall also that $\HQuot_P({\cal E}^n)$ is a smooth, irreducible,
  projective variety of dimension
  $\sum_{i=1}^I\,(n-r_i)(r_i-r_{i-1})+\sum_{i=1}^I\,d_i(n_{i+1}-n_{i-1})$.
 The $S^1$-action on $\CP^1$ induces an $S^1$-action on
  $\Hom(\CP^1,\Fl_{r_1,\,\ldots\,,r_I}({\Bbb C}^n))$ and
  $\HQuot_P({\cal E}^n)$ respectively.
 The two actions coincide under the natural embedding of
  $\Hom(\CP^1, \Fl_{r_1,\,\ldots\,,r_I}({\Bbb C}^n))$
  in $\HQuot({\cal E}^n)$.

\bigskip

\noindent $\bullet$
{\bf Mirror principle diagram for flag manifolds.}
For the details of Mirror Principle, readers are referred to
 [L-L-Y1$\,$: I, II, III, IV]. Some survey is given in [L-L-Y2].
To avoid digressing too far away, here we shall take
 [L-L-Y1, III$\,$: Sec.\ 5.4] as our starting point and
 restrict to the case that the target manifold of stable maps is
 $X=\Gr_r({\Bbb C}^n)$.
Recall the embedding of
 $$
  \tau:X:=\Fl_{r_1,\,\ldots\,,r_I}({\Bbb C}^n) \hookrightarrow
  Y:=\CP^{{n\choose r_1}-1}\times\,\cdots\,\times\CP^{{n\choose r_I}-1}
 $$
 induced from the Pl\"{u}cker embeddings
 $\tau_i:\Gr_{r_i}({\Bbb C}^n)\,\rightarrow\,\CP^{\,{n\choose r_i}-1}$.
$\tau$ induces an isomorphism between the divisor class groups
 $\tau^{\ast}:A^1(Y)\stackrel{\sim}{\rightarrow} A^1(X)$.

Recall next the Mirror Principle diagram for
 $X=\Fl_{r_1,\,\ldots\,,r_I}({\Bbb C}^n)$.
The geometric objects involved are contained in the following
 diagram$\,$:

{\footnotesize
$$
 \begin{array}{cccccccclcl}
  V  &  & U_d & & V_d & & {\cal U}_d        \\
  \downarrow   & & \downarrow   & & \downarrow   & & \downarrow \\
  X  & \stackrel{ev}{\longleftarrow}
     & \overline{M}_{0,1}(X,d)  & \stackrel{\rho}{\longrightarrow}
     & \overline{M}_{0,0}(X,d)  & \stackrel{\pi}{\longleftarrow}
     & M_d             & \stackrel{\varphi}{\longrightarrow}
     & W_d             & \stackrel{\psi}{\longleftarrow}
     & {\cal Q}_d:= HQuot_P{{\cal E}^n}                             \\
   & & & & & & \cup    & & \hspace{.6ex}\cup  & & \hspace{.6ex}\cup \\
   & & & & & & F_0     & \stackrel{ev^Y}{\longrightarrow}
     & Y_0\,(\supset X_0=X)   & \stackrel{g}{\longleftarrow}
     & E_0\;=\;\cup_s\,E_{0s}\,,\\
   & & & & & & \hspace{2.2em}\downarrow\,\mbox{\scriptsize $ev^X$}
     & & \mbox{\scriptsize $|$}\wr \\
   & & & & & & X  & \stackrel{\tau}{\longrightarrow} & Y
 \end{array}
$$
{\normalsize where}} 
\begin{itemize}
 \item [(1)] {\it Moduli spaces}$\,$:
  $\overline{M}_{0,0}(X,d)$ is the moduli space of genus-$0$ stable
   maps of degree $d=(d_1,\,\ldots\,,d_I)$ into $X$,
  $\overline{M}_{0,1}(X,d)$ is the moduli space of genus-$0$,
   $1$-pointed stable maps of degree $d$ into $X$,
  $M_d=\overline{M}_{0,0}(\CP^1\times X, (1,d))$,
  $W_d$ is the linearized moduli space at degree $d$,
   which can be chosen to be the product of projective spaces:
   $\prod_{i=1}^I\,
     {\Bbb P}(H^0(\CP^1,\,{\cal O}_{{\scriptsizeBbb C}{\rm P}^1}(d_i))
                                   \otimes \Lambda^{r_i}{\Bbb C}^n)$
   for $X=\Fl_{r_1,\,\ldots\,,r_I}({\Bbb C}^n)$, and
  ${\cal Q}_d=\HQuot_P({\cal E}^n)$ with
   $P=(P_1(t),\,\cdots\,,P_I(t))$, where $P_i(t)=(n-r_i)t+d_i+(n-r_i)$.

 \item [(2)] {\it Group actions}$\,:$
  there are ${\Bbb C}^{\times}$-actions on
   $M_d$, $W_d$, and ${\cal Q}_d$ respectively that are compatible
   with the morphisms among these moduli spaces;
  these ${\Bbb C}^{\times}$-actions induce $S^1$-actions on these
   moduli spaces by taking the subgroup
   $U(1)\subset{\Bbb C}^{\times}$.

 \item [(3)] {\it Morphisms}$\,$:
   $\ev$ is the evaluation map,
   $\rho$ is the forgetful map,
   $\pi$ is the contracting morphism,
   $\varphi$ is the collapsing morphism, and
   $\psi$ is an $S^1$-equivariant resolution of singularities
    of $\varphi(M_d)$.
   $\varphi$ and$\psi$ are discussed in detail in [L-L-L-Y: Sec.\ 3.1]
    when $X$ is a Grassmannian manifold.
   Their generalization to flag manifolds will be discussed in Sec.\ 3.1.

 \item [(4)] {\it Bundles}$\,$:
   $V$ is a vector bundle over $X$,
   $V_d=\rho_!\ev^{\ast}V$, $U_d=\rho^{\ast}V_d$,
   and ${\cal U}_d=\pi^{\ast}V_d$.

 \item [(5)] {\it Special $S^1$-fixed-point locus}$\,$:
   $F_0\simeq \overline{M}_{0,1}(X,d)$ is the special $S^1$-fixed-point
    component in $M_d$ that corresponds to gluing stable maps
    $(C^{\prime}, f^{\prime}, x^{\prime})$ to $\CP^1$ at
    $x^{\prime}\in C^{\prime}$ and $\infty\in \CP^1$,
   $Y_0$ is the special $S^1$-fixed-point component in $W_d$
    such that $\varphi^{-1}(Y_0)=F_0$, and
   $E_0$ is the $S^1$-fixed-point locus in $\psi^{-1}(Y_0)$
   and is called the distinguished $S^1$-fixed-point locus
   or components in ${\cal Q}_d$.
  There is a natural smooth morphism $p$ from each component $E_{0s}$
   of $E_0$ onto the flag manifold $X$.

 \item[(6)] {\it Relation of $\psi$ and $\phi$.}
  It will be shown in Sec.\ 3.1 that
   $\varphi(M_d)=\psi({\cal Q}_d)$ and that
  $\psi$ is a resolution of singularities of $\varphi(M_d)$.
  This implies that [L-L-Y1, III: Lemma 5.5] holds.
\end{itemize}

Associated to each $(V,b)$, where $b$ is a multiplicative
characteristic class, is the Euler series
$A(t)\in A^{\ast}(X)(\alpha)[t]\,$:
$$
 \begin{array}{lllll}
  A(t)  & :=  &  A^{V,b}  & :=
    & e^{-H\cdot t/\alpha}\,\sum_d\,A_d\,e^{d\cdot t}\,, \\[1.2ex]
  A_d   & =  & i_0^{\ast}\,b({\cal U}_d) & :=
    & \ev^X_{\ast}\,\left(
            \frac{\rho^{\ast}b(V_d)\cap[M_{0,1}(d,X)]}{
                    e_{{\tinyBbb C}^{\times}}(F_0/M_d)}
                     \right)\;
      =\; \frac{ (i_{X_0}^{\ast}\varphi_{\ast}b({\cal U}_d))\,
                  \cap [X_0] }{ e_{{\tinyBbb C}^{\times}}(X_0/W_d) }\,,
           \hspace{1ex}\mbox{denoted}\hspace{1ex}
            \frac{\Theta_d}{e_{{\tinyBbb C}^{\times}}(X_0/W_d)}\,,
              \\[2ex]
   & & & =
     &  g_{\ast}\left(
           \sum_s\,
             \frac{ (\,
               i_{E_{0s}}^{\ast}\,
                 g^{\ast}\,i_{X_0}^{\ast}\varphi_{\ast}b({\cal U}_d)\,
                    )\,
                \cap\,[E_{0s}] }{
                 e_{{\tinyBbb C}^{\times}}(E_{0s}/{\cal Q}_d) }
                  \right)\,,
         \hspace{1ex}\mbox{denoted}\hspace{1ex}
         g_{\ast}\left(
           \sum_s\,
           \frac{\Xi_{d,s}}{
             e_{{\tinyBbb C}^{\times}}(E_{0s}/{\cal Q}_d)}
                  \right)\,,
 \end{array}
$$
where $\alpha=c_1({\cal O}_{\CPscriptsize^{\infty}})(1)$ is
 the generator for $H^{\ast}_{{\scriptsizeBbb C}^{\times}}(\pt)$.
On the other hand, one has the intersection numbers and their generating
function
$$
 \begin{array}{lllll}
  K_d   & =  & K_d^{V,\,b}   & =
    & \int_{M_{0,0}(d,X)}\,b(V_d)\,,  \\[1.2ex]
  \Phi  & =  & \Phi^{V,\,b}   &  =  &  \sum_d\,K_d\,e^{d\cdot t}\,.
 \end{array}
$$
In the good cases, $K_d$ and $\Phi$ can be obtained from $A_d$ and $A(t)$
 by appropriate integrals of the form $\int_X\,e^{-H\cdot t/\alpha}\,A_d$,
 where $H=(H_1,\,\cdots,\,H_I)$ is the restriction to $X$ of
 the hyperplane classes, also denoted by $H$, on $Y$ from its product
 projective space components, e.g.\ [L-L-Y1, III$\,$: Theorem 3.12].
This integral can be turned into an integral on $E_0\,$:

{\small
$$
 \int_X\,\tau^{\ast}\,e^{H\cdot t}\cap A_d\;
  =\; \int_{Y_0}\,
        e^{H\cdot t}\,
         \cap\,
        g_{\ast}\left(
         \sum_s\,
          \frac{\Xi_{d,s}}{
              e_{{\tinyBbb C}^{\times}}(E_{0s}/{\cal Q}_d)}
                \right)\;
  =\; \sum_s\,
        \int_{E_{0s}}\,
          \frac{g^{\ast}e^{H\cdot t}\,\cap\,\widehat{\Xi}_{d,s}}{
                 e_{{\tinyBbb C}^{\times}}(E_{0s}/{\cal Q}_d)}\,,
 $$
{\normalsize where}} 
$\widehat{\Xi}_{d,s}$ is the Poincar\'{e} dual of $\Xi_{d,s}$
with respect to $[E_{0,s}]$.
As will be discussed in Sec.\ 3.1, $E_{0s}$ is a flag manifold
fibred over $X$ and, hence, $g^{\ast}e^{H\cdot t}$ can be read
off from the natural fibration of flag manifolds
$E_{0s}\rightarrow X$.

Following [L-L-Y1, III$\,$: Sec.\ 5.4], in the case that $b=1$
the above integral is reduced to the integral
$$
 \int_{X}\tau^{\ast}e^{H\cdot t}\cap {\mathbf 1}_d\;
 =\; \sum_s\,\int_{E_{0s}}\,
       \frac{g^{\ast}\psi^{\ast}e^{\kappa\cdot\zeta}}{
              e_{{\tinyBbb C}^{\times}}(E_{0s}/{\cal Q}_d)}\,,
$$
where $\kappa=(\kappa_1,\,\cdots\,,\kappa_I)$ is the tuple of
 hyperplane classes in $W_d$ from its product projective space
 components.
Via the natural smooth morphism $p:E_{0s}\rightarrow X$, one can
 integrate out the fiber of $p$ and lead to an integral over $X$.

In this article, we work out all the equivariant Euler classes
 $e_{{\tinyBbb C}^{\times}}(E_{0s}/{\cal Q}_d)$
 and hence an exact expression of this integral.
This determines $A(t)$ with $b=1$ by [L-L-Y1, II: Lemma 2.5].
In the case of Grassmannian manifolds, the discussion gives also
 the known expression of $A(t)$ in [B-CF-K], via which the Hori-Vafa
 formula for Grassmannian manifolds was checked. For general $b$
 induced from a concavex bundle on $X$, our method gives an
 explicit formula for the hypergeometric series in the mirror
 formulas.

\bigskip

\noindent $\bullet$
{\bf Conventions and notation.}
\begin{itemize}
 \item [(1)]
  For historical reason, due to the relation of the Euler series with
   hyper-geometric series when $X$ is a toric manifold, $A(t)$ will also
   be denoted by $\HG[b]^{X,V}(t)$ and be called a hypergeometric series.

 \item [(2)]
  All the dimensions are {\it complex} dimensions unless otherwise
  noted.

 \item [(3)]
  The $S^1$-actions involved in this article are induced from
  ${\Bbb C}^{\times}$-actions and both have the same fixed-point
  locus. In many places, it is more convenient to phrase things
  in term of ${\Bbb C}^{\times}$-action and we will not distinguish
  the two actions when this ambiguity causes no harm.

 \item [(4)]
  A locally free sheaf and its associated vector bundle are denoted
  the same.

 \item [(5)]
 An $I\times J$  matrix whose $(i,j)$-entry is $a_{ij}$ is denoted
  by $(a_{ij})_{i,j}$ when the position of an entry is
  emphasized and
  by $[a_{ij}]_{I\times J}$ when the size of the matrix is emphasized.

 \item [(6)]
  From Section 2 on, the smooth curve $C$ will be $\CP^1$
  unless other noted.

 \item [(7)]
  All the products of ${\Bbb C}$-schemes are products over
  $\Spec{\Bbb C}$.

 \item [(8)]
  For notation simplicity, the structure sheaf of a scheme is denoted
  also by ${\cal O}$ when the scheme is clear from the contents.
\end{itemize}

\bigskip

\section{The $S^1$-fixed-point components on hyper-Quot-schemes.}

The $S^1$-fixed-point components on the hyper-Quot-scheme
 $\HQuot_P({\cal E}^n)$ and their topology are studied
 in this section.

\bigskip

\subsection{Inclusion pairs of $S^1$-invariant subsheaves of ${\cal E}^n$.}

A characterization of $S^1$-invariant subsheaves of ${\cal E}^n$
 is given in [L-L-L-Y$\,$: Sec.\ 2.1] and is summarized into
 the following fact$\,$:

\bigskip

\noindent
{\bf Fact 2.1.1 [$S^1$-invariant subsheaf].}
(Cf.\ [L-L-L-Y$\,$: Sec.\ 2.1].)
{\it
 \begin{itemize}
  \item[{\rm (1)}]
   An $S^1$-invariant subsheaf ${\cal V}$ of ${\cal E}^n$ with Hilbert
    polynomial of ${\cal E}^n\!/\mbox{\raisebox{-.4ex}{${\cal V}$}}$
    being $P(t)=(n-r)\,t+d+(n-r)$ is characterized by the following
    data
    $(\,V^{(0)}_{\bullet}, \alpha_{\bullet}\,;\,
        V^{(\infty)}_{\bullet},\beta_{\bullet}\,):$
    \begin{itemize}
     \item [{\rm ({\it i})}]
      A pair of flags $(V_{\bullet}^{(0)},V_{\bullet}^{(\infty)})$
      of ${\Bbb C}^n\,:$
      $$
       V_{\bullet}^{(0)}\,:\,
         V^{(0)}_1\hookrightarrow\,\cdots\,\hookrightarrow\, V_k^{(0)}
           \hookrightarrow {\Bbb C}^n
       \hspace{2ex}\mbox{and}\hspace{2ex}
       V_{\bullet}^{(\infty)}\,:\,
         V^{(\infty)}_1\hookrightarrow\,\cdots\,\hookrightarrow\,
         V_l^{(\infty)}\hookrightarrow {\Bbb C}^n
      $$
      with $V_k^{(0)}=V_l^{(\infty)}$, both of dimension $r$.

     \item [{\rm ({\it ii})}]
      A pair of integer sequences
       {\rm (cf.\ [L-L-L-Y: Definition 2.1.6])}
       $$
        \alpha_{\bullet}: 0\le \alpha_1\le \cdots \le \alpha_r
         \hspace{2em}\mbox{and}\hspace{2em}
        \beta_{\bullet}: 0\le \beta_1\le \cdots \le \beta_r
       $$
       that satisfy
       $\;(\alpha_1+\cdots+\alpha_r)+(\beta_1+\cdots+\beta_r)=d$.
    \end{itemize}

  \item[{\rm (2)}]
   For any $S^1$-invariant coordinate system on $\,\CP^1${\rm :}
    $$
     \CP^1=U_0\cup U_{\infty}\,,
      \hspace{2ex}\mbox{where}\hspace{2ex}
      U_0=\Spec{\Bbb C}[z] \hspace{2ex}\mbox{and}\hspace{2ex}
      U_{\infty}=\Spec{\Bbb C}[w]
    $$
    with
     the gluing
      $$
       \Spec{\Bbb C}[z]\,
        \hookleftarrow\, \Spec{\Bbb C}[z,z^{-1}]\,
        \stackrel{\hspace{2ex}z\leftrightarrow w^{-1}}{\simeq}\,
         \Spec{\Bbb C}[w,w^{-1}]\,
        \hookrightarrow\, \Spec{\Bbb C}[w]\,.
      $$
     with the $S^1$-action$\,:$
      $z\mapsto e^{i\theta }z$ and $w\mapsto e^{-i\theta}z$,
    a data
     $(V^{(0)}_{\bullet},\alpha_{\bullet}\,;\,
       V^{(\infty)}_{\bullet},\beta_{\bullet})$ in Item {\rm (1)}
     determines local subsheaves
     $({\cal V}^{(0)},{\cal V}^{(\infty)})$ of
     $({\cal E}^n|_{U_0},{\cal E}^n|_{U_{\infty}})$,
     which automatically glue together over $U_0\cap U_{\infty}$
     via the canonical isomorphism
     $$
      {\cal V}^{(0)}|_{U_0\cap U_{\infty}}\;
       \simeq\; {\cal E}^n|_{U_0\cap U_{\infty}}\;
       \simeq\; {\cal E}^n|_{U_0\cap U_{\infty}}
     $$
     as ${\cal O}_{U_0\cap U_{\infty}}$-modules, and hence
     an $S^1$-invariant subsheaf ${\cal V}$ of ${\cal E}^n$
     with the required Hilbert polynomial for
     ${\cal E}^n\!/\mbox{\raisebox{-.4ex}{${\cal V}$}}$.

  \item[{\rm (3)}]
   Given the data
    $(V^{(0)}_{\bullet}, \alpha_{\bullet}\,;\,
                         V^{(\infty)}_{\bullet},\beta_{\bullet})$
    in Item {\rm (1)}, write
    $(\alpha_{\bullet}\,;\,\beta_{\bullet})$
    as $(a_{\bullet}, m_{\bullet}\,;\, b_{\bullet}, n_{\bullet})$:
    {\small
    $$
     \begin{array}{ccccccccccccccc}
      0 & \le  & a_1\,(=\alpha_1) & <  & \cdots & < & a_k\,(=\alpha_r)\,;
        & &   0 & \le  & b_1\,(=\beta_1) & <
                                       & \cdots & < & b_l\,(=\beta_r)
                                                        \\[.6ex]
        &      & m_1 &    & \cdots &   & m_k
        & &     &      & n_1 &    & \cdots &   & n_l
     \end{array}
    $$
    {\normalsize\it with}} 
    the multiplicity of $a_i$, $b_j$ indicated.
   Recall the notation $M^{\sim}$ for the coherent sheaf on $\Spec A$
    associated to an $A$-module $M$, cf. {\rm [Ha]}.
   Then, in Item {\rm (2)},
    $$
     {\cal V}^{(0)}\;
      =\; \left(\,
             z^{a_1}\,W_1^{(0)}\, +\, \cdots\,
                                +\, z^{a_k}\, W_k^{(0)}\,
          \right)^{\sim}\,,
    $$
    where $W_i^{(0)}$ is any subspace in $V_i^{(0)}-V_{i-1}^{(0)}$
    of rank $m_i$, and
    $$
     {\cal V}^{(\infty)}\;
      =\; \left(\,w^{b_1}\,W_1^{(\infty)}\, +\, \cdots\,
                                +\, w^{b_l}\, W_l^{(\infty)}\,
           \right)^{\sim}\,,
    $$
    where $W_j^{(\infty)}$ is any subspace in
    $V_j^{(\infty)}-V_{j-1}^{(\infty)}$ of rank $n_j$.
 \end{itemize}
} 
(Cf.\ {\sc Figure} 2-1-1; see also {\sc Figure} 2-1-2.)

\begin{figure}[htbp]
 \setcaption{{\sc Figure} 2-1-1.
  \baselineskip 14pt
   An $S^1$-invariant subsheaf ${\cal V}$ of a trivialized trivial
   bundle ${\cal E}^n$ is characterized by a pair of flags with
   identical last element, together with integral labels on elements
   in the flags.
 } 
 \centerline{\psfig{figure=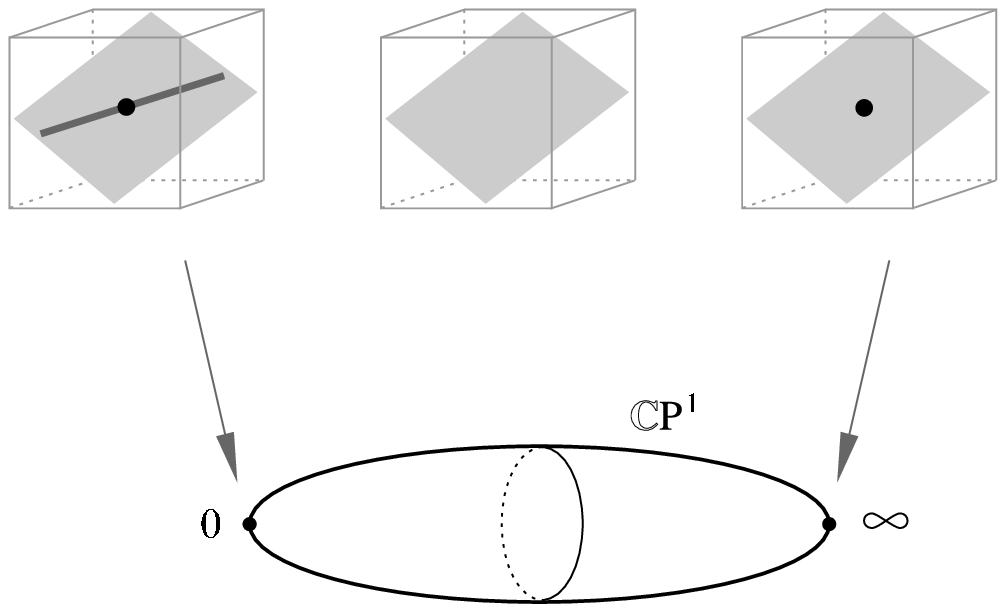,width=13cm,caption=}}
\end{figure}

\bigskip

The goal of this subsection is to generalize the above result to the case
 of inclusion sequences of $S^1$-invariant subsheaves of ${\cal E}^n$.

\bigskip

\noindent
{\bf Definition 2.1.2 [(integrally) labelled flag].} {\rm
 An ({\it integrally}) {\it labelled flag} of ${\Bbb C}^n$
 $$
  V_{\bullet}(s_{\bullet})\;:\; V_1(s_1)\; \hookrightarrow\;
    \cdots\; \hookrightarrow\; V_k(s_k)
 $$
 is an ordinary flag
  $V_1\hookrightarrow\cdots\hookrightarrow V_k\hookrightarrow{\Bbb C}^n$
  together with a label $s_i\in{\Bbb Z}$ attached to each $V_i$ such that
  $s_1 < \cdots < s_k$.
} 

\bigskip

With the same notation, Fact 2.1.1 Item (1) can be rephrased as follows.

\bigskip

\noindent
{\bf Fact 2.1.1$^{\prime}$ [$S^1$-invariant subsheaf].} {\it
 An $S^1$-invariant subsheaf ${\cal V}$ of
  ${\cal E}^n={\cal O}_{\CPscriptsize^1} \otimes {\Bbb C}^n$
  is characterized by a pair of labelled flags
  $(V_{\bullet}^{(0)}(s_{\bullet})\,,\,
    V_{\bullet}^{(\infty)}(t_{\bullet}))$,
 where the label $s_i$ of $V^{(0)}_i$ is $\alpha_{dim\,V^{(0)}_i}$
  and the label $t_j$ of $V^{(\infty)}_j$ is
  $\beta_{dim\,V^{(\infty)}_j}$.
} 

\bigskip

\noindent
{\bf Definition 2.1.3 [admissible inclusion].} {\rm
 Given a labelled flag
  $$
   V_{\bullet}\;:\; V_1(s_1)\; \hookrightarrow\; \cdots\;
    \hookrightarrow\; V_k(s_k) \;
    \hookrightarrow\; V_{k+1}(\infty):={\Bbb C}^n(\infty)
  $$
  and a labelled vector subspace $\Pi(s)\subset {\Bbb C}^n$,
  we say that $\Pi(s)$ is {\it admissibly contained} in $V_i(s_i)$
  for some $i$ if $\Pi\subset V_i$ and $s\ge s_i$.
 Since the sequence of the labels of the flag is non-decreasing,
  there is a maximal $i\le k$ such that $\Pi(s)$ is admissibly
  contained in $V_i(s_i)$ but not in $V_{i+1}(s_{i+1})$.
  In this case, we say that $\Pi(s)$ is {\it admissibly and
  critically} contained in $V_i(s_i)$.
} 

\bigskip

The following lemma is the inductive step in understanding an
 inclusion sequence of $S^1$-invariant subsheaves of ${\cal E}^n$.

\bigskip

\noindent
{\bf Lemma 2.1.4 [$S^1$-invariant pair].} {\it
  Let ${\cal V}_1$ and ${\cal V}_2$ be $S^1$-invariant subsheaves of
  ${\cal E}^n$ characterized by the data
  $(V_{1,\bullet}^{(0)}(s_{1,\bullet})\,;\,
    V_{1,\bullet}^{(\infty)}(t_{1,\bullet}))$
  and
  $(V_{2,\bullet}^{(0)}(s_{2,\bullet})\,;\,
    V_{2,\bullet}^{(\infty)}(t_{2,\bullet}))$
  respectively.
 Then ${\cal V}_1$ is a subsheaf of ${\cal V}_2$
  if and only if each labelled subspace in the sequence
  $V_{1,\bullet}^{(0)}(s_{1,\bullet})$
  {\rm (}resp.\
   $V_{1,\bullet}^{(\infty)}(t_{1,\bullet})${\rm )}
  can be admissibly contained in a labelled subspace
  in the sequence $V_{2,\bullet}^{(0)}(s_{2,\bullet})$
  {\rm (}resp.\
   $V_{2,\bullet}^{(\infty)}(t_{2,\bullet})${\rm )}.
} 

\bigskip

\noindent
{\bf Definition 2.1.5 [order/precedence].} {\rm
 When the condition in Lemma 2.1.4 is met, we shall say that
  $(V_{1,\bullet}^{(0)}(s_{1,\bullet})\,;\,
    V_{1,\bullet}^{(\infty)}(t_{1,\bullet}))$
  {\it precedes}
   $(V_{2,\bullet}^{(0)}(s_{2,\bullet})\,;\,
     V_{2,\bullet}^{(\infty)}(t_{2,\bullet}))$.
 In notation,
  $(V_{1,\bullet}^{(0)}(s_{1,\bullet})\,;\,
    V_{1,\bullet}^{(\infty)}(t_{1,\bullet}))\,
    \preccurlyeq\,
      (V_{2,\bullet}^{(0)}(s_{2,\bullet})\,;\,
       V_{2,\bullet}^{(\infty)}(t_{2,\bullet}))$.
} 

\bigskip

\noindent
{\it Proof of Lemma 2.1.4.}
 Given a pair of flags in ${\Bbb C}^n$
  $$
   (V_{\bullet},V^{\prime}_{\bullet^{\prime}})\;
    :=\; (\,V_1\hookrightarrow \cdots \hookrightarrow V_k
               \hookrightarrow {\Bbb C}^n\,,\,
            V^{\prime}_1\hookrightarrow \cdots
               \hookrightarrow V_{k^{\prime}}^{\prime}
                \hookrightarrow {\Bbb C}^n\,)\,,
  $$
  by considering either the double filtration or the double graded
  object of ${\Bbb C}^n$ associated to the pair of flags, one can show
  that there exists a direct-sum decomposition
  ${\Bbb C}^n=\oplus_m E_m$ of ${\Bbb C}^n$ such that any $V_i$,
  $V_{i^{\prime}}^{\prime}$ is a sum of some direct summands in this
  decomposition$\,:$
  $$
   V_i\; =\; \oplus_j\,E_{i_j} \hspace{2em}\mbox{and}\hspace{2em}
   V_{i^{\prime}}^{\prime}\;
     =\; \oplus_{j^{\prime}}\,E_{i^{\prime}_{j^{\prime}}}\,.
  $$
 Such a decomposition of ${\Bbb C}^n$ is said to be {\it compatible}
  with the pair of flags $(V_{\bullet},V_{\bullet^{\prime}})$.

 Apply this to our problem first with
  $$
   (V_{\bullet}(s_{\bullet}),
    V^{\prime}_{\bullet}(s_{\bullet}^{\prime}))\;
   =\; \left(\, V^{(0)}_{1,\bullet}(s_{1,\bullet})\;,\;
                 V^{(0)}_{2,\bullet}(s_{2,\bullet})\,\right)
  $$
  and choose $W_i \subset V_i-V_{i-1}$ and
  $W^{\prime}_{i^{\prime}}
   \subset V^{\prime}_{i^{\prime}}-V^{\prime}_{i^{\prime}-1}$,
  as defined in Fact 2.1.1 (3), to be also direct sums
  with the summands some $E_m$'s:
  $$
   W_i\;=\; \oplus_j\, E_{i_j}
    \hspace{1em}\mbox{and}\hspace{1em}
   W^{\prime}_{i^{\prime}}\;
        =\; \oplus_{j^{\prime}}\, E_{{i^{\prime}}_{j^{\prime}}}\,.
  $$
 Recall Fact 2.1.1 (3), Fact 2.1.1$^{\prime}$ and the notations therein.
 Then
  $$
   {\cal V}^{(0)}\; :=\; {\cal V}_1|_{U_0}\;
   = \; \left(\,\oplus_i\, z^{a_i} W_i\,\right)^{\sim}\;
   = \; \left(\,\oplus_i\, \oplus_j\,
                   z^{a_i} E_{i_j}\,\right)^{\sim}
  $$
  while
  $$
   {\cal V}^{\prime\,(0)}\; :=\; {\cal V}_2|_{U_0}\;
   = \; \left(\,\oplus_{i^{\prime}}\,
                 z^{a^{\prime}_{i^{\prime}}}
                             W^{\prime}_{i^{\prime}}\,\right)^{\sim}\;
   = \; \left(\,\oplus_{i^{\prime}}\,\oplus_{j^{\prime}}\,
                 z^{a^{\prime}_{i^{\prime}}} E_{i^{\prime}_{j^{\prime}}}\,
        \right)^{\sim}\,.
  $$
 Each $E_{i_j}$ in $W_i$ appears exactly once as some
  $E_{i^{\prime}_{j^{\prime}}}$ in some $W_{i^{\prime}}$.
 Consequently, ${\cal V}^{(0)}$ is a subsheaf of ${\cal V}^{\prime\:(0)}$
  if and only if
  $$
   a_i\; \ge\; a^{\prime}_{i^{\prime}}
    \hspace{1em}\mbox{whenever
     $\;E_{i_j}=E_{i^{\prime}_{j^{\prime}}}$}\,.
  $$
 But this means precisely that
  $V^{(0)}_{1,\bullet}(s_{1,\bullet})
    \preccurlyeq\, V^{(0)}_{2,\bullet}(s_{2,\bullet})$.
 (Cf.\ {\sc Figure} 2-1-2.)

 \begin{figure}[htbp]
  \setcaption{{\sc Figure} 2-1-2.
   \baselineskip 14pt
    The relation of the characterization data of successive
     $S^1$-invariant subsheaves ${\cal V}_1\hookrightarrow{\cal V}_2$
     of ${\cal E}^n$.
  } 
  \centerline{\psfig{figure=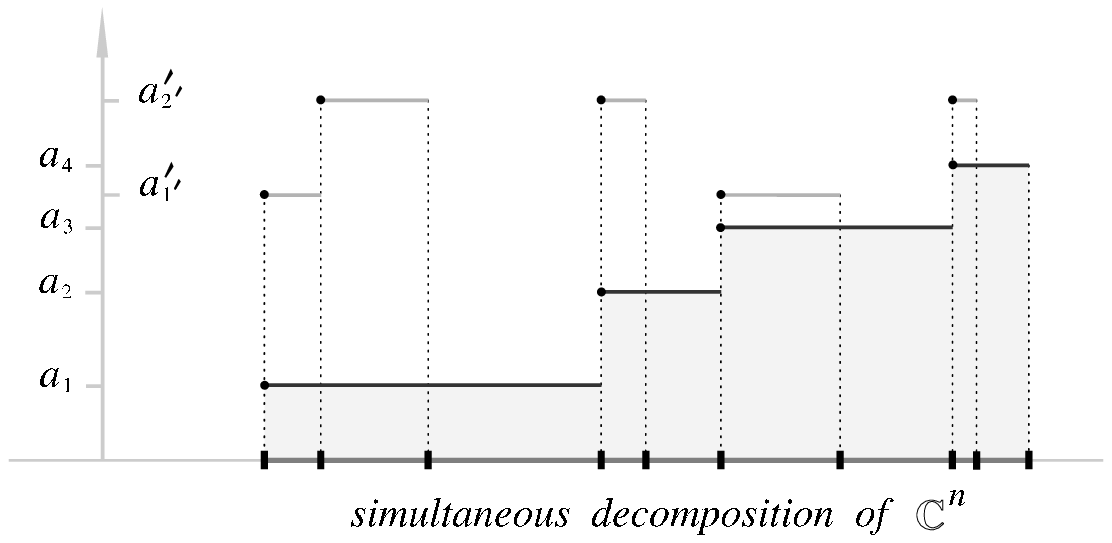,width=13cm,caption=}}
 \end{figure}

 Apply the same argument next to
  $$
   (V_{\bullet}(s_{\bullet}),
    V_{\bullet}^{\prime}(s_{\bullet}^{\prime}))\;
   =\; \left(\,
         V^{(\infty)}_{1,\bullet}(t_{1,\bullet})\;,\;
         V^{(\infty)}_{2,\bullet}(t_{2,\bullet})\,\right)
  $$
  to conclude that
  $V^{(\infty)}_{\bullet}(t_{1,\bullet})
      \preccurlyeq\, V^{(\infty)}_{\bullet}(t_{2,\bullet})$
  as well.
 This completes the proof.

\noindent\hspace{14cm}$\Box$

\bigskip

To better describe the structure of $S^1$-fixed-point components
 in $\HQuot_P({\cal E}^n)$, we introduce a couple of definitions
 in the passing.

\bigskip

\noindent
{\bf Definition/Lemma 2.1.6
     [admissible pair of $(\alpha_{\bullet};\beta_{\bullet})$].}
{\it
  Let
  $(\alpha_{2,\bullet};\beta_{2,\bullet})$ be from the
   characterization data of the $S^1$-invariant subsheaf ${\cal V}_2$
   of ${\cal E}^n$,
  $P_1=P_1(t)=(n-r_1)t+ d_1 + (n-r_1)$ be an integral polynomial and
  $(\alpha_{1,\bullet};\beta_{1,\bullet})$ be admissible to $P_1$.
 Then
  there exists an $S^1$-invariant subsheaf ${\cal V}_1$ of
  ${\cal V}_2$ whose labels in its characterization data comes from
  $(\alpha_{1,\bullet};\beta_{1,\bullet})$
 if and only if
  $$
   r_1\; \le\; r_2\,, \hspace{1em}
    \alpha_{1,i}\;\ge\;\alpha_{2,i}
      \hspace{1em}\mbox{and}\hspace{1em}
    \beta_{1,i}\; \ge\; \beta_{2,i}
     \hspace{1em}\mbox{for $i=1,\,\ldots,\,r_1$}\,.
  $$
 This condition depends only on the data
  $(\alpha_{2,\bullet};\beta_{2,\bullet})$,
  not on any other detail of ${\cal V}_2$.
 We say that $(\alpha_{1,\bullet};\beta_{1,\bullet})$
  is admissible to $(\alpha_{2,\bullet};\beta_{2,\bullet})$.
 In notation,
  $(\alpha_{1,\bullet}\,;\,\beta_{1,\bullet})
    \rightarrow (\alpha_{2,\bullet}\,;\,\beta_{2,\bullet})$.
} 

\bigskip

\noindent
{\bf Definition/Lemma 2.1.7 [characteristic chain of subspaces].}
{\it
 Following Definition/Lemma 2.1.6,
 given
  $(\alpha_{1,\bullet}\,;\,\beta_{1,\bullet})
    \rightarrow (\alpha_{2,\bullet}\,;\,\beta_{2,\bullet})$,
 let ${\cal V}_1$ be an $S^1$-invariant subsheaf of ${\cal V}_2$
  with characterization data
  $(V^{(0)}_{1,\bullet}(s_{1,\bullet})\,;\,
    V^{(\infty)}_{1,\bullet}(t_{1,\bullet}))$
  equivalent to
  $(V^{(0)}_{1,\bullet}, \alpha_{1,\bullet}\,;\,
    V^{(\infty)}_{1,\bullet}, \beta_{1,\bullet})$.
 Recall that each element in $V^{(0)}_{1,\bullet}$
  {\rm (}resp.\ $V^{(\infty)}_{1,\bullet}${\rm )}
  is admissibly and critically contained in a unique element
  in $V^{(0)}_{2,\bullet}(s_{2,\bullet})$
  {\rm (}resp.\
   $V^{(\infty)}_{2,\bullet}(t_{2,\bullet})${\rm )}.
 These latter elements form a sequence of successive subspaces
  $$
   \Pi^{(0)}_{\bullet} \hspace{2em}
     (\,\mbox{resp.}\hspace{2ex}
     \Pi^{(\infty)}_{\bullet}\,)\,.
  $$
 Then the pair $(\Pi^{(0)}_{\bullet}\,;\,\Pi^{(\infty)}_{\bullet})$
  depends only on $(\alpha_{1, \bullet}\,;\,\beta_{1,\bullet})$,
  not on the choice of ${\cal V}_1$.
} 

\bigskip

\noindent
{\bf Definition 2.1.8 [restrictive flag manifold].} {\rm
 Given
  an inclusion sequence of subspaces
   $\Pi_{\bullet}:\Pi_1\subset \cdots \subset \Pi_s \subset {\Bbb C}^n$
   (with $\Pi_i=\Pi_{i+1}$ allowed) and
  a strictly increasing sequence of integers
   $0<k_1<\cdots <k_s$ with $k_i\le \dimm \Pi_i =: l_i$,
 the subspace $\Fl_{k_1,\cdots,k_s}({\Bbb C}^n,\Pi_{\bullet})$
  of the flag manifold $\Fl_{k_1,\cdots,k_s}({\Bbb C}^n)$ defined by
 $$
  \Fl_{k_1,\cdots,k_s}({\Bbb C}^n,\Pi_{\bullet})\;
  :=\; \left\{\,
         \mbox{flags
          $V_{\bullet}:V_1\hookrightarrow\cdots\hookrightarrow V_s
                            \hookrightarrow {\Bbb C}^n$}\;
        \left|\;
         \mbox{$\dimm V_i=k_i\,$ and $\,V_i\subset \Pi_i$}\,
                \right. \right\}
 $$
  is called
  a {\it restrictive flag manifold associated to $\Pi_{\bullet}$}.
} 

\bigskip

\noindent
{\bf Lemma 2.1.9
 [$\Fl_{k_1,\cdots,k_s}({\Bbb C}^n,\Pi_{\bullet})$ smooth].}
{\it
 $\Fl_{k_1,\cdots,k_s}({\Bbb C}^n,\Pi_{\bullet})$
 is a projective, connected, smooth manifold of dimension
 $k_1(l_1-k_1)+(k_2-k_1)(l_2-k_2)+\, \cdots\,+(k_s-k_{s-1})(l_s-k_s)$.
} 

\bigskip

\noindent
{\it Proof.}
 Projectivity of $\Fl_{k_1,\cdots,k_s}({\Bbb C}^n,\Pi_{\bullet})$
  follows from the prejectivity of $\Fl_{k_1,\cdots,k_s}({\Bbb C}^n)$.
 By construction, $\Fl_{k_1,\cdots,k_s}({\Bbb C}^n,\Pi_{\bullet})$
  admits a tower of fibrations by Grassmannian manifolds and hence
  is connected.
 The fibrations in the tower is not topologically locally trivial,
  in particular the fibers of a fixed fibration in the tower can vary;
  so it is not immediate from this fibration tower that
  $\Fl_{k_1,\cdots,k_s}({\Bbb C}^n,\Pi_{\bullet})$ is smooth.

 To prove the last statement, let $F_{\bullet}{\Bbb C}^n$ be the
  filtration of ${\Bbb C}^n$ by $\Pi_{\bullet}$, then one can show
  that at each point $[V_{\bullet}]$ of
  $\Fl_{k_1,\cdots,k_s}({\Bbb C}^n,\Pi_{\bullet})$, the space
  $T_{[V_{\bullet}]}\Fl_{k_1,\cdots,k_s}({\Bbb C}^n,\Pi_{\bullet})$
  of the first order deformations of $[V_{\bullet}]$ as a restrictive
  flag fits into an exact sequence of complex vector spaces of the form
  $$
   0\;\longrightarrow\;
    \Hom(V_{\bullet}, V_{\bullet})\;
     \longrightarrow\;
    \Hom(V_{\bullet}, F_{\bullet}{\Bbb C}^n)\;
     \longrightarrow\;
    T_{[V_{\bullet}]}\Fl_{k_1,\,\ldots,\,k_s}({\Bbb C}^n, \Pi_{\bullet})\;
     \longrightarrow\; 0\,.
  $$
 (Cf.\ See Sec.\ 3.5 for more related details in the construction of
       an Euler sequence for
       $T_{\ast}\Fl_{k_1,\cdots,k_s}({\Bbb C}^n,\Pi_{\bullet})$.)
 It follows that
  $$
   \dimm T_{[V_{\bullet}]}
           \Fl_{k_1,\,\ldots,\,k_s}({\Bbb C}^n,\Pi_{\bullet})\;
   =\; k_1(l_1-k_1)\,+\, (k_2-k_1)(l_2-k_2)\,+\, \cdots\,
      +\, (k_s-k_{s-1})(l_s-k_s)\,.
  $$
 Since this is independent of $[V]$ and all elements in
  $T_{[V_{\bullet}]}
    \Fl_{k_1,\,\ldots,\,k_s}({\Bbb C}^n,\Pi_{\bullet})$
  are realizable from a family of restrictive flags over a small disc,
  $\Fl_{k_1,\cdots,k_s}({\Bbb C}^n,\Pi_{\bullet})$, as a scheme, must
  be reduced everywhere and hence is smooth of the above dimension.
 This concludes the lemma.

\noindent\hspace{14cm}$\Box$

\bigskip

With these preparations, we can now describe first the topology of
 the connected components of the $S^1$-fixed-point locus in the special
 Quot-scheme
 $\Quot_{P_1}({\cal V}_2\hookrightarrow {\cal E}^n)$ and
then the topology of the connected components of the
 $S^1$-fixed-point components in the general hyper-Quot-scheme
 $\HQuot_P({\cal E}^n)$.

\bigskip

\noindent
{\bf Lemma 2.1.10
     [$S^1$-invariant subsheaves of a fixed $S^1$-invariant sheaf].}
{\it
 Fix an $S^1$-invariant subsheaf ${\cal V}_2$ of ${\cal E}^n$
  with characterization data
  $(\,V^{(0)}_{2,\bullet}(s_{2,\bullet})\,;\,
      V^{(\infty)}_{2,\bullet}(t_{2,\bullet})\,)$
  equivalent to
  $(\,V^{(0)}_{2,\bullet}, \alpha_{2,\bullet}\,;\,
      V^{(\infty)}_{2,\bullet},\beta_{2,\bullet}\,):$
  and the Hilbert polynomial of
   ${\cal E}^n\!/\mbox{\raisebox{-.4ex}{${\cal V}_2$}}$
   being $P_2=P_2(t)=(n-r_2)\,t + d_2 + (n-r_2)$.
 \begin{itemize}
  \item[{\rm (1)}]
   The space
    $F^{S^1}_{P_1}({\cal V}_2\hookrightarrow{\cal E}^n)$
    of $S^1$-invariant subsheaves ${\cal V}_1$ of ${\cal V}_2$
    with the Hilbert polynomial of
    ${\cal E}^n\!/\mbox{\raisebox{-.4ex}{${\cal V}_1$}}$ being
    $P_1=P_1(t)=(n-r_1)\,t+d_1+(n-r_1)$ is non-empty
    if and only if there exists
    $(\alpha_{1,\bullet}\,;\,\beta_{1,\bullet})$ admissible both
    to $P_1$ and $(\alpha_{2,\bullet}\,;\,\beta_{2,\bullet})$.

  \item[{\rm (2)}]
   The set of connected components of
    $F^{S^1}_{P_1}({\cal V}_2\hookrightarrow{\cal E}^n)$
    is in one-to-one correspondence with the set of pairs
    $(\alpha_{1,\bullet}\,;\,\beta_{1,\bullet})$
    that is admissible both to $P_1$ and
    $(\alpha_{2,\bullet}\,;\,\beta_{2,\bullet})$.
   Let $(A;B)$ be a pair of incomplete matrices defined by
   {\scriptsize
    $$
     A\;=\;\left[\,
             \begin{array}{lllll}
               \alpha_{1,1} & \cdots & \alpha_{1,r_1}\\[.6ex]
               \alpha_{2,1} & \cdots & \alpha_{2,r_1}
                                           & \cdots & \alpha_{2,r_2}
             \end{array}
           \,\right]_{2\times r_2}
        \hspace{1em}\mbox{\normalsize and}\hspace{1em}
     B\;=\;\left[\,
             \begin{array}{lllll}
               \beta_{1,1} & \cdots & \beta_{1,r_1}\\[.6ex]
               \beta_{2,1} & \cdots & \beta_{2,r_1}
                                           & \cdots & \beta_{2,r_2}
             \end{array}
           \,\right]_{2\times r_2}\,.
    $$
   {\normalsize\it $($Note}} 
    that $\alpha_{1,j}$ and $\beta_{1,j}$ for $j>r_1$
     is left undefined/blank;
    this is why we call $A$ and $B$ incomplete matrices.$)$
   The corresponding component of
    $F^{S^1}_{P_1}({\cal V}_2\hookrightarrow {\cal E}^n)$
    will be denoted by $F_{(A;B)}$.

  \item[{\rm (3)}]
   Each $(A;B)$ in Item {\rm (2)}
    determines a pair
    $(\Pi^{(0)}_{2,\bullet}\,;\,\Pi^{(\infty)}_{2,\bullet})$
    of chains of subspaces by Definitio/Lemma 2.1.7
   and hence a pair of restrictive flag manifolds
    $$
     \Fl_{(\alpha_{1,\bullet})}({\Bbb C}^n,\Pi^{(0)}_{\bullet})\;
      :=\; \Fl_{m_{1,1}, m_{1,1}+m_{1,2},\,\ldots,\, r_1}
               ({\Bbb C}^n,\Pi^{(0)}_{2,\bullet})
    $$
    and
    $$
     \Fl_{(\beta_{1,\bullet})}({\Bbb C}^n,\Pi^{(\infty)}_{\bullet})\;
      :=\; \Fl_{n_{1,1}, n_{1,1}+n_{1,2},\,\ldots,\, r_1}
                ({\Bbb C}^n,\Pi^{(\infty)}_{2,\bullet})\,.
    $$
   Notice that both
    $\Fl_{(\alpha_{1,\bullet})}({\Bbb C}^n,\Pi^{(0)}_{2,\bullet})$ and
    $\Fl_{(\beta_{1,\bullet})}({\Bbb C}^n,\Pi^{(\infty)}_{2,\bullet})$
    fiber over $\Gr_{r_1}({\Bbb C}^n)$.
   The topology of $F_{(A;B)}$ is then given by
    $$
     F_{(A;B)}\;
     =\;
      \Fl_{(\alpha_{1,\bullet})}({\Bbb C}^n, \Pi^{(0)}_{2,\bullet})\,
       \times_{Gr(n,r_1)}\,
      \Fl_{(\beta_{1,\bullet})}({\Bbb C}^n, \Pi^{(\infty)}_{2,\bullet})\,.
    $$
   It depends only on the connected $S^1$-fixed-point component
    $V^{\prime}$ belongs to in the $\Quot$-scheme
    $\Quot_{P^{\prime}}({\cal E}^n)$.
   In other words, it depends only on the admissible incomplete
   matrix $(A;B)$. {\rm(}This justifies our notation.{\rm)}
 \end{itemize}
} 

\bigskip

\noindent
{\it Proof.}
 These follow immediately from our discussion in this subsection.

\noindent\hspace{14cm}$\Box$

\bigskip

\noindent
{\bf Corollary 2.1.11
     [$S^1$-fixed-point locus in $\HQuot_{P_1,P_2}({\cal E}^n)$].} {\it
 The set of connected components of $S^1$-fixed-point locus in
  $\HQuot_{P_1,P_2}({\cal E}^n)$ is in a natural one-to-one
  correspondence with the set of pairs
  {\small
  $$
   \left\{\,
      \left(
       \begin{array}{c}
         (\alpha_{1,\bullet}\,;\,\beta_{1,\bullet}) \\[.6ex]
         (\alpha_{2,\bullet}\,;\,\beta_{2,\bullet})
       \end{array}
      \right) \,\left|\,
        \begin{array}{l}
           \mbox{$(\alpha_{1,\bullet}\,;\,\beta_{1,\bullet})$
                  admissible to $P_1$,
                 $(\alpha_{2,\bullet}\,;\,\beta_{2,\bullet})$
                  admissible to $P_2$}  \\[.6ex]
           \mbox{and
             $(\alpha_{1,\bullet}\,;\,\beta_{1,\bullet})\;
                \rightarrow\;
                 (\alpha_{2,\bullet}\,;\,\beta_{2,\bullet})$}
        \end{array}
              \right.\,\right\}\,,
  $$
  {\normalsize i.e.}}\ 
  the set of pairs of incomplete matrices $(A;B)$ given in
  Lemma 2.1.10, Item $(2)$.
 Denote the $S^1$-fixed-point component corresponding to $(A;B)$
  by $E_{(A;B)}$. Then $E_{(A;B)}$ is smooth.
} 

\bigskip

\noindent
{\it Proof.}
 Recall from [L-L-L-Y: Sec.\ 2.1] that the pair
  $(\alpha_{1,\bullet};\beta_{1,\bullet})$
  (resp.\ $(\alpha_{2,\bullet};\beta_{2,\bullet})$) admissible to
  $P_1$ (resp.\ $P_2$) determines a unique $S^1$-fixed-point component
  $E_{(\alpha_{1,\bullet};\beta_{1,\bullet})}$
  (resp.\ $E_{(\alpha_{2,\bullet};\beta_{2,\bullet})}$) in
  $\Quot_{P_1}({\cal E}^n)$ (resp.\ $\Quot_{P_2}({\cal E}^n)$),
  whose topology is a fibered product of flag manifolds and
  hence is connected and smooth.

 Consider the natural inclusion
 $\iota: \HQuot_{P_1,P_2}({\cal E}^n) \hookrightarrow
         \Quot_{P_1}({\cal E}^n)\times\Quot_{P_2}({\cal E}^n)$.
 Let $E_{(A;B)}$ be the $S^1$-fixed-point sublocus in
  $\HQuot_{P_1,P_2}({\cal E}^n)$ that consists of points associated
  to pairs of $S^1$-invariant subsheaves of ${\cal E}^n$ whose
  discrete part of the characterization data is given by $(A;B)$.
 Then
  $$
   E_{(A;B)}\;
    =\; (\,E_{(\alpha_{1,\bullet};\beta_{1,\bullet})}\,
           \times\, E_{(\alpha_{2,\bullet};\beta_{2,\bullet})}\,)\,
     \cap\, \iota\,\left(\,HQuot_{P_1,P_2}^{S^1}({\cal E}^n)\,\right)\,.
  $$
 With respect to the ambient product structure via $\iota$,
  $E_{(A;B)}$ fibers over $E_{(\alpha_{2,\bullet};\beta_{2,\bullet})}$
  with fiber $F_{(A;B)}$ in Lemma 2.1.10, Item (2).
 Since both $E_{(\alpha_{2,\bullet};\beta_{2,\bullet})}$ and $F_{(A;B)}$
  are connected and smooth, so is $E_{(A;B)}$.
 This completes the proof.

\noindent\hspace{14cm}$\Box$

\bigskip

\bigskip

\subsection{The $S^1$-fixed-point locus in $\HQuot_P({\cal E}^n)$.}

Successive applications of
 Lemma 2.1.4, Lemma 2.1.10 and Corollary 2.1.11
 give rise to the following description of the topology of
 the $S^1$-fixed-point components in $\HQuot_P({\cal E}^n)$.

\bigskip

\noindent
{\bf Proposition 2.2.1
     [$S^1$-fixed-point component in $\HQuot_P({\cal E}^n)$].} {\it
 Recall the sequence of Hilbert polynomials
  $$
   P\;:\; P_1\,,\, \cdots\,,\, P_I\,.
  $$
 \begin{itemize}
  \item[{\rm (1)}]
   For each admissible sequence of pairs of finite sequences:
   $$
    (\alpha_{1,\bullet}\,;\,\beta_{1,\bullet})\;
      \rightarrow\; \cdots\;
      \rightarrow\; (\alpha_{I, \bullet}\,;\,\beta_{I,\bullet})
    \hspace{2em}
    \mbox{with $\;(\alpha_{i,\bullet}\,;\,\beta_{i,\bullet})\;$
          admissible to $\;P_i$}\,,
   $$
   define the $I\times r_I$ incomplete matrices
   {\footnotesize
   $$
    A\;=\;\left[\,
             \begin{array}{lllll}
               \alpha_{1,1} & \cdots & \alpha_{1,r_1} \\
                \cdots      & \cdots & \cdots          & \cdots \\
               \alpha_{I,1} & \cdots & \cdots
                                     & \cdots & \alpha_{I,r_I}
             \end{array}
           \,\right]_{I\times r_I}
        \hspace{1em}\mbox{\normalsize\it and}\hspace{1em}
     B\;=\;\left[\,
             \begin{array}{lllll}
               \beta_{1,1} & \cdots & \beta_{1,r_1} \\
                \cdots     & \cdots & \cdots        & \cdots \\
               \beta_{I,1} & \cdots & \cdots
                                    & \cdots & \beta_{I,r_I}
             \end{array}
           \,\right]_{I\times r_I}\,.
   $$
   {\normalsize\it Then}} 
   there is a natural one-to-one correspondence between
   the set of $(A;B)$ defined above and the set of $S^1$-fixed-point
   components of $\HQuot_P({\cal E}^n)$.

  \item[{\rm (2)}]
   Let
   $$
    (A_1;B_1)=(A;B)\,,\; (A_2;B_2)\,,\, \cdots\,,\, (A_I;B_I)
   $$
   be a sequence of pairs of incomplete matrices defined by
   {\footnotesize
   $$
    A_i\;=\;\left[\,
             \begin{array}{lllll}
               \alpha_{i,1} & \cdots & \alpha_{i,r_i} \\
                \cdots      & \cdots & \cdots          & \cdots \\
               \alpha_{I,1} & \cdots & \cdots
                                     & \cdots & \alpha_{I,r_I}
             \end{array}
           \,\right]_{I\times r_I}
        \hspace{1em}\mbox{\normalsize\it and}\hspace{1em}
     B_i\;=\;\left[\,
             \begin{array}{lllll}
               \beta_{i,1} & \cdots & \beta_{i,r_i} \\
                \cdots     & \cdots & \cdots        & \cdots \\
               \beta_{I,1} & \cdots & \cdots
                                    & \cdots & \beta_{I,r_I}
             \end{array}
           \,\right]_{I\times r_I}\,.
   $$
   {\normalsize\it Then}} 
   $E_{(A,B)}$ admits a tower of fibrations that is compatible
    with the tower of fibratiosn of the flag manifold
    $X=\Fl_{r_1,\,\ldots,r_I}({\Bbb C}^n)$:
   $$
    \begin{array}{ccccccccc}
      E_{(A,B)}\,=\,E_{(A_1;B_1)}
       & \stackrel{f_1}{\longrightarrow}      & \cdots
       & \stackrel{f_{i-1}}{\longrightarrow}  & E_{(A_i;B_i)}
       & \stackrel{f_i}{\longrightarrow}      & \cdots
       & \stackrel{f_{I-1}}{\longrightarrow}  & E_{(A_I;B_I)} \\[.6ex]
      \downarrow\,\mbox{\scriptsize $g_1$}    & & \cdots &
       & \downarrow\,\mbox{\scriptsize $g_i$} & & \cdots  &
       & \downarrow\,\mbox{\scriptsize $g_I$}   \\
      X\,=\,\Fl_{r_1,\,\ldots,\,r_I}({\Bbb C}^n)
       & \stackrel{\underline{f}_1}{\longrightarrow}  & \cdots
       & \stackrel{\underline{f}_{i-1}}{\longrightarrow}
       & \Fl_{r_i,\,\ldots,\,r_I}({\Bbb C}^n)
       & \stackrel{\underline{f}_i}{\longrightarrow}  & \cdots
       & \stackrel{\underline{f}_{I-1}}{\longrightarrow}
       & \Gr_{r_I}({\Bbb C}^n)\,.
    \end{array}
   $$
  The fiber $F_{i,i+1}$ of
   $f_i:E_{(A_i;B_i)}\rightarrow E_{(A_{i+1;B_{i+1}})}$
   is the fibered product obtained in Lemma 2.1.10 (3)
   with $(A;B)$ in Lemma 2.1.10 (3) consisting of the $i$-th
   and the $(i+1)$-th rows of $(A;B)$ in Item {\rm (1)} above.
  We shall denote $E_{A_i;B_i}$ by $E_{(A;B)}^{(i)}$.

 \item[{\rm (3)}]
  In particular, $E_{(A,B)}$ is smooth for all $(A;B)$ in Item $(1)$
  above.
\end{itemize}
} 

\bigskip

\noindent
{\it Remark 2.2.2 $[\,$tower of fibrations$\,]$.}
 The tower of fibrations in Proposition 2.2.1, Item (2),
 is exactly the one induced from the tower of trivial fibrations
 $$
  \prod_{i=1}^I\,\Quot_{P_i}({\cal E}^n)\;
  \longrightarrow\;\prod_{i=2}^I\,\Quot_{P_i}({\cal E}^n)\;
  \longrightarrow\; \cdots\;
  \longrightarrow\; \Quot_{P_I}({\cal E}^n)
 $$
 and the inclusion
 $$
  \HQuot_P({\cal E}^n)
          \hookrightarrow \prod_{i=1}^I\,\Quot_{P_i}({\cal E}^n)\,.
 $$

\bigskip

\bigskip

\section{An exact computation of
    $\int_{X}\tau^{\ast}e^{H\cdot t}\cap {\mathbf 1}_d$
    from the mirror principle diagram.}

With the $S^1$-fixed-points locus in $\HQuot_P({\cal E}^n)$ understood
 in Sec.\ 2, we proceed now to compute the fundamental hypergeometric
 series $\HG[{\mathbf 1}]^X(t)$ reviewed in Sec.\ 1.
There are many technical details involved in the process and
 we study them in Sec.\ 3.1 - Sec.\ 3.5, following the logic order
 toward an exact expression in Sec.\ 3.6.

\bigskip

\bigskip

\subsection{\bf The extended Mirror Principle diagram and
     the distinguished $S^1$-fixed-point components in
     the hyper-Quot-scheme ${\cal Q}_d$.}

To make the comparison immediate, here we follow the notations
 in [L-L-Y1$\,$: III, Sec.\ 5.4].
Recall the following approach ibidem to compute $A(t)$
 when there is a commutative diagram$\,$:
$$
 \begin{array}{cccccl}
  F_0  & \stackrel{e^Y}{\longrightarrow}  & Y_0
       & \stackrel{g}{\longleftarrow}    & E_0  & \\[.6ex]
  \hspace{1ex}\downarrow\mbox{\scriptsize $i$}
     & & \hspace{1ex}\downarrow\mbox{\scriptsize $j$}
     & & \hspace{1ex}\downarrow\mbox{\scriptsize $k$} &  \\
  M_d  & \stackrel{\varphi}{\longrightarrow}  & W_d
       & \stackrel{\psi}{\longleftarrow} & {\cal Q}_d  & ,
 \end{array}
$$
where
 ${\cal Q}_d$ is an $S^1$-manifold,
 $\psi:{\cal Q}_d\rightarrow W_d$ is an $S^1$-equivariant resolution
  of singularities of $\varphi(M_d)$,
 $E_0$ is the set of fixed-points in $\psi^{-1}(Y_0)$ and is
  called the distinguished $S^1$-fixed-point locus, and
 $\varphi_{\ast}[M_d]=\psi_{\ast}[{\cal Q}_d]$ in $A_{\ast}^{S^1}(W_d)\,$.

In the case that $X$ is the Grassmannian manifold $\Gr_r({\Bbb C}^n)$,
 ${\cal Q}_d$ is the Quot-scheme $\Quot_{P(t)=(n-r)t+(d+n-r)}({\cal E}^n)$,
 and the linearized moduli space $W_d$ for $X$ is the projective space
 ${\Bbb P}(H^0(C,\,{\cal O}_C(d))\otimes \Lambda^r{\Bbb C}^n)$
 of $(\!\begin{array}{c}
         \mbox{\scriptsize $n$} \\[-1.2ex]
         \mbox{\scriptsize $r$}
        \end{array}\!)$-tuple of degree-$d$ homogeneous polynomials
 on $C$.
 This is a linearized moduli space for ${\Bbb P}(\Lambda^r{\Bbb C}^n)$
 that is turned into a linearized moduli space for $X$ via
 the Pl\"{u}cker embedding
 $\Gr_r({\Bbb C}^n)\rightarrow{\Bbb P}(\Lambda^r{\Bbb C}^n)$.
An element
 $[{\cal E}^n
   \rightarrow {\cal E}^n\!/\mbox{\raisebox{-.4ex}{${\cal V}$}}]$
 in ${\cal Q}_d$ can be represented by an $n\times r$-matrix $A_{\cal V}$
 of homogeneous polynomials in $z_0,\,z_1$ of degree $d$.
The map
 $\psi:\Quot_P({\cal E}^n) \rightarrow
    W_d={\Bbb P}(H^0(\CP^1,{\cal O}(d))\otimes\Lambda^r{\Bbb C}^n)$
 is given by taking the
 $(\!\begin{array}{c}
      \mbox{\scriptsize $n$} \\[-1.2ex]
      \mbox{\scriptsize $r$}
     \end{array}\!)$-tuple of $r\times r$-minors of
 $A_{\cal V}$.
 From this we deduce that the distinguished $S^1$-fixed-point
 components are exactly those labelled by admissible
 $(\alpha_{\bullet}\,;\,0_{\bullet})$.
(Cf.\ See [L-L-L-Y: Sec.\ 3.1] for more details and
      some related references.)

In the current case,
 $X$ is the flag manifold $\Fl_{r_1,\,\ldots,\,r_I}({\Bbb C}^n)$
 and the following two natural embeddings
 $$
   \iota_1\,:\,
          \Fl_{r_1,\,\ldots,\,r_I}({\Bbb C}^n)\;
             \longrightarrow\;
              \Gr_{r_1}({\Bbb C}^n)\,\times\,\cdots\,
                                \times\,\Gr_{r_I}({\Bbb C}^n)
 $$
 and
 $$
  \iota_2\,:\,
   \HQuot_P({\cal E}^n)\;
    \longrightarrow\;
      \Quot_{P_1}({\cal E}^n)\,
        \times\; \cdots\;\times\, \Quot_{P_I}({\cal E}^n)
 $$
 give rise to the following choices of spaces and morphisms for
 the diagram at the beginning of the subsection:
 (To save burden of notations, the projection map of a product to
  its $i$-th component will be denoted by $\pr_i$, regardless of
  which space is in question.)
\begin{itemize}
 \item[(1)]
  The embedding
  $$
   \tau\,
    =\,(\tau_1\circ\pr_1\circ\iota_1,\,\ldots,\,
        \tau_I\circ\pr_I\circ\iota_1)\;:\;
     X\,=\,\Fl_{r_1,\,\ldots,\,r_I}({\Bbb C}^n)\;
     \longrightarrow\;
      Y\,=\,\prod_{i=1}^I\CP^{{n\choose r_i}-1}
  $$
  where
  $\tau_i:\Gr_{r_i}({\Bbb C}^n)
          \rightarrow \CP^{{n\choose r_i}-1 }$
  is the Pl\"{u}cker embedding.

 \item[(2)]
  ${\cal Q}_d$ is the hyper-Quot-scheme $\HQuot_P({\cal E}^n)$,
  where $P=(P_1,\,\ldots,P_I)$ with $P_i=P_i(t)=(n-r_i)t+d_i+(n-r_i)$.

 \item[(3)]
  The linearized moduli space $W_d$ for $X$ is the product of
   projective spaces
   $$
    W_d\;
     =\;\prod_{i=1}^I\, W_{d_i}\;
     =\;\prod_{i=1}^I\,
     {\Bbb P}(H^0(C,\,{\cal O}_C(d_i))\otimes \Lambda^{r_i}{\Bbb C}^n)
   $$
   of $(\!\begin{array}{c}
           \mbox{\scriptsize $n$} \\[-1.2ex]
           \mbox{\scriptsize $r_i$}
         \end{array}\!)$-tuple
   of degree-$d_i$ homogeneous polynomials on $C$.
 This is a linearized moduli space for
  $\prod_{i=1}^I{\Bbb P}(\Lambda^{r_i}{\Bbb C}^n)$
  that is turned into a linearized moduli space for $X$ via
  $\tau:X\rightarrow Y$.
\end{itemize}

\begin{itemize}
 \item[(4)]
  The collapsing morphism
   $\varphi=(\varphi_1\circ\iota_3,\,\ldots,\, \varphi_I\circ\iota_3)$,
   where
   $$
    \iota_3\,:\, M_d(X)\;
     \hookrightarrow\;
      M_{d_1}(\Gr_{r_1}({\Bbb C}))\times\,
            \cdots\,\times M_{d_I}(\Gr_{r_I}({\Bbb C}^n))
   $$
   is the embedding induced by $\iota_1$ and
   $$
    \varphi_i\,:\, M_d(\Gr_{r_i}({\Bbb C}^n))\;
     \longrightarrow\;
     {\Bbb P}(H^0(C,\,{\cal O}_C(d_i))\otimes \Lambda^{r_i}{\Bbb C}^n)
   $$
   is the collapsing morphism when $X=\Gr_{r_i}({\Bbb C}^n)$.

 \item[(5)]
  The morphism
   $\psi=(\psi_1\circ\pr_1\circ\iota_2,\,\ldots,\,
          \psi_I\circ\pr_I\circ\iota_2)$,
   where
   $$
    \psi_i\,:\,\Quot_{P_i}({\cal E}^n)\;\longrightarrow\;
    {\Bbb P}(H^0(C,\,{\cal O}_C(d_i))\otimes \Lambda^{r_i}{\Bbb C}^n)
   $$
   is the map $\psi$ when $X$ is $\Gr_{r_i}({\Bbb C}^n)$ and the degree
   of curves in question is $d_i$, as reviewed in the beginning of
   this subsection.
\end{itemize}

\bigskip

\noindent
{\bf Lemma 3.1.1 [identical image].} {\it
 $\,\varphi(M_d)=\psi({\cal Q}_d)\,$ in $W_d$ and
 $\psi:{\cal Q}_d\rightarrow \varphi(M_d)$
  is a resolution of singularities of $\varphi(M_d)$.
} 

\bigskip

\noindent
{\it Proof.}
 On the stable map side,
  $\CP^1\times \Fl_{r_1,\,\ldots,\,r_I}({\Bbb C}^n)$ is nonsingular,
   projective, and convex;
  thus, $M_d$ contains the space
   $W_d^0
    := \Mor_{(1,d)}(\CP^1,\CP^1\times\Fl_{r_1,\ldots,r_I}({\Bbb C}^n))$
   of morphisms from $\CP^1$ to
   $\CP^1\times \Fl_{r_1,\,\ldots,\,r_I}({\Bbb C}^n)$
   of degree $(1,d)$ as an open dense subset
   (e.g. [F-P]).
 On the coherent sheaf side, ${\cal Q}_d$ contains an open dense subset
  ${\cal Q}_d^0$ that consists of sequences of vector bundle inclusions
  ${\cal V}_1\hookrightarrow\,\cdots\,\hookrightarrow {\cal V}_I
                                      \hookrightarrow {\cal E}^n$
  (e.g.\ [CF1]).
 Consequently, we only need to show that
  $\varphi(M_d^0)=\psi({\cal Q}_d^0)$ in $W_d$.
 Furthermore, since $M_d^0$ is a smooth Deligne-Mumford stack and
  $W_d^0$ is a smooth scheme, we only need to check the above identity
  at the level of atlas variety and, hence, only on the related sets of
  closed points.

 Let
  $\Fl=\Fl_{r_1,\,\ldots,\,r_I}({\Bbb C}^n)$ and
  $S_1\hookrightarrow \,\cdots\,\hookrightarrow S_I
      \hookrightarrow {\cal O}_{Fl}\otimes{\Bbb C}^n$
  be the tautological subbundles on $\Fl$.
 The map $\varphi: M_d^0\rightarrow W_d$ can be described as follows:
 The dual vector bundle morphisms:
  ${\cal O}_{Fl}\otimes{\Bbb C}^n\rightarrow S_i^{\vee}$
  induces a morphism
  ${\cal O}_{Fl}\times\Lambda^{r_i}({\Bbb C}^n)\rightarrow\Det S_i^{\vee}$
  produce the Pl\"{u}cker embedding of $\Gr_{r_i}({\Bbb C}^n)$ in
   $\CP^{\,{n\choose r_i}-1}$.
 The universal $\Delta$-collection of $\CP^{\,{n\choose r_i}-1}$
  as a toric variety ([Cox]) pulls back to a $\Delta$-collection on
  $\Gr_{r_i}({\Bbb C}^n)$ and then on
  $\Fl_{r_1,\,\ldots,\,r_I}({\Bbb C}^n)$.
 This $\Delta$-collection is given exactly by the morphism
  ${\cal O}_{Fl}\otimes \Lambda^{r_i}{\Bbb C}^n$ above.
  (Here, the comparison data $\{\,c_m\,\}_m$ in [Cox] of different line
   bunles in a $\Delta$-collection is implicit by consider only
   $\Det S_i^{\vee}$ instead of $n\choose r_i$ isomorphic copies of it.)
 Given a morphism
  $f:\CP^1
    \rightarrow \CP^1\times\Fl_{r_1,\,\ldots,\,r_I}({\Bbb C}^n)$
  of degree $(1,d)$, one has then a pull-back $\Delta$-collection
  on $\CP^1$ with $\Delta$ the fan associated to the toric manifold
  $\CP^1\times\CP^{{n\choose r_1}-1}\times\,\cdots\,
                                    \times\CP^{{n\choose r_I}-1}$.
 The morphism $\varphi$ is then the collapsing morphism,
  as constructed in [L-L-Y1, II] for toric variety
  and elaborated in [L(CH)-L-Y], that relates maps to a toric variety
  to the linearized moduli space $W_d$ via $\Delta$-collections of [Cox].

 On the other hand, the map $\psi:{\cal Q}_d^0\rightarrow W_d$ can be
  described as follows.
 Fix a presentation $\CP^1=\Proj{\Bbb C}[z_0,z_1]$ and identify
  ${\cal O}_{Fl}\otimes{\Bbb C}^n$ with ${\Bbb C}[z_0,z_1]^{\oplus n}$.
 Then an element $[{\cal V}_{\bullet}] \in {\cal Q}_d$ is identified
  with an $I$-tuple of ${\Bbb C}[z_0,z_1]$-valued matrices
  $(A_{{\cal V}_1},\,\cdots\,A_{{\cal V}_I})$
  as reviewed in the beginning of this subsection.
 Taking the minors of these matrices with appropriate order
  gives the map $\psi:{\cal Q}_d\rightarrow W_d$.

 Now, the matrices $A_{{\cal V}_i}$ that appear in the definition of
  $\psi$ corresponds to the inclusion
  ${\cal V}_i\hookrightarrow {\cal O}\otimes{\Bbb C}^n$ on $\CP^1$.
 For $[{\cal V}_{\bullet}]\in {\cal Q}_d^0$,
  taking minors corresponds to a morphism
   $\Det{\cal V}_i\rightarrow {\cal O}\otimes\Lambda^{r_i}{\Bbb C}^n$
   on $\CP^1$.
  The dual of these,
   ${\cal O}\otimes\Lambda^{r_i}{\Bbb C}^n
                            \rightarrow \Det{\cal V}_i^{\vee}$,
   is associated to taking the minors of the transpose $A_{{\cal V}_i}^T$
   of $A_{{\cal V}_i}$ in the corresponding order.
  They constitute a $\Delta$-collection on $\CP^1$ that corresponds
   to a map
   $\CP^1\rightarrow
     \CP^{{n\choose r_1}-1}\times\,\cdots\,\CP^{{n\choose r_I}-1}$
   that factors through the morphism
    $f_{{\cal V}_{\bullet}}:\CP^1\rightarrow \Fl$ associated
    to $[{\cal V}_{\bullet}]$.
 It follows from the key steps and ingredients reviewed above
  in the construction of $\varphi$ and $\psi$ that
  $\varphi(f_{{\cal V}_{\bullet}}) = \psi([{\cal V}_{\bullet}])$ in $W_d$
  since the minors of $A_{{\cal V}_{\bullet}}$ and
  $A_{{\cal V}_{\bullet}}^T$ give the same tuple when taken in
  a consistent respective order.
 Conversely, an $f\in M_d^0$ induces a
  $[{\cal V}^f_{\bullet}]\in {\cal Q}_d^0$ by pulling back the
  tautological subbundles $S_i$ on $\Fl$ and one can check that
  $\varphi(f)=\psi([{\cal V}^f_{\bullet}])$ in $W_d$ as well.

 Since the correspondence
  $M_d^0 \rightarrow {\cal Q}_d^0$ with
  $f\mapsto [{\cal V}^f_{\bullet}]$ is surjective,
 this shows that $\varphi(M_d^0)=\psi({\cal Q}_d^0)$ at the set of
  closed points of the domain stack/variety of $\varphi$ and $\psi$
  respectively and, hence, at the whole domain stack/variety.

 Since ${\cal Q}_d$ is smooth, to show that
  $\psi:{\cal Q}_d\rightarrow \varphi(M_d)$ is a resolution
   of singularities of $\varphi(M_d)$, one only have to show that
   the morphism $\psi:{\cal Q}_d\rightarrow \varphi(M_d)$ is birational.
  But this follows from the details of $\varphi$ and $\psi$ reviewed
   above that both ${\cal Q}_d$ and $\varphi(M_d)$ contains
   a copy of $\Mor(\CP^1,\Fl_{r_1,\,\ldots\,r_I}({\Bbb C}^n))$
   as an open dense subset and that the restriction of $\psi$ to this
   subset is the identity map.

 This concludes the proof.

\noindent\hspace{14cm}$\Box$

\bigskip

Consequently, Lemma 5.5 of [L-L-Y1, III] holds and one can compute
 the fundamental hypergeometric series $\HG[{\mathbf 1}]^X(t)$ for
 $X=\Fl_{r_1,\,\ldots,\,r_I}({\Bbb C}^n)$ via the localization method
 on the hyper-Quot-scheme ${\cal Q}_d$ rather than on $W_d$ directly.

The first task now is to identify the distinguished
 $S^1$-fixed-point components $E_0$ in ${\cal Q}_d$ that appear in
 the extended Mirror Principle diagram..
Since the linearized moduli space $W_d$ is a product of $W_{d_i}$,
 its $S^1$-fixed-point components $Y_u$ must also come from products
 of $S^1$-fixed-point components in $W_d$:
 $$
  Y_u\;=\; Y_{u_1}\times\,\cdots\,\times Y_{u_I}\,,
 $$
 where $u=(u_1,\,\ldots,\,u_I)$, $0\le u_i\le d_i$, and
 $Y_{u_i}$ is the $S^1$-fixed-point component in $W_{d_i}$ labelled
 by $u_i$, as given in [L-L-L-Y: Sec.\ 3.1].
The $Y_0$ of the Mirror Principle diagram, given in the beginning of
 this subsection corresponds to $u=(0,\,\ldots,\,0)$.

\bigskip

\noindent
{\bf Lemma 3.1.2 [image of $E_{(A;B)}$].}  {\it
 $\;\psi(E_{(A;B)})\,\subset\, Y_{u}\;$ with
 $\;u_i=\beta_{i,1}+\,\cdots\, +\beta_{i,r_i}$.
} 

\bigskip

\noindent
{\it Proof.}
 Under the embedding
  $$
   \iota_2\,:\,
    \HQuot_P({\cal E}^n)\;
     \longrightarrow\;
       \Quot_{P_1}({\cal E}^n)\,
         \times\; \cdots\;\times\, \Quot_{P_I}({\cal E}^n)\,,
  $$
  $E_{(A;B)}$ (is the unique $S^1$-fixed-point component in
  $\HQuot_P({\cal E}^n)$ that) goes into the $S^1$-fixed-point component
  $E_{(\alpha_{1,\bullet}\,;\,\beta_{1,\bullet})}\times\,
   \cdots\,\times E_{(\alpha_{I,\bullet}\,;\,\beta_{I,\bullet})}$
  in
  $\Quot_{P_1}({\cal E}^n)\,
    \times\; \cdots\;\times\, \Quot_{P_I}({\cal E}^n)$.
 From [L-L-L-Y: Lemma 3.1.1],
  $\psi_i(\Quot_{P_i}({\cal E}^n))\,\subset\,
              Y_{\beta_{i,1}+\,\cdots\,+\beta_{i,r_i}}$.
 This implies the lemma.

\noindent\hspace{14cm}$\Box$

\bigskip

\noindent
{\bf Corollary 3.1.3 [distinguished components].} {\it
 The distinguished $S^1$-fixed-point locus $E_0$ in
  the Mirror Principle diagram is given by
  $$
   E_0\;=\;\coprod_A\, E_{(A;0)}\,,
  $$
  where $A$ runs over all the incomplete matrices $(\alpha_{i,j})_{i,j}$
  with the $i$-th row $(\alpha_{i,\bullet})$ being a partition
  $0\le \alpha_{i,1}\le \,\cdots\,\le\alpha_{i,r_i}$ of $d_i$
  into non-negative integers of length $r_i$ and
  $\alpha_{i,j}\ge \alpha_{i+1,j}$ for $j=1,\,\ldots,\,r_i$.
 The set of such $A$ is always non-empty.
  {\rm (}Cf.\ {\sc Figure} {\rm 3-1-1}.{\rm )}
} 

\begin{figure}[htbp]
 \setcaption{{\sc Figure} 3-1-1.
  \baselineskip 14pt
   The set of distinguished $S^1$-fixed-point component $E_{(A;0)}$
    in $\HQuot_P({\cal E}^n)$, with $P=(P_1,\,\ldots,\,P_I)$,
    $P_i=P_i(t)=(n-r_i)t+d_i+(n-r_i)$
    is the same as the set of Young tableaux whose entries satisfy
    some monotonous conditions both horizontally and vertically.
   Illustrated here is such a set for $\Fl_{2,3}({\Bbb C}^n)$
    with $g=0$ stable maps of multiple degree $(d_1,d_2)=(2,6)$
    in consideration.
 } 
 \centerline{\psfig{figure=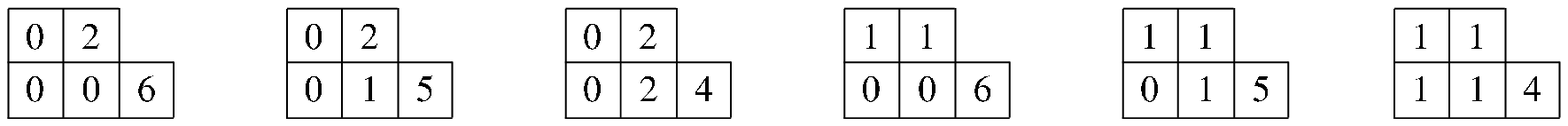,width=13cm,caption=}}
\end{figure}

\bigskip

We will come back in Sec.\ 3.4 to work out the hyperplane-induced
 classes on $E_{(A;0)}$ from the map $\psi$ after understanding
 $E_{(A;0)}$ better.

\bigskip

\bigskip

\subsection{Bundles with fiber restrictive flag manifolds
            and the class $\Omega({\cal P}_{\bullet})$.}

The discussion in the previous subsection together with
 Proposition 2.2.1 implies that a distinguished $S^1$-fixed-point
 component $E_{(A;0)}$ in $\HQuot_P({\cal E}^n)$ admits
 a tower of fibrations by restrictive flag manifolds.
For this reason and quantities that will be needed in later
 subsections, we digress in this subsection to take a look at bundles
 with fiber restrictive flag manifolds and their associated bundle
 with fiber ordinary flag manifolds.

Let
 ${\cal E}$ be a vector bundle of rank $n$ on a smooth variety $Y$,
 $k_{\bullet}: 1\le k_1< \,\cdots\,< k_s< n$ with $k_i$ integers, and
 ${\cal P}_{\bullet}: {\cal P}_1\hookrightarrow\,\cdots\,\hookrightarrow
  {\cal P}_s\hookrightarrow {\cal E}$ be a (not necessarily strict)
   inclusion sequence of vector subbundles of ${\cal E}$.
Let $l_i$ be the rank of ${\cal P}_i$.
Then associated to the triple $({\cal E},k_{\bullet},{\cal P}_{\bullet})$
 is a bundle $g:W\rightarrow Y$ over $Y$, whose fiber $W_y$ over $y\in Y$
 is the restrictive flag manifold
 $\Fl_{k_1,\,\ldots,\,k_s}({\cal E}_y,{\cal P}_{\bullet,y})$,
  where $(\,\cdot\,)_y$ denotes the fiber of $(\,\cdot\,)$ at $y$.
By construction, $W\rightarrow Y$ is naturally a subbundle of
 a flag manifold bundle $g^{\prime}:W^{\prime}\rightarrow Y$,
 whose fiber $W^{\prime}_y$ at $y$ is the flag manifold
 $\Fl_{k_1,\,\ldots,\,k_s}({\cal E}_y)$.
Denote the tautological inclusion $W\hookrightarrow W^{\prime}$
 over $Y$ by $\iota$. By construction, $g=g^{\prime}\circ\iota$.
Let ${\cal S}_j$ (resp.\ ${\cal S}_j^{\prime}$), $j=1,\,\ldots,\,s$,
 be the tautological bundles on $W$ (resp.\ $W^{\prime}$), whose fiber
 at $w\in g^{-1}(y)$ (resp.\ $w^{\prime}\in {g^{\prime}}^{-1}(y)$)
 is the $j$-th element of the flag
 $w\in\Fl_{k_1,\,\ldots,\,k_s}({\cal E}_y,{\cal P}_{\bullet,y})$
 (resp.\ $w^{\prime}\in\Fl_{k_1,\,\ldots,\,k_s}({\cal E}_y)$).
Note that $\iota^{\ast}{\cal S}^{\prime}_{\bullet}={\cal S}_{\bullet}$.

Define $\Omega({\cal P}_{\bullet})$ to be (the Poincar\'{e} dual of)
 $[W]$ in $A_{\ast}(W^{\prime})$.
We will discuss below how to express $\Omega({\cal P}_{\bullet})$
 in terms of the Chern roots of ${\cal S}_j^{\prime}$ and ${\cal P}_j$.

\bigskip

\begin{flushleft}
{\bf A pedagogic discussion of a basic excess example.}
\end{flushleft}
Consider the sequence of bundle morphisms obtained from composition
 $$
  \varphi_j\; :\; {\cal S}_j^{\prime}\;
     \hookrightarrow\; g^{\prime\ast}{\cal E}\;
     \longrightarrow\;
     g^{\prime\ast}{\cal E}/
              \mbox{\raisebox{-.4ex}{$g^{\prime\ast}{\cal P}_j$}}
 $$
 of ${\cal O}_{W^{\prime}}$-modules.
Then $\Ker\varphi_j$ is an ${\cal O}_{W^{\prime}}$-module,
 which is not locally free in general.
The reduced scheme associated to the intersection $Z:=\cap_{j=1}^s\,Z_j$
 of the minimal stratum $Z_j$ of the flattening stratification of
 $\Ker\varphi_j$ is exactly $W$ and
 $(\Ker\varphi_j)|_W=\iota^{\ast}{\cal S}^{\prime}_j$.

\bigskip

\noindent
{\bf Lemma 3.2.1 [$Z=W$ as schemes].} {\it
 $Z$ is reduced and hence $Z=W$ as schemes.
} 

\bigskip

\noindent {\it Proof.}
 Since $Z_j$ is the minimal stratum of the flattening stratification
  of $\Ker\varphi_j$, it is a closed subscheme of $W^{\prime}$ defined
  by the Fitting ideal sheaf ${\cal I}_j$ whose local sections are
  generated by the set of all entries of any local presentation of
  the ${\cal O}_{W^{\prime}}$-module moprhism $\varphi_j$.
 Since the problem is local, we may assume
  that $Y=\Spec R$, where $R$ is a local ring with the residue field
   extending ${\Bbb C}$, and
  that ${\cal E}$ and ${\cal P}_{\bullet}$ are free $R$-modules
   with ${\cal E}=\oplus R$ and ${\cal P}_j$ being
   the direct sum of the first $l_j$ direct summands in
   the decomposition of ${\cal E}$.
 Under such specification, the quotients
  ${\cal E}/\mbox{\raisebox{-.4ex}{${\cal P}_j$}}$ are realized
  as the direct sum of the the last $(n-l_j)$ direct summand in the
  decomposition of ${\cal E}$.

 Given a point $w^{\prime}$ in $W^{\prime}$, after
  a $\GL_n(R)$-transformation $M=(r_{\mu\nu})_{n\times n}$
  one can represent $w^{\prime}$ by the $n\times k_s$ matrix
  $[w^{\prime}]$ obtained from enlarging the $k_s\times k_s$
  diagonal matrix
  $\Diag[\,I_{k_1},\, I_{k_2-k_1},\,\cdots,\, I_{k_s-k_{s-1}}\,]$,
  where $I_{\bullet}$ is the $\bullet\times\bullet$ identity matrix,
  to an $n\times k_s$ matrix by adding zero entries.

 In terms of this presentation, $w^{\prime}$ is contained in an affine
  chart
  $$
   U\; :
    =\; \Spec R\,[x_{\mu\nu}\,:\, 1\le \nu\le k_s,\,
                  k_j< \mu\le n\;\;
                      \mbox{if $k_{j-1}< \nu \le k_j$,
                               $j=1,\,\ldots,\,s$}\, ]\,,
  $$
  where $x_{\mu\nu}$ are indeterminants that appear as the entries of
  the block lower triangular matrix determined by $[w^{\prime}]$,
  as indicated in {\sc Figure} 3-2-1.
 \begin{figure}[htbp]
  \setcaption{{\sc Figure} 3-2-1.
   \baselineskip 14pt
   The lower-triangular matrices determined by $[w^{\prime}]$.
  } 
  \centerline{\psfig{figure=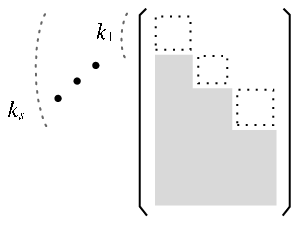,width=13cm,caption=}}
 \end{figure}
 Let $\Xi$ be the $n\times k_s$ matrix
  $([w^{\prime}]_{\mu\nu})+(x_{\mu\nu})$.
 Then, over $U$,
  ${\cal S}^{\prime}_j$ is generated by the column vectors of $\Xi$,
  $g^{\prime\ast}{\cal E}/
             \mbox{\raisebox{-.4ex}{$g^{\prime\ast}{\cal P}_j$}}$
   is generated by the column vectors of $M^{(j)}$,
   the lower $(n-l_j)\times n$ part of the matrix $M$, and
  $\varphi_j$ are represented by the matrix,
  (cf.\ the projection to the quotient
     $g^{\prime\ast}{\cal E}/g^{\prime\ast}{\cal P}_j$)
   $$
      (M\,\Xi)^{(j)}\,,
   $$
   the lower $(n-l_j)\times k_s$ part of the $n\times k_s$ matrix $M\Xi$.
 Thus on $U$, the ideal sheaf of $Z$ is generated by
  $k_1(n-l_1)+k_2(n-l_2)+\,\cdots\,+k_s(n-l_s)$-many
  degree-$1$ elements in the polynomial ring
  $$
   R\,[\,x_{\mu\nu}\,:\, 1\le \nu\le k_s,\, k_j< \mu\le n\;\;
                      \mbox{if $k_{j-1}< \nu \le k_j$,
                               $j=1,\,\ldots,\,s$}\, ]
  $$
  over $R$.
 Since the intersection of linear subvarieties is always smooth,
  this implies that $Z\cap U$ is smooth over $Y$.
 Since $w^{\prime}$ is arbitrary, this shows that $Z$ is smooth over
  $Y$; in particular, it is reduced. Consequently, $Z=W$ as schemes.

\noindent\hspace{14cm}$\Box$

\bigskip

\noindent
In this way, $W$ is realized as the degeneracy locus of the bundle
 morphism
 $$
  \phi\,:=\,(\varphi_1,\,\cdots,\,\varphi_s)\;:\;
  {\cal S}^{\prime}\,
    :=\,\oplus_{j=1}^s\,{\cal S}_j^{\prime}\;
  \longrightarrow\;
  {\cal P}^{\prime}\,
    :=\,\oplus_{j=1}^s\,g^{\prime\ast}{\cal E}/g^{\prime\ast}{\cal P}_j
 $$
 on $W^{\prime}$, over which the rank of $\phi$ is $0$.

Let us perform some dimension count:
 if $\phi$ were generic while sending ${\cal S}_j^{\prime}$ into
  $g^{\prime\ast}{\cal E}/g^{\prime\ast}{\cal P}_j$, its minimal
  degeneracy locus has codimension in $W^{\prime}$ equal to
  $$
   k_1(n-l_1)+k_2(n-l_2)+\,\cdots\,+k_s(n-l_s);
  $$
  on the other hand the codimension of $W$ in $W^{\prime}$
  is the same as the codimension of
  $\Fl_{k_1,\,\ldots,\, k_s}({\Bbb C}^n,\Pi_{\bullet})$
  in $\Fl_{k_1,\,\ldots,\,k_s}({\Bbb C}^n)$, which is
  $$
   k_1(n-l_1)+(k_2-k_1)(n-l_2)+\,\cdots\,+(k_s-k_{s-1})(n-l_s)\,.
  $$
Thus we see that $W$ is an excess degeneracy locus of
 $\phi$ and we cannot apply the Thom-Porteous formula directly
 to represent $[W]$ in $W^{\prime}$ in terms of Chern classes
 of ${\cal S}^{\prime}$ and ${\cal P}^{\prime}$.
We now discuss how to remedy this.

\bigskip

\begin{flushleft}
{\bf Removal of excess via nesting restrictions.}
\end{flushleft}
Consider the following sequence of morphisms
 $$
  \phi^{(j)}
   :=(\varphi_1,\,\cdots,\,\varphi_j)\;:\;
   {\cal S}^{\prime,(j)}\,
    :=\,\oplus_{j^{\prime}=1}^j\,{\cal S}_{j^{\prime}}^{\prime}\;
   \longrightarrow\;
   {\cal P}^{\prime,(j)}\,
    :=\,\oplus_{j^{\prime}=1}^j\,
        g^{\prime\ast}{\cal E}/g^{\prime\ast}{\cal P}_{j^{\prime}}
 $$
 and the associated sequence of minimal degeneracy locus $W_{(j)}$.
Then
 $$
  W^{\prime} =: W_{(0)}\,\supset\,
  W_{(1)}\,\supset\, W_{(2)}\,\supset\,\cdots\,\supset W_{(s)}\,=\,W
 $$
and, similar to the previous discussion, all the $W_{(j)}$ are smooth;
indeed they are all restrictive flag manifold bundle over $Y$.
Observe that the codimension of $W_{(j)}$ in $W_{(j-1)}$ is
 $(k_j-k_{j-1})(n-l_j)$.

Consider the morphism
 $$
  \overline{\varphi}_j\; :\;
     {\cal S}_j^{\prime}/
              \mbox{\raisebox{-.4ex}{${\cal S}_{j-1}^{\prime}$}}\;
    \longrightarrow\;
     g^{\prime\ast}{\cal E}/
              \mbox{\raisebox{-.4ex}{$g^{\prime\ast}{\cal P}_j$}}
  \hspace{1em}\mbox{on $\;\;W_{(j-1)}$}
 $$
 induced from the restriction of $\varphi_j$ to related sheaves over
 $W_{(j-1)}$.
Though not well-defined on any bigger domain, $\overline{\varphi}_j$
 is well-defined on $W_{(j-1)}$ since $\varphi_j$ maps
 ${\cal S}_{j-1}^{\prime}$ into
 $g^{\prime\ast}{\cal P}_{j-1}\subset g^{\prime\ast}{\cal P}_j$
 when restricted to $W_{(j-1)}$.
The minimal degeneracy locus of
 $\overline{\varphi}_j$ on $W_{(j-1)}$ is exactly $W_{(j)}$.
Since the codimension of $W_{(j)}$ in $W_{(j-1)}$ is the same as
 $\rank(\,
      {\cal S}_j^{\prime}/
              \mbox{\raisebox{-.4ex}{${\cal S}_{j-1}^{\prime}$}}\,)\,
     \cdot\,
     \rank(\,
        g^{\prime\ast}{\cal E}/
              \mbox{\raisebox{-.4ex}{$g^{\prime\ast}{\cal P}_j$}}\,)$,
 $W_{(j)}$ is now a proper degeneracy locus of $\overline{\varphi}_j$,
 cf.\ {\sc Figure 3-2-2}.
 \begin{figure}[htbp]
  \setcaption{{\sc Figure} 3-2-2.
   \baselineskip 14pt
   A restrictive flag manifold bundle $W$ is not directly realizable
    as a critical degeneracy locus of a bundle morphism on
    the associated flag manifold bundle $W^{\prime}$.
   The nesting construction removes the excess degeneracy.
  } 
  \centerline{\psfig{figure=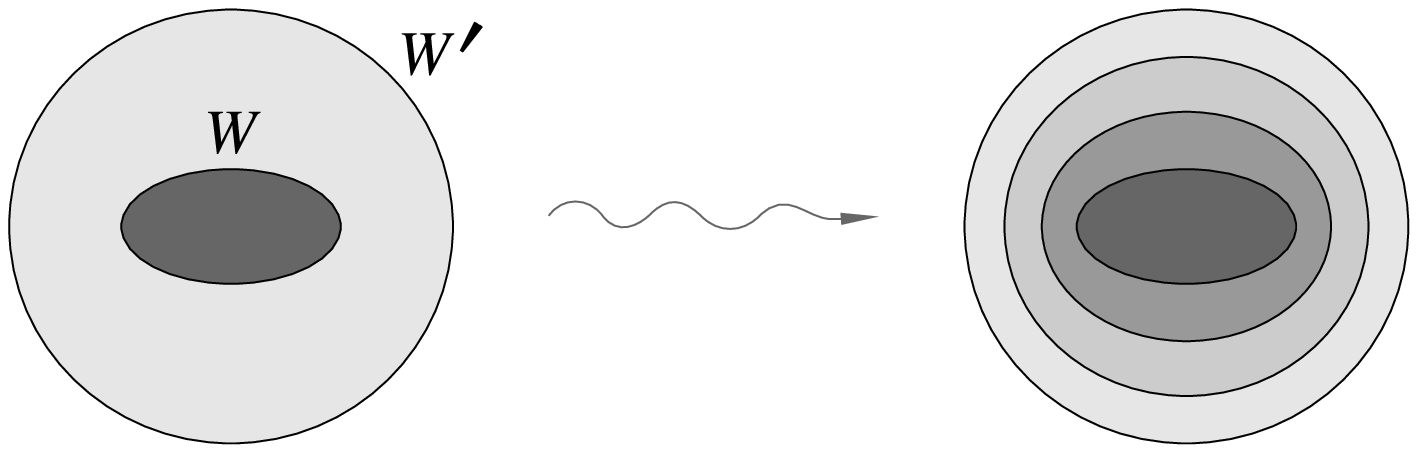,width=13cm,caption=}}
 \end{figure}
Thus, the Thom-Porteous formula applies and one can express
 (the Poincar\'{e} dual $\Omega_{(j)}$ of) $[W_{(j)}]$ in
 $A_{\ast}(W_{(j-1)})$ in terms of the Chern classes of
 ${\cal S}_j^{\prime}/
            \mbox{\raisebox{-.4ex}{${\cal S}_{j-1}^{\prime}$}}$
 and
 $g^{\prime\ast}{\cal E}/
              \mbox{\raisebox{-.4ex}{$g^{\prime\ast}{\cal P}_j$}}$
 and, hence, in terms of Chern roots of
 ${\cal S}_j^{\prime}/
            \mbox{\raisebox{-.4ex}{${\cal S}_{j-1}^{\prime}$}}$
 and ${\cal E}/{\cal P}_j$ via determinantal identities.
(Cf.\ [Fu1].)

Let
 $$
  \begin{array}{rcl}
   c({\cal S}_j^{\prime}/
            \mbox{\raisebox{-.4ex}{${\cal S}_{j-1}^{\prime}$}})(t)
     & =
     & \prod_{j^{\prime\prime}=1}^{k_j-k_{j-1}}
              (1+y_{j;j^{\prime\prime}}t)
       \hspace{3em}\mbox{and}            \\[1ex]
   c({\cal E}/\mbox{\raisebox{-.4ex}{${\cal P}_j$}})(t)
     & = & \prod_{j^{\prime}=1}^{n-l_j}(1+q_{j;j^{\prime}}t)
  \end{array}
 $$
 be the Chern polynomials of bundles involved in terms of their
 Chern roots.
It follows from [Fu1: Chapter 14 and Lemma A.9.1] that

 \begin{eqnarray*}
  \Omega_{(j)}
    & = & \Delta_{n-l_j}^{(k_j-k_{j-1})}
          \left(
           \frac{
            c\left(g^{\prime\ast}{\cal E}/
             \mbox{\raisebox{-.4ex}{$g^{\prime\ast}{\cal P}_j$}}\right)
                }{
              c\left({\cal S}_j^{\prime}/
               \mbox{\raisebox{-.4ex}{${\cal S}_{j-1}^{\prime}$}}\right)
               }
          \right)\;
    =\;  \Delta_{n-l_j}^{(k_j-k_{j-1})}
          \left(
           \frac{\prod_{j^{\prime}=1}^{n-l_j}(1+q_{j;j^{\prime}})
               }{\prod_{j^{\prime\prime=1}}^{k_j-k_{j-1}}
                                      (1+y_{j;j^{\prime\prime}})   }
          \right)            \\[.6ex]
  & = & \prod_{j^{\prime}=1}^{n-l_j}\,
          \prod_{j^{\prime\prime}=1}^{k_j-k_{j-1}}\,
           (q_{j;j^{\prime}}-y_{j;j^{\prime\prime}} )\,.
 \end{eqnarray*}

Recall now the natural morphism
 $$
  \cdot\,\Omega_{(j)}\;:\:  A^{\ast}(W_{(j)})\;
   \longrightarrow\; A^{\ast+ (k_j-k_{j-1})(n-l_j)}(W_{(j-1)})
 $$
 dual to the intersection product
 $$
  \cdot\, [W_{(j)}]\;:\;
    A_{\ast}(W_{(j-1)})\;
      \longrightarrow\; A_{\ast-(k_j-k_{j-1})(n-l_j)}(W_{(j)})\,,
 $$
 for $j=1,\,\ldots,\, s$.
Thus, the Poincar\'{e} dual of $[W_{(j)}]$ in $A^{\ast}(W_{(j-2)})$
 is given by $\Omega_{(j)}\cdot\Omega_{(j-1)}$.
Iterating this procedure, one concludes the following lemma:

\bigskip

\noindent
{\bf Lemma 3.2.2 [the class $\Omega({\cal P}_{\bullet})$].} {\it
 $($The Poincar\'{e} dual of$\,)$ $[W]$ in $W^{\prime}$ is given by
 $$
  \Omega({\cal P}_{\bullet})\;
  =\; \Omega_{(s)}\cdot\Omega_{(s-1)}\cdot\,\cdots\,\cdot
                                             \Omega_{(1)} \;
  =\; \prod_{j=1}^s\,
      \prod_{j^{\prime}=1}^{n-l_j}\,
      \prod_{j^{\prime\prime}=1}^{k_j-k_{j-1}}\,
       (q_{j;j^{\prime}}-y_{j;j^{\prime\prime}})
 $$
 in terms of the Chern roots of
  ${\cal S}_j^{\prime}/
         \mbox{\raisebox{-.4ex}{${\cal S}_{j-1}^{\prime}$}}$
  and
  ${\cal E}/\mbox{\raisebox{-.4ex}{${\cal P}_j$}}$,
  $j=1,\,\ldots,\, s$..
} 

\bigskip

\noindent
{\it Remark 3.2.3 $[$rational presentation of $\Omega$$]$.}
 Let
  $c({\cal E})(t) = \prod_{j^{\prime}=1}^n(1+e_{j^{\prime}}t)$ and
  $c({\cal P}_j)(t) = \prod_{j^{\prime}=1}^{l_j}(1+p_{j;j^{\prime}}t)$
  be the Chern polynomials of bundles involved in terms of their
  Chern roots.
 When the ratio makes sense, a rational presentation of $\Omega_{(j)}$
  (resp.\ $\Omega({\cal P}_{\bullet})$)
  is given by
  $$
   \Omega_{(j)}\;
    =\; \prod_{j^{\prime\prime}=1}^{k_j-k_{j-1}}\,
             \frac{ \prod_{j^{\prime}=1}^n\,
                      (e_{j^{\prime}}-y_{j;j^{\prime\prime}})
                 }{ \prod_{j^{\prime}=1}^{l_j}\,
                      (p_{j;j^{\prime}}-y_{j;j^{\prime\prime}}) }
   \hspace{1em}\mbox{(resp.}\hspace{2ex}
   \Omega({\cal P}_{\bullet})\;
     =\; \prod_{j=1}^s\,
         \prod_{j^{\prime\prime}=1}^{k_j-k_{j-1}}\,
            \frac{ \prod_{j^{\prime}=1}^n\,
                     (e_{j^{\prime}}-y_{j;j^{\prime\prime}})
                }{ \prod_{j^{\prime}=1}^{l_j}\,
                     (p_{j;j^{\prime}}-y_{j;j^{\prime\prime}}) }\,).
 $$

\bigskip

\bigskip

\subsection{Tautological sheaves on $E_{(A;0)}$ and
            $E_{(A;0)}\times\CPlarge^1$.}

\begin{flushleft}
{\bf Tautological sheaves and filtrations on $E_{(A;0)}$.}
\end{flushleft}
Recall from Proposition 2.2.1 the tower of fibrations of $E_{(A;0)}$
 $$
  E_{(A;0)}=E^{(1)}_{(A;0)}\;
   \stackrel{f_1}{\longrightarrow}\; \cdots\;
   \stackrel{f_{i-1}}{\longrightarrow}\;
  E_{(A;0)}^{(i)}\;\stackrel{f_i}{\longrightarrow}\; \cdots\;
   \stackrel{f_{I-1}}{\longrightarrow}\;
  E_{(A;0)}^{(I)}=\Fl_{m_{I,\bullet}}({\Bbb C}^n)
 $$
 with the fiber $F^{(i)}$ of $f_i$ being the restrictive flag manifold
 $\Fl_{m_{i,\bullet}}({\Bbb C}^n,\Pi_{i+1,\bullet})$.
The two systems of tautological vector bundles on
 $\Fl_{m_{i,\bullet}}({\Bbb C}^n,\Pi_{i+1,\bullet})$,
 one comes from the restriction of the sequence of tautological
 subbundles on $\Fl_{m_{i.\bullet}}({\Bbb C}^n)$ and the other
 from the data $\Pi_{i+1,\bullet}$,
 gives rise to two inclusion sequences of locally free sheaves
 ${\cal S}_{i,\bullet}$ and $f_i^{\ast}{\cal P}_{i+1,\bullet}$ on
 $E_{(A;0)}^{(i)}$
 with ${\cal S}_{i,j}\hookrightarrow f_i^{\ast}{\cal P}_{i+1,j}$,
 $j=1,\,\ldots,\,K_i$.
 From the discussion in Sec.\ 2.1, each ${\cal P}_{i+1,j}$
 is an ${\cal S}_{i+1,j^{\prime}}$ for some ${\cal S}_{i+1,j^{\prime}}$
 on $E_{(A;0)}^{(i+1)}$.
The pull-back of ${\cal S}_{i,\bullet}$ and
 $f_i^{\ast}{\cal P}_{i+1,\bullet}$ on $E_{(A;0)}^{(i)}$ to the whole
 $E_{(A;0)}$ will be denoted by $\widehat{\cal S}_{i,\bullet}$ and
 $\widehat{\cal P}_{i+1,\bullet}$ respectively.
Recall also the tautological embedding of $E_{(A;0)}$ into a product
 flag manifolds (cf.\ Sec.\ 3.1, the discussion before Lemma 3.1.2).
 It is good to keep in mind that
 both $\widehat{\cal S}_{i,j}$ and $\widehat{\cal P}_{i+1,j}$ are
 the restriction to $E_{(A;0)}$ of tautological bundles on this product.

In terms of Sec.\ 3.2,
 given $E_{(A;0)}^{(i+1)}$ with the tautological subbundles
  ${\cal S}_{i+1,\bullet}$,
 the smooth bundle map $f_i:E_{(A;0)}^{(i)}\rightarrow E_{(A;0)}^{(i+1)}$
  with fiber restrictive flag manifolds can be constructed as in that
  subsection from
  a sequence of integers $r_{i,1}<\,\cdots\,<r_{i,K_i}=r_i$
   determined by $(\alpha_{i,1},\,\cdots,\,\alpha_{i,r_i})$ and
  a sequence of subbundles
   ${\cal P}_{i+1,j}={\cal S}_{i+1,j^{\prime}}$ on $E_{(A;0)}^{(i+1)}$
   with $j^{\prime}$ determined by the submatrix from $A$:
   $$
    \left[\,
     \begin{array}{lllll}
       \alpha_{i,1}   & \cdots   & \alpha_{i,r_i}       \\[.6ex]
       \alpha_{i+1,1} & \cdots   & \alpha_{i+1,r_i}
                      & \cdots   & \alpha_{i+1,r_{i+1}}
     \end{array}\,
    \right]_{2\times r_{r+1}}\,.
   $$
For later use, denote this $j^{\prime}$ as the value $I_A(i,j)$
 of a function, denoted by $I_A$, on index pairs $(i,j)$.

\bigskip

\noindent
{\bf Definition 3.3.1 [index function asociated to $A$].} {\rm
 $I_A$ will be called the ({\it first}) {\it index function}
 associated to the Young tableau $A$.
} 

\bigskip

\noindent
{\it Remark 3.3.2 $[\,$more index functions$\,]$.}
 Later in Sec.\ 3.5, there will also be the {\it second} and
 the {\it third index function} that appear in the discussion.
 They are all determined by the Young tableau $A$ and will be
 denoted by $I^{\prime}_A$ and $I^{\prime\prime}_A$.

\bigskip

Introduce also the Chern roots of the components of the associated
 graded vector bundles on $E_{(A;0)}\,$:
 $$
  c\left(\widehat{\cal S}_{i,j}/\widehat{\cal S}_{i,j-1}\right)
   =\; \prod_{k=1}^{m_{i,j}}(1+y_{i,j; k})
  \hspace{1em}\mbox{
    for $i=1,\,\ldots,\,I$}\,,
 $$
 where
  $\widehat{\cal S}_{i,0}:=0$,
  $\widehat{\cal S}_{i,K_i+1}:={\cal O}_{E_{(A;0)}}\otimes{\Bbb C}^n$,
  and $m_{i,K_i+1};= n-r_i$.

\bigskip

\begin{flushleft}
{\bf Tautological sheaves and filtrations on $E_{(A;0)}\times\CP^1$.}
\end{flushleft}
Let
 $\pi_1: E_{(A;0)}\times \CP^1\rightarrow E_{(A;0)}$ and
  $\pi_2: E_{(A;0)}\times \CP^1\rightarrow \CP^1$
  be the two projection maps and
 ${\cal E}:=\pi_2^{\ast}\,{\cal E}^n$.
Then one has the following filtrations by locally free sheaves:
(Caution that this is {\it not} a double filtration.)
$$
 \begin{array}{ccccccccccccc}
   &&
   {\cal E}_{1,\bullet}             &    & \cdots  &
      & {\cal E}_{i,\bullet}        &    & \cdots  &
      & {\cal E}_{I,\bullet}        &    &                   \\[-.6ex]
   &&
   \cdot\cdot             &         &    &
      & \cdot\cdot        &         &    &
      & \cdot\cdot        &              &                   \\[1ex]
   {\cal E}_{\bullet}\;:
     && {\cal E}_1       & \hookrightarrow   & \cdots  & \hookrightarrow
      & {\cal E}_i  & \hookrightarrow   & \cdots  & \hookrightarrow
      & {\cal E}_I  & \hookrightarrow   & {\cal E}\,,          \\[.6ex]
   &&
   \|             &                   & \vdots  &
      & \|        &                   & \vdots  &
      & \|        &                   &                      \\[.6ex]
   &&
   {\cal E}_{1,K_1} &                   & \vdots  &
      & {\cal E}_{i,K_i}              & & \vdots  &
      & {\cal E}_{I,K_I}              & &                    \\[.6ex]
   &&
   \cup             &                   & \vdots  &
      & \cup        &                   & \vdots  &
      & \cup        &                   &                    \\[.6ex]
   &&
   \vdots           &                   & \vdots  &
      & \vdots      &                   & \vdots  &
      & \vdots      &                   &                    \\[.6ex]
   &&
   \cup             &                   & \vdots  &
      & \cup        &                   & \vdots  &
      & \cup        &                   &                    \\[.6ex]
   &&
   {\cal E}_{1,1}   &                   & \vdots  &
      & {\cal E}_{i,1}                & & \vdots  &
      & {\cal E}_{I,1}                & &                    \\[.6ex]
 \end{array}
$$
where
 the horizontal filtration $F_{\bullet}{\cal E} := {\cal E}_{\bullet}$
  on the first line comes from the universal filtration of ${\cal E}$
  on $\HQuot_P({\cal E}^n)\times\CP^1$
 while the vertical filtration
  $F_{\bullet}{\cal E}_i :={\cal E}_{i,\bullet}$ of ${\cal E}_i$
  is determined by the labelled flag $V_{i,\bullet}(s_{i,\bullet})$
  in ${\Bbb C}^n$ that characterizes ${\cal E}_i$.
 (In particular, the length of the vertical filtration
  $F_{\bullet}{\cal E}_i$ of ${\cal E}_i$ depends on $i$.)

 From the above diagram of various tautological sheaves on
 $E_{(A;0)}\times\CP^1$, one has the following two types of associated
 graded objects on $E_{(A;0)}\times\CP^1$:
\begin{itemize}
 \item[(1)]
  From the horizontal filtration:
  $\oplus_{i=1}^{I+1}{\cal E}_i/\mbox{\raisebox{-.4ex}{${\cal E}_{i-1}$}}$.
  \begin{itemize}
   \item[(1.1)]
    The quotient
     ${\cal E}_i/\mbox{\raisebox{-.4ex}{${\cal E}_{i-1}$}}$
     in general is not a direct sum of a locally free and a torsion
     ${\cal O}_{E_{(A;0)}\times\CPscriptsize^1}$-module.

   \item[(1.2)]
    For each $i$, there is a natural stratification of $E_{(A;0)}$
     by locally closed subsets in $E_{(A;0)}$.
     The strata are labelled by the isomorphism type of
      ${\cal E}_i/\mbox{\raisebox{-.4ex}{${\cal E}_{i-1}$}}$
      on each fiber $\CP^1$.
     When restricted to each stratum,
      ${\cal E}_i/\mbox{\raisebox{-.4ex}{${\cal E}_{i-1}$}}$
      is of the form of a direct sum of a locally free shaef and
      a torsion sheaf.
    There is a unique open stratum in this stratification.

   \item[(1.3)]
    By taking intersections of strata of the stratifications associated
     to different $i$, one obtains a stratification of of $E_{(A;0)}$
     that gives a minimal common refinement of all the stratifications
     in Item (1.2).
    From Item (1.2) this common refinement contains a unique open
     stratum.
    Set ${\cal E}_0=0$ and ${\cal E}_{I+1}={\cal E}$, then
     the graded sheaf
     $\oplus_{i=1}^{I+1}
      {\cal E}_i/\mbox{\raisebox{-.4ex}{${\cal E}_{i-1}$}}$
     is of the form of a direct sum of a locally free sheaf and a
     torsion sheaf over each stratum of the common stratification.
     (Cf.\ See [CF1], [CF2], and [M-M] for more related studies.)
 \end{itemize}
\end{itemize}

\begin{itemize}
 \item[(2)]
  For each vertical filtration,
   $$
    {\cal E}_{i,j}/\mbox{\raisebox{-.4ex}{${\cal E}_{i,j-1}$}}\;
    =\;
    \left(\,\pi_1^{\ast}\left(
            \widehat{\cal S}_{i,j}/
             \mbox{\raisebox{-.4ex}{$\widehat{\cal S}_{i,j-1}$}}\right)\,
    \right)
    (-a_{i,j}\,z)\,,
   $$
   where $z$ here represents the divisor $[E_{(A;0)}\times\{0\}]$
   on $E_{(A;0)}\times\CP^1$ (and will be omitted in the following
   discussion).
  The $S^1$-action on $\CP^1$ induces a natural $S^1$-action on
   the trivialized trivial bundle ${\cal E}$, which induces in turn
   an $S^1$-action on ${\cal E}_{i,\bullet}$ and hence on the graded
   bundle
   $\oplus_{j=1}^{K_i}
     {\cal E}_{i,j}/\mbox{\raisebox{-.4ex}{${\cal E}_{i,j-1}$}}$,
   where we set ${\cal E}_{i,0}=0$.
  These graded objects from the vertical filtration in the diagram
   will play crucial roles in our later discussions.
\end{itemize}

\bigskip

Recall also the restricting bundles $\Pi_{\bullet}$ in the discussion
 that are related to restrictvie flag manifolds.
Since $\widehat{\cal P}_{i+1,j}=\widehat{\cal S}_{i+1,I_A(i,j)}$,
 they are covered in the above discussion.
In particular,
 $\widehat{\cal P}_{i+1,K_i}=\widehat{\cal S}_{i+1,K_{i+1}}$,
 ($\,$i.e.\ $I_A(i,K_i)=K_{i+1}\,$) always holds.

\bigskip

\bigskip

\subsection{The hyperplane-induced classes on $E_{(A;0)}$.}

Recall from Sec.\ 3.1 the commutative diagram
 $$
  \begin{array}{ccccc}
   E_{(A;0)} \
    & \stackrel{k}{\hookrightarrow}
    & {\cal Q}_d   & \stackrel{\psi}{\longrightarrow} & W_d \\
   \downarrow & & \downarrow & & \| \\
   \prod_{i=1}^I\,E_{(\alpha_{i,\bullet};0)}
    & \longrightarrow & \prod_{i=1}^I\,Quot_{P_i}({\cal E}^n)
    & \longrightarrow & \prod_{i=1}^I\, W_{d_i}(\Gr_{r_i}({\Bbb C}^n))\,.
  \end{array}
 $$
The following lemma follows from the result of [L-L-L-Y: Sec.\ 3.1]
 for the hyperplane-induced class in the case of Grassmannian manifolds
 and the discussion and notations of the tautological bundles on
 $E_{(A;0)}$ in Sec.\ 3.3.

\bigskip

\noindent
{\bf Lemma 3.4.1 [hyperplane-induced classes].} {\it
 Let $\kappa_1,\,\ldots,\,\kappa_I$ be the hyperplane classes
  on $W_d$ from the product structure
  $W_d=\prod_{i=1}^I\,W_{d_i}(\Gr_{r_i}({\Bbb C}^n))$.
 Then
  $$
   k^{\ast}\psi^{\ast}\kappa_i\;
    =\; - c_1(\widehat{\cal S}_{i,K_i}).
  $$
 In terms of Chern roots of $\widehat{\cal S}_{i,j}$,
  $$
   k^{\ast}\psi^{\ast}\kappa_i\;
    =\; -(\,    y_{i,1;1} + \,\cdots\,+ y_{i,1; m_{i,1}}\,
            +\, \cdots\,
            +\, y_{i,K_i;1} + \,\cdots\,+ y_{i,K_i; m_{i,K_i}}\,
          \,)\,.
  $$
} 

\bigskip

\subsection{An exact computation of
   $e_{{\tinyBbb C}^{\times}}(E_{(A;0)}/\HQuot_P({\cal E}^n))$.}

In this subsection we work out an exact expresion of
 $e_{{\tinyBbb C}^{\times}}(E_{(A;0)}/\HQuot_P({\cal E}^n))$
 in terms of Chern roots $y_{i,j;k}$ and $S^1$-weights of
 $\widehat{\cal S}_{i,j}/\widehat{\cal S}_{i,j-1}$.

\bigskip

\begin{flushleft}
{\bf An Euler sequence for $(T_{\ast}\HQuot_P({\cal E}^n))|_{E_{(A;0)}}$.}
\end{flushleft}
Recall the projection maps
  $\pi_1: E_{(A;0)}\times \CP^1\rightarrow E_{(A;0)}$ and
  $\pi_2: E_{(A;0)}\times \CP^1\rightarrow \CP^1$,
 the subsheaves
 $$
  0\; \hookrightarrow\; {\cal E}_1\;
      \stackrel{j_1}{\hookrightarrow}\; \cdots\;
      \stackrel{j_{I-1}}{\hookrightarrow}\; {\cal E}_I\;
      \stackrel{j_I}{\hookrightarrow}\;
   {\cal E}\;=\pi_2^{\ast}({\cal E}^n)
 $$
 and the quotient sheaves
 $$
  {\cal E}\;
     \stackrel{p_1}{\twoheadrightarrow}\;
      {\cal E}\!/\mbox{\raisebox{-.4ex}{${\cal E}_1$}}\;
     \stackrel{p_2}{\twoheadrightarrow}\; \cdots\;
     \stackrel{p_I}{\twoheadrightarrow}\;
      {\cal E}\!/\mbox{\raisebox{-.4ex}{${\cal E}_I$}}\;
    \rightarrow 0\,.
 $$
 on $E_{(A;0)}\times \CP^1$.
Let ${\cal K}$ be the kernel of the following morphism of
 ${\cal O}_{E_{(A;0)}\times\CPscriptsize^1}$-modules
 $$
   \oplus_{i=1}^I\,
     \Homsheaf(\,{\cal E}_i\,,\,
      {\cal E}\!/\mbox{\raisebox{-.4ex}{${\cal E}_i$}}\,)\;
    \longrightarrow\;
   \oplus_{i=1}^{I-1}
      \Homsheaf(\,{\cal E}_i\,,\,
       {\cal E}\!/\mbox{\raisebox{-.4ex}{${\cal E}_{i+1}$}}\,)
 $$
 defined by
 $$
  (\varphi_i)_{i=1}^I \;   \longmapsto\;
    (p_{i+1}\circ\varphi_i - \varphi_{i+1}\circ j_i)_{i=1}^{I-1}
 $$
 on each open subset $U$ of $E_{(A;0)}\times\CP^1$,
 cf.\ [CF1: Appendix].
Then

\bigskip

\noindent
{\bf Lemma 3.5.1
     [$(T_{\ast}\HQuot_P({\cal E}^n))|_{E_{(A;0)}}$ as push-forward].}
{\it
 $$
  \pi_{1\ast}{\cal K}=(T_{\ast}\HQuot_P({\cal E}^n))|_{E_{(A;0)}}\,.
 $$
} 

\bigskip

\noindent
{\it Proof.}
 In [CF1: Appendix], a fiberwise statement is given. Here we strengthen
  his result to a global statement.

 Let us outline first the approach of the proof.
 Let $\HQ$ be the stack associated to the hyper-Quot scheme
  $\HQuot_P({\cal E}^n)$.
 Then the tangent stack $T_{\ast}\HQ$ is the stackification of
  the prestack that associates to each affine ${\Bbb C}$-scheme $U$
  the groupoid $\HQ(U_{\varepsilon})$, where
  $U_{\varepsilon} := U\times_{\scriptsizeBbb C}{\Bbb C}[\varepsilon]$,
  with $\,\varepsilon^2=0$.
  \begin{figure}[htbp]
  \setcaption{{\sc Figure} 3-5-1.
   The functor that gives rise to the tangent space $T_{\ast}{\cal M}$
    of a moduli stack ${\cal M}$.
   A  morphism from $U$ (resp.\ $U_{\varepsilon}$) to the moduli stack
    ${\cal M}$ is the same as a flat family over $U$
    (resp.\ $U_{\varepsilon}$) of the objects ${\cal M}$ parameterizes
    (cf.\ the right third of the figure).
   \baselineskip 14pt
  } 
  \centerline{\psfig{figure=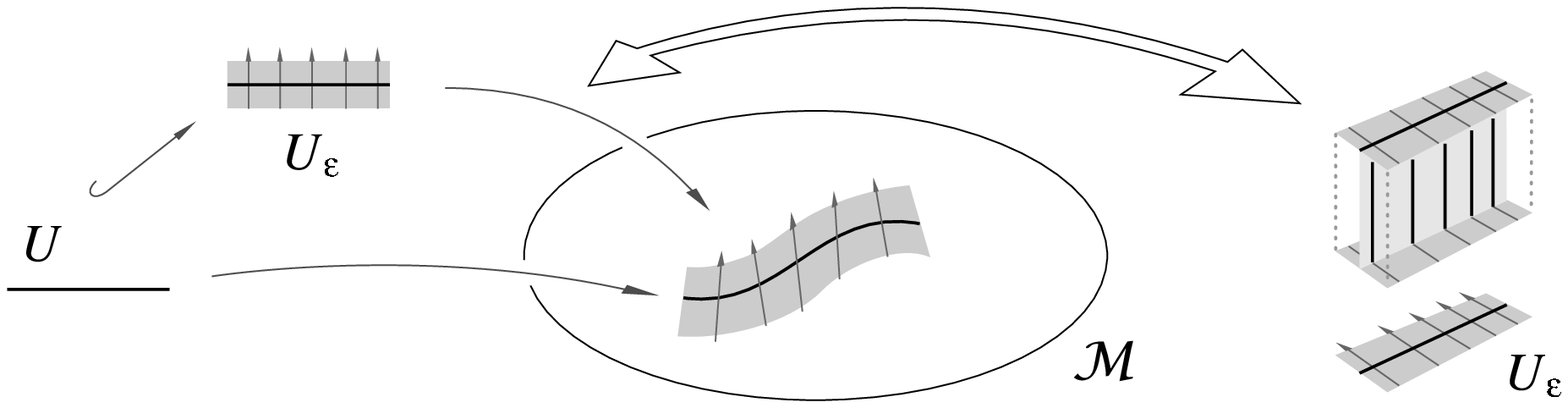,width=13cm,caption=}}
 \end{figure}
 $T_{\ast}\HQ$ is represented by the scheme associated to the tangent
   sheaf (a vector bundle in our case) $T_{\ast}\HQuot_P({\cal E}^n)$.
 Associated to the scheme $T_{\ast}\HQuot_P({\cal E}^n)|_{E_{(A;0)}}$
  is then the stack $T_{\ast}\HQ|_{E_{(A;0)}}$ over the subcategory of
  $E_{(A;0)}$-schemes from the stackification of the prestack that
  associates to each $U\rightarrow E_{(A;0)}$, where $U$ is affine,
  the groupoid $T_{\ast}\HQ(U) = \HQ(U_{\varepsilon})$.
 In our current case, $\HQuot_P({\cal E}^n)$ and hence $E_{(A;0)}$
  are projective. Thus $E_{(A;0)}$ can be covered by an atlas that
  consists of finitely many affine schemes such that any of their
  intersections are also affine.
 One shows first that the statement of the lemma holds
  for any affine open subset $U$ of $E_{(A;0)}$.
 Since all the isomorphisms of the groupoids $\HQ(U)$ and
  $\HQ(U_{\varepsilon})$ are identity maps, both $\HQ(U)$ and
  $\HQ(U_{\varepsilon})$ are sets canonically. Thus, to say that
  the lemma holds for $U$ means that
  $$
   (T_{\ast}\HQ|_{E_{(A;0)}})(U\hookrightarrow E_{(A;0)})\; \simeq\;
    H^0(\pi_1^{-1}(U), {\cal K}|_{\pi_1^{-1}(U)})
  $$
  as sets,
  where $(T_{\ast}\HQ|_{E_{(A;0)}})(U\hookrightarrow E_{(A;0)})$
  is the groupoid (set) of union of all groupoids
  $\HQ(U_{\varepsilon}\rightarrow E_{(A;0)})$ with
  $U_{\varepsilon}\rightarrow E_{(A;0)}$ extending the given inclusion
  $U\hookrightarrow E_{(A;0)}$.
 This set isomorphism will be constructed in a canonical/functorial
  way.
 Once this is achieved, then
  since the collection of groupoids
   $T_{\ast}\HQ|_{E_{(A;0)}}(U\hookrightarrow E_{(A;0)})$
   glue to give the stack $T_{\ast}\HQ|_{E_{(A;0)}}$ via the Isom-functor
   construction and the Grothendieck descent,
  the collection $\pi_{1\ast}({\cal K}|_{\pi_1^{-1}(U)})$ must glue
   to give the restriction of tangent bundle
   $(T_{\ast}\HQuot_P({\cal E}^n))|_{E_{(A;0)}}$, which represents
   the stack $T_{\ast}\HQ|_{E_{(A;0)}}$.

 Note that one may prove the Lemma first for the whole
  $\HQuot_P({\cal E}^n)$ and then discuss the restriction to $E_{(A;0)}$.
  But then one has to deal with the issue of commutativity of push-forward
  and restriction to a closed subscheme, which in general does not hold
  but has to be checked case by case.
 The above setting incorporates this issue into the discussion directly.

 \bigskip
 \noindent
 {\it Case $(a)\,:$ $X=\Gr_r({\Bbb C}^n)$.}
 In this case $\HQuot_P({\cal E}^n)$ is the Quot-scheme
  $\Quot_P({\cal E}^n)$.
 Let $U$ be an open affine subscheme of $E_{(A;0)}$.
 Since $E_{(A;0)}$ is smooth, we will assume that $U$ is smooth and
  quasi-projective.
 Let
   $U\stackrel{i}\hookrightarrow U_{\varepsilon}
                                      \stackrel{\pi}{\rightarrow} U$
   be the natural morphisms whose composition is the identity map on
   $U$.
 (The corresponding morphisms
    $U\times\CP^1 \rightarrow U_{\varepsilon}\times\CP^1
                  \rightarrow U\times\CP^1$
    will be denoted the same.)
  Let
    ${\cal V}$ be the tautological subbundle of ${\cal E}$ on
     $U\times\CP^1$,
    ${\cal E}_{\varepsilon} = \pi^{\ast}{\cal E}$ on
     $U_{\varepsilon}\times\CP^1$, and
    ${\cal V}^{\prime}$ be a subsheaf of ${\cal E}_{\varepsilon}$
     of Hilbert polynomial $P$ with its restriction to $U\times\CP^1$
     being ${\cal V}$.
  Then
   $\pi_{\ast}{\cal V}^{\prime}$ is a locally free subsheaf of
   $\pi_{\ast}{\cal E}_{\varepsilon}
    ={\cal E}\oplus {\cal E}\otimes\varepsilon$
    (canonically) on $U\times\CP^1$ with the Hilbert polynomial of the
     associated quotient sheaf being $2P$ and
   one has a canonical exact sequence
   $$
    0\;\longrightarrow\; {\cal V}\otimes{\varepsilon}\;
       \longrightarrow\; \pi_{\ast}{\cal V}^{\prime}\;
       \longrightarrow\; {\cal V}\; \longrightarrow \; 0\,.
   $$
  (This sequence splits non-canonically; thus
   $\pi_{\ast}{\cal V}^{\prime}
    \simeq {\cal V}\oplus {\cal V}\otimes\varepsilon$ non-canonically.)

  The above sequence together with projection of the locally free
   subsheaves $\pi_{\ast}{\cal V}^{\prime}$ of
   $\pi_{\ast}{\cal E}_{\varepsilon}$ into the direct summands,
   ${\cal E}$ and ${\cal E}\otimes\varepsilon$,
   of $\pi_{\ast}{\cal E}_{\varepsilon}$ induces the following diagram
   of canonical morphisms and isomorphisms
   $$
    \begin{array}{ccc}
     \pi_{\ast}{\cal V}^{\prime} & \longrightarrow
       & {\cal E}\otimes\varepsilon \simeq {\cal E} \\
     \downarrow   & & \downarrow  \\
    {\cal V}
      & &  ({\cal E}\otimes\varepsilon)/
              \mbox{\raisebox{-.4ex}{$({\cal V}\otimes\varepsilon)$}}
           \simeq {\cal E}\!/\mbox{\raisebox{-.4ex}{${\cal V}$}}
    \end{array}
   $$
   with both of the vertical arrows epimorphisms.
  Any local section $s$ of ${\cal V}$ can be lifted to a local section
   $s^{\prime}$ in
   $\pi_{\ast}{\cal V}^{\prime}$. The latter then maps to a local
    section in
    $\Hom({\cal V},
             {\cal E}\!/\mbox{\raisebox{-.4ex}{${\cal V}$}})$
    by following the above diagram.
  The image of $s$ in
   $\Hom({\cal V},
           {\cal E}\!/\mbox{\raisebox{-.4ex}{${\cal V}$}})$
   depends only on $s$, not on the choice of its lifting $s^{\prime}$.
  Thus, one obtains a canonical homomorphism
   $$
    \varphi_{{\cal V}^{\prime}}\;:\;
     {\cal V}\; \longrightarrow\;
     {\cal E}\!/\mbox{\raisebox{-.4ex}{${\cal V}$}}\,.
   $$
  The correspondence
   ${\cal V}^{\prime}\mapsto \varphi_{{\cal V}^{\prime}}$
   gives a map from
   $(T_{\ast}\HQ|_{E_{(A;0)}})(U\hookrightarrow E_{(A;0)})$
   to
   $\Hom({\cal V},
          {\cal E}\!/\mbox{\raisebox{-.4ex}{${\cal V}$}})
   $
   (which is $H^0(\pi_1^{-1}(U),{\cal K}|_{\pi_1^{-1}(U)})$
   in the current case).

  Conversely, given a
   $\varphi:{\cal V}\; \longrightarrow\;
     {\cal E}\!/\mbox{\raisebox{-.4ex}{${\cal V}$}}$,
   let ${\cal V}^{\prime\prime}_{\varphi}$
   be the locally free subsheaf of
    $\pi_{\ast}{\cal E}_{\varepsilon}
      ={\cal E}\oplus {\cal E}\otimes\varepsilon$,
   whose elements in fibers of ${\cal V}^{\prime\prime}_{\varphi}$
    are given by
    $$
    \{\,(v_0, v^{\prime})\,|\,
           v_0\in {\cal V}\,,\;
           v^{\prime}\; \mbox{is mapped to $\varphi(v_0)$ under
            ${\cal V}\otimes \varepsilon \rightarrow
               ({\cal E}\otimes\varepsilon)/
                 \mbox{\raisebox{-.4ex}{$({\cal V}\otimes\varepsilon)$}}
               \simeq {\cal E}\!/\mbox{\raisebox{-.4ex}{${\cal V}$}}
             $} \,\}\,.
    $$
  Then ${\cal V}^{\prime\prime}_{\varphi}$ fits into the exact sequence
   of ${\cal O}_{\pi_1^{-1}(U)}$-modules
   $$
    0\; \longrightarrow\; {\cal V}\otimes\varepsilon \;
        \longrightarrow\; {\cal V}^{\prime\prime}_{\varphi}\;
        \longrightarrow\; {\cal V}\; \longrightarrow\; 0
   $$
   and is invariant under the action of ${\varepsilon}$ on
    $\pi_{\ast}{\cal E}_{\varepsilon}$ induced from the multiplication
    of $\varepsilon$ on ${\cal E}_{\varepsilon}$.
  Thus
   ${\cal V}^{\prime\prime}_{\varphi}
                       =\pi_{\ast}{\cal V}^{\prime}_{\varphi}$
   for a unique locally free (and hence flat)
   ${\cal O}_{U_{\varepsilon}\times\CPscriptsize^1}$-submodule
   ${\cal V}^{\prime}_{\varphi}$ of ${\cal E}_{\varepsilon}$ on
   $U_{\varepsilon}\times\CP^1$.
  This gives a map from
   $\Hom({\cal V},
          {\cal E}\!/\mbox{\raisebox{-.4ex}{${\cal V}$}})$
   to  $(T_{\ast}\HQ|_{E_{(A;0)}})(U\hookrightarrow E_{(A;0)})$.
  One can check that the correspondences,
   ${\cal V}^{\prime}\mapsto\varphi_{{\cal V}^{\prime}}$ and
   $\varphi\mapsto{\cal V}^{\prime}_{\varphi}$,
   are inverse to each other.
  These constructions are canonical and functorial; thus
   $$
    (T_{\ast}\HQuot_P({\cal E}^n))|_{E_{(A;0)}}\;
     =\; \pi_{1\ast}
          \Homsheaf({\cal S},
             {\cal E}\!/\mbox{\raisebox{-.4ex}{${\cal S}$}})\,,
   $$
  where ${\cal S}$ is the tautological subsheaf of ${\cal E}$ on
   $E_{(A;0)}\times\CP^1$.

 \bigskip
 \noindent
 {\it Case $(b)\,:$ $X=\Fl_{r_1,\,\ldots,\,r_I}({\Bbb C}^n)$.}
 Repeat the same discussion for nested sequence of subsheaves
   over $U_{\varepsilon}\times\CP^1$ gives an embedding
  $$
   (T_{\ast}\HQ|_{E_{(A;0)}})(U\hookrightarrow E_{(A;0)})\;
    \hookrightarrow\;
    \oplus_{i=1}^I\, \Hom_{U\times\CPscriptsize^1}\,
                      (\,{\cal E}_i,
                  {\cal E}\!/\mbox{\raisebox{-.4ex}{${\cal E}_i$}}\,)\,.
  $$
  We shall now check that its image coincide with the set
  $H^0(U\times\CP^1,{\cal K}|_{U\times\CPscriptsize^1})$ as the subset
  $$
   \{\, (\varphi_i)_{i=1}^I \;|\;
     (p_{i+1}\circ\varphi_i-\varphi_{i+1}\circ j_i)_{i=1}^{I-1}\;=\;0\,\}
  $$
  of
  $\oplus_{i=1}^I\,
     \Hom(\,{\cal E}_i\,,\,
      {\cal E}\!/\mbox{\raisebox{-.4ex}{${\cal E}_i$}}\,)$.
 By induction, we only need to consider the case $I=2$.

 Let
  ${\cal E}_1^{\prime}\hookrightarrow {\cal E}_{\varepsilon}$
  (resp.\ ${\cal E}_2^{\prime}\hookrightarrow {\cal E}_{\varepsilon}$)
  be a flat subsheaf extension of ${\cal E}_1\hookrightarrow{\cal E}$
   (resp.\ ${\cal E}_2\hookrightarrow {\cal E}$) to
   $U_{\varepsilon}\times \CP^1$.
 Suppose that ${\cal E}_1^{\prime}$ is contained in
  ${\cal E}_2^{\prime}$.
 Then
  $\varphi_{{\cal E}_2^{\prime}}\circ j_1$ is the restriction of
   $\varphi_{{\cal E}_2^{\prime}}$ (defined on ${\cal E}_2$) to
   ${\cal E}_1$.
 In choosing the lifting sections $s^{\prime}$ in
  $\pi_{\ast}{\cal E}_2^{\prime}$ for local sections $s$ in
  ${\cal E}_1$ to define $\varphi_{{\cal E}_2^{\prime}}(s)$,
  one may choose $s^{\prime}$ a local section in
  $\pi_{\ast}{\cal E}_1\hookrightarrow \pi_{\ast}{\cal E}_2$
  since $\varphi_{{\cal E}_2^{\prime}}(s)$ is independent of the choice
  of liftings.
 Consequently, for $s$ a local section of ${\cal E}_1$
  $\varphi_{{\cal E}_2^{\prime}}(s) = \varphi_{{\cal E}_1^{\prime}}(s)$
  modulo ${\cal V}_2$, which is exactly
  $p_2\circ\varphi_{{\cal E}_1^{\prime}}(s)$.
 This shows that
  $\varphi_{{\cal E}_2^{\prime}}\circ j_1
               = p_2\circ \varphi_{{\cal E}_1^{\prime}}$
 and hence that
  $(T_{\ast}\HQ|_{E_{(A;0)}})(U\hookrightarrow E_{(A;0)})$
  embeds in
  $H^0(U\times\CP^1,{\cal K}|_{U\times\CPscriptsize^1})$.

 Conversely, suppose that
  $\varphi_{{\cal E}_2^{\prime}}\circ j_1
               = p_2\circ \varphi_{{\cal E}_1^{\prime}}$.
 From Part (a) of the proof, consider the canonical quotients
  $\delta_1:\pi_{\ast}{\cal E}_1^{\prime}\rightarrow {\cal E}_1$ and
  $\delta_2:\pi_{\ast}{\cal E}_2^{\prime}\rightarrow {\cal E}_2$.
 Treating all these locally free sheaves as vector bundles, then
  for a given $v\in {\cal E}_1\subset {\cal E}_2$,
  the assumption that
  $\varphi_{{\cal E}_2^{\prime}}\circ j_1
                 = p_2\circ \varphi_{{\cal E}_1^{\prime}}$.
  implies that the projection of $\delta_1^{-1}(v)$ in
  ${\cal E}\otimes\varepsilon$ is contained in the projection of
  $\delta_2^{-1}(v)$ in ${\cal E}\otimes\varepsilon$.
 Since the projection of $\delta_1^{-1}(v)$ (resp.\ $\delta_2^{-1}(v)$)
  to ${\cal E}\otimes\varepsilon$ is injective, this shows that
  $\delta_1^{-1}(v)\hookrightarrow \delta_2^{-1}(v)$ for all
  $v\in {\cal E}_1$ and hence
  that ${\cal E}_1^{\prime}$ is contained in ${\cal E}_2^{\prime}$.
 Together with the previous discussion, this proves that
  $$
   (T_{\ast}\HQ|_{E_{(A;0)}})(U\hookrightarrow E_{(A;0)})\;
   =\;
   H^0(U\times\CP^1,{\cal K}|_{U\times\CPscriptsize^1})\,.
  $$

 Consequently, the collection $\pi_{1\ast}({\cal K}|_{\pi_1^{-1}(U)})$
  glue to give $(T_{\ast}\HQuot_P({\cal E}^n))|_{E_{(A;0)}}$ and
 we conclude the proof.

\noindent\hspace{14cm}$\Box$

\bigskip

Recall from [CF1: Appendix] that there is a sheaf morphism
 $$
  \Homsheaf_{E_{(A;0)}\times\CPscriptsize^1}({\cal E},{\cal E})\;
   \longrightarrow\; \Homsheaf_{E_{(A;0)}\times\CPscriptsize^1}(
      {\cal E}_i,{\cal E}\!/\mbox{\raisebox{-.4ex}{${\cal E}_i$}})
 $$
 given by (denote ${\cal E}_i\hookrightarrow {\cal E}$ and
  ${\cal E}\rightarrow {\cal E}\!/\mbox{\raisebox{-.4ex}{${\cal E}_i$}}$
  also by $j_i$ and $p_i$ respectively)
 $$
  \psi\;\longmapsto\; ( p_i\circ\psi\circ j_i)_{i=1}^I
 $$
 that factors through ${\cal K}$ and is generically surjective
 (i.e.\ surjective at the stalk - or equivalently the fiber -
 at the generic point in the Zariski topology) onto ${\cal K}$.
An investigation of the non-surjectivity onto ${\cal K}$ of this
 morphism at the stalk at points on $E_{(A;0)}\times\{0\}$ motivates
 the following construction.

Consider the sheaf morphism
$$
 \Psi\,:\,
  \oplus_{i=1}^I
  \Homsheaf_{E_{(A;0)}\times\CPscriptsize^1}({\cal E}_i,{\cal E})\;
   \longrightarrow\;
    \oplus_{i=1}^I
     \Homsheaf_{E_{(A;0)}\times\CPscriptsize^1}({\cal E}_i,
     {\cal E}\!/\mbox{\raisebox{-.4ex}{${\cal E}_i$}})
$$
given by the following map of local sections on any open subset of
$E_{(A;0)}\times\CP^1$:
$$
 (\psi_i)_{i=1}^I\;\longmapsto\;(p_i\circ\psi_i)_{i=1}^I\,.
$$
Let ${\cal G}$ be the subsheaf of
 $\oplus_{i=1}^I
    \Homsheaf_{E_{(A;0)}\times\CPscriptsize^1}({\cal E}_i,{\cal E})$
 defined by the local sections $(\psi_i)_{i=1}^I$ such that
 the image sheaf $(\psi_i-\psi_{i+1}\circ j_i)({\cal E}_i)$
 of ${\cal E}_i$ in ${\cal E}$ lies in ${\cal E}_{i+1}$
 (on open subsets of $E_{(A;0)}\times\CP^1$).

\bigskip

\noindent
{\bf Lemma 3.5.2 [locally free resolution of ${\cal K}$].} {\it
 \begin{itemize}
  \item[$(1)$]
   The morphism $\Psi$ maps ${\cal G}$ surjectively onto ${\cal K}$.

  \item[$(2)$]
   ${\cal G}$ is locally free of rank
    $(r_1r_2+\,\cdots\,+ r_{I-1}r_I+ r_In)$.
   Along each $\CP^1$-fiber, ${\cal G}$ is a direct sum of
    non-negative line bundles.

  \item[$(3)$]
   The kernel of ${\cal G}\stackrel{\Psi}{\longrightarrow}{\cal K}$
    is given by $\oplus_{i=1}^I\Homsheaf({\cal E}_i,{\cal E}_i)$,
    which is locally free of rank $r_1^2+\,\cdots\,+r_I^2$.
 \end{itemize}
 Thus, in particular,
 $$
   0\;
    \longrightarrow\;\oplus_{i=1}^I\Homsheaf({\cal E}_i,{\cal E}_i)\;
    \longrightarrow\;{\cal G}\;\longrightarrow{\cal K}\;
    \longrightarrow\;0
 $$
 is a locally free resolution of ${\cal K}$.
} 

\bigskip

\noindent
{\it Proof.} Item (1) and Item (3) follow immediately by construction.
 For Item (2), we only need to check that ${\cal G}$ is locally free
  of the rank claimed along each $\CP^1$-fiber over a closed point of
  $E_{(A;0)}$. Since both $E_{(A;0)}\times\CP^1$ and $E_{(A;0)}$ are
  smooth and the rank is independent of the $\CP^1$-fibers, this then
  implies that ${\cal G}$ is locally free.

 Consider the sheaves ${\cal E}_i$ restricted to a $\CP^1$-fiber.
 Fix a realization $\CP^1=\Proj{\Bbb C}[z_0,z_1]$ that is compatible with
 the $S^1$-action with $0\in \CP^1$ corresponding to $[0:1]$.
 Recall from the proof of
  Lemma 2.1.4  
  the simultaneous $S^1$-weight subspace decomposition of an
  adjacent pair ${\cal E}_i\hookrightarrow {\cal E}_{i+1}$ on $\CP^1$.
 Incorporating these into presentation, one can identify ${\cal E}_i$
  and ${\cal E}_{i+1}$ as graded ${\Bbb C}[z_0,z_1]$-modules:
  $$
   {\cal E}_i\;
    =\; \left(\,\oplus_{j=1}^{r_i}\;
                z_0^{a_j}\!\cdot{\Bbb C}[z_0,z_1]\,
        \right)^{\sim}\,,\;\;
   {\cal E}_{i+1}\;
    =\; \left(\,
         \oplus_{j^{\prime}=1}^{r_{i+1}}\;
          z_0^{a^{\prime}_{j^{\prime}}}\!\cdot{\Bbb C}[z_0,z_1]\,
        \right)^{\sim}\,,\;\;
   {\cal E}^n\;=\; \left(\,{\Bbb C}[z_0,z_1]^{\oplus n}\,\right)^{\sim}
  $$
  such that
   $a_j\ge a^{\prime}_j$, $j=1,\,\ldots,\,r_i$, and that
   the inclusion ${\cal E}_i\hookrightarrow {\cal E}_{i+1}$
    is induced from the natural inclusions of graded modules
    $z_0^{a_j}\!\cdot{\Bbb C}[z_0,z_1]
     \hookrightarrow z_0^{a^{\prime}_j}\!\cdot{\Bbb C}[z_0,z_1]$,
    $j=1,\,\ldots,\, r_i$, from the identity map
    ${\Bbb C}[z_0,z_1]\hookrightarrow {\Bbb C}[z_0,z_1]$.
 In terms of these, the local sections $\psi_i$ and $\psi_{i+1}$
  of the Hom-sheaves $\Homsheaf({\cal E}_i,{\cal E})$ and
  $\Homsheaf({\cal E}_{i+1},{\cal E})$ are represented respectively as
  (degree-$0$ part of the localization of)
  ${\Bbb C}[z_0,z_1]$-valued matrices of the following block form:
  $$
   \psi_i\; =\; B_i
    \hspace{2em}\mbox{and}\hspace{2em}
   \psi_{i+1}\; =\; [B_{i+1},\ast\,]_{n\times r_{i+1}}\,,
  $$
  where both $B_i$ and $B_{i+1}$ are $n\times r_i$ matrices.
  $B_{i+1}$ corresponds to the composition $\psi_{i+1}\circ j_i$.
 Thus,
 $$
  (\psi_i-\psi_{i+1}\circ j_i)({\cal E}_i)\subset {\cal E}_{i+1}
   \hspace{2em}\Longleftrightarrow\hspace{2em}
  B_i-B_{i+1}\;
    =\; \left[
          \begin{array}{c}
            C_{i,i+1} \\
            0
          \end{array}
        \right]_{\,n\times r_i}\,,
 $$
 where
 $C_{i,i+1}=(c_{kl}(z_0,z_1))_{k,l}$ is an $r_{i+1}\times r_i$-matrix
 that satisfies
 $$
   \degree_{z_0}\,c_{kl}(z_0,z_1)\; \ge\; a^{\prime}_l-a_l\,,
 $$
Since, for a fixed $\psi_{i+1}$, the space of $C_{i,i+1}$ that satisfy
 the above condition is a free graded ${\Bbb C}[z_0,z_1]$-module of
 rank $r_ir_{i+1}$.
Let $i$ run from $1$ to $I$, this show that the restriction of ${\cal G}$
 to each $\CP^1$-fiber is locally free of rank
 $r_1r_2+\,\cdots\,+r_{I-1}r_I+ r_In$ and hecne the first half of
 Item (2).

Since $a^{\prime}_l-a_l\le 0$ for each $l=1,\,\ldots,\, r_i$, this proves
 the second half of Item (2).

This completes the proof.

\noindent\hspace{14cm}$\Box$

\bigskip

The above lemma, Item (2), implies that
$R^1\pi_{1\ast}{\cal G}=0$ by Grauert theorem,
 cf.\ [Ha: III.\ Corollary 12.9].
Consequently,

\bigskip

\noindent
{\bf Corollary 3.5.3 [Euler sequence of
 $(T_{\ast}\HQuot_P({\cal E}^n))|_{E_{(A;0)}}$].} {\it
 The restriction \newline
 $(T_{\ast}\HQuot_P({\cal E}^n))|_{E_{(A;0)}}=\pi_{1\ast}{\cal K}$
 fits into the following exact sequence
 {\small
  $$
   0\;\longrightarrow\;
      \pi_{1\ast}
        \left(\oplus_{i=1}^I\Homsheaf({\cal E}_i,{\cal E}_i)\right)\;
      \longrightarrow\; \pi_{1\ast}{\cal G}\;
      \longrightarrow\; \pi_{1\ast}{\cal K}\;
      \longrightarrow\;
      R^1\pi_{1\ast}
       \left(\oplus_{i=1}^I\Homsheaf({\cal E}_i,{\cal E}_i)\right)\;
      \longrightarrow\; 0\,.
  $$
 } 
} 

\bigskip

\begin{flushleft}
{\bf An Euler sequence for the vertical tangent bundle
     $T_{\ast}^{(vert,i)}E_{(A;0)}$.}
\end{flushleft}

\noindent
{\bf Lemma 3.5.4 [Euler sequence for
     $T_{\ast}\Fl_{k_1,\,\ldots,\,k_s}({\Bbb C}^n,\Pi_{\bullet})$].}
{\it
 Let $M:=\Fl_{k_1,\,\ldots,\,k_s}({\Bbb C}^n,\Pi_{\bullet})$ for notation,
  ${\cal E}={\cal O}_M\otimes{\Bbb C}^n$,
  $F_{\bullet}S: S_1\hookrightarrow \cdots \hookrightarrow S_s=:S$
   be the tautological filtration of subsheaves on
   $\Fl_{k_1,\,\ldots,\,k_s}({\Bbb C}^n,\Pi_{\bullet})$, and
  $F_{\bullet}{\cal E}$ be the filtration of ${\cal E}$ by
   ${\cal O}_M\otimes\Pi_{\bullet}$.
 Then $T_{\ast}\Fl_{k_1,\,\ldots,\,k_s}({\Bbb C}^n,\Pi_{\bullet})$
  fits into the following exact sequences of locally free
  ${\cal O}_M$-modules
  \begin{itemize}
   \item[$(1)$] $($compact form$):$
    {\small
    $$
     0\;\longrightarrow\;
      \Homsheaf_{{\cal O}_M}(F_{\bullet}S, F_{\bullet}S)\;
        \longrightarrow\;
      \Homsheaf_{{\cal O}_M}(F_{\bullet}S, F_{\bullet}{\cal E})\;
        \longrightarrow\;
      T_{\ast}\Fl_{k_1,\,\ldots,\,k_s}({\Bbb C}^n, \Pi_{\bullet})\;
        \longrightarrow\; 0\,;
    $$
    } 

  \vspace{-1em}
  \item[$(2)$] $($splitted form$):$
    $$
     0\;\longrightarrow\;
       \oplus_{j=1}^s\Homsheaf_{{\cal O}_M}(S_j,S_j)\;
        \longrightarrow\;
       {\cal H}\;
        \longrightarrow\;
       T_{\ast}\Fl_{k_1,\,\ldots,\,k_s}({\Bbb C}^n, \Pi_{\bullet})\;
        \longrightarrow\; 0\,,
    $$
   where ${\cal H}$ is a subsheaf of
   $\oplus_{j=1}^s\Homsheaf_{{\cal O}_M}(S_j,{\cal O}_M\otimes\Pi_j)$
   consisting of local sections  $(\psi_j)_{j=1}^s$ of sheaf morphisms
   such that the image sheaf
   $(\psi_j-\psi_{j+1}\circ\iota_j)(S_j)$ in ${\cal O}_M\otimes\Pi_j$
   lies in $S_{j+1}$.
  $($Here $\iota_j:S_j\rightarrow S_{j+1}$ is
   the inclusion morphism.$)$
 \end{itemize}
} 

\bigskip

\noindent
{\it Proof.}
 Recall the deformation of flags in ${\Bbb C}^n$ and the construction
  of the Euler sequence for the tangent bundle
  of the usual flag manifolds and that for
  $(T_{\ast}\HQuot_P({\cal E}^n))|_{E_{(A;0)}}$.
 Statement (1) follows by the case of ordinary flag manifolds but
  take into account the restriction that the subspace $S_j$
   are now restricted to move only in $\Pi_j$.
 (Note that the morphisms in the splitted form of the Euler sequence
  already apear in the construction of Euler sequence for
  $(T_{\ast}HQuot_P({\cal E}^n))|_{E_{(A;0)}}$. )

 Locally-freeness of the sheaf of modules involved in both Statement (1)
  and Statement (2) follows from the fact that all the sheaves involved
  are locally free and hence can be identified as vector bundles and that
  all the filtration involved are filtrations by subbundles.

 This concludes the proof.

\noindent\hspace{14cm}$\Box$

\bigskip

Recall the tower of fibrations of $E_{(A,0)}$
 $$
  E_{(A;0)}=E^{(1)}_{(A;0)}\;\longrightarrow\; \cdots\;
   \longrightarrow E_{(A;0)}^{(i)}\;\longrightarrow\;\cdots\;
   \longrightarrow\; E_{(A;0)}^{(I)}=\Fl_{m_{I,\bullet}}({\Bbb C}^n)
 $$
 with the fiber of $E_{(A;0)}^{(i)}\rightarrow E_{(A;0)}^{(i+1)}$
 being the restrictive flag manifold
 $\Fl_{m_{i,\bullet}}({\Bbb C}^n,\Pi_{i+1,\bullet})$.
Denote the vertical tangent bundle of the fibration
 $E_{(A;0)}^{(i)}\rightarrow E_{(A;0)}^{(i+1)}$
 by $T_{\ast}^{(vert)}E_{(A;0)}^{(i)}$
 (i.e.\ the subbundle of $T_{\ast}E_{(A;0)}^{(i)}$ consisting of
 tangent vectors along the fibers of the fibration $f_i$)
 and its pull-back to $E_{(A;0)}$ by
 $T_{\ast}^{(vert,i)}E_{(A;0)}$.
Recall also the various tautological subsheaves
  $\widehat{\cal S}_{i,\bullet}$ and $\widehat{\cal P}_{i+1,\bullet}$
  on $E_{(A;0)}$, $i=1,\,\ldots,\,I$.
Let
 $\widehat{\cal S}_i=\widehat{\cal S}_{i,K_i}$,
 $F_{\bullet}\widehat{\cal S}_i$ be the filtration
  $\widehat{\cal S}_{i,\bullet}$ of $\widehat{\cal S}_i$, and
 $F^{(i)}_{\bullet}{\cal E}$ be the filtration
  $\widehat{\cal P}_{i+1,\bullet}$
  of ${\cal E}={\cal O}_{E_{(A;0)}}\otimes{\Bbb C}^n$.
Then Lemma 3.5.4
 together with the locally-freeness of the sheaf of modules involved
 imply immediately the following$\,$:

\bigskip

\noindent
{\bf Corollary 3.5.5
     [Euler sequence for $T_{\ast}^{(vert,i)}E_{(A;0)}$].}
{\it
 The $i$-th vertical tangent bundle
  $T_{\ast}^{(vert,i)}E_{(A;0)}$ fits into the following
  exact sequence of ${\cal O}_{E_{(A;0)}}$-modules$\,$;
  \begin{itemize}
   \item[$(1)$] $($compact form$)\,:$
    {\small
    $$
     0\;\longrightarrow\;
      \Homsheaf(F_{\bullet}\widehat{\cal S}_i,
                F_{\bullet}\widehat{\cal S}_i)\;
        \longrightarrow\;
      \Homsheaf(F_{\bullet}\widehat{\cal S}_i,
                F^{(i)}_{\bullet}{\cal E})\;
        \longrightarrow\;
      T_{\ast}^{(vert,i)}E_{(A;0)}\; \longrightarrow\; 0\,;
    $$
    } 

  \vspace{-1em}
  \item[$(2)$] $($splitted form$)\,:$
    $$
     0\;\longrightarrow\;
       \oplus_{\j=1}^{K_i}\Homsheaf(\widehat{\cal S}_{i,j},
                                    \widehat{\cal S}_{i,j})\;
        \longrightarrow\; {\cal H}^{(i)}\;
        \longrightarrow\; T_{\ast}^{(vert,i)}E_{(A;0)}\;
        \longrightarrow\; 0\,,
    $$
   where ${\cal H}^{(i)}$ is a subsheaf of
   $\oplus_{j=1}^{K_i}\Homsheaf(\widehat{\cal S}_{i,j},
                                \widehat{\cal P}_{i+1,j})$
   consisting of local sections  $(\psi_j)_{j=1}^{K_i}$
   of sheaf morphisms such that the image sheaf
   $(\psi_j-\psi_{j+1}\circ\iota_j)(\widehat{\cal S}_{i,j})$ in
   $\widehat{\cal P}_{i+1,j}$ lies in $\widehat{\cal S}_{i,j+1}$.
  $($Here $\iota_j:\widehat{\cal S}_{i,j}
           \rightarrow \widehat{\cal S}_{i,j+1}$
   is the inclusion morphism.$)$
 \end{itemize}
} 

\bigskip

\begin{flushleft}
{\bf A decomposition of $\nu(E_{(A;0)}/\HQuot_P({\cal E}^n))$
     in the $K$-group of $E_{(A;0)}$.}
\end{flushleft}

\noindent
We give a decomposition of $T_{\ast}E_{(A;0)}$ and
 $(T_{\ast}\HQuot_P({\cal E}^n))|_{E_{(A;0)}}$
 in the $K$-group of $E_{(A;0)}$.
The decomposition of $\nu(E_{(A;0)}/HQuot_P({\cal E}))$ follows
 from
 $$
  [\nu(E_{(A;0)}/HQuot_P({\cal E}))]\;
  =\; [(T_{\ast}\HQuot_P({\cal E}^n))|_{E_{(A;0)}}]\,
       -\, [T_{\ast}E_{(A;0)}]\,.
 $$

\bigskip

\noindent $(a)$
{\it The $[T_{\ast}E_{(A;0)}]$-part}$\,$:
Recall that associated to a smooth morphism of smooth variety
 $f:X\rightarrow Y$ is the exact sequence of ${\cal O}_X$-modules
 $$
  0\;\longrightarrow\; T^{vert}_{\ast}(X/Y)\;
     \longrightarrow\; T_{\ast}X\; \longrightarrow f^{\ast}T_{\ast}Y\;
     \longrightarrow\; 0\,.
 $$
This gives rise to the decomposition
 $$
  [T_{\ast}X]\;=\; [T_{\ast}^{vert}(X/Y)] + [f^{\ast}T_{\ast}Y]
 $$
 in the $K$-group of $X$.
Since $f_i:E_{(A;0)}^{(i)}\rightarrow E_{(A;0)}^{(i+1)}$ is projective,
 $f_i^{\ast}$ is an exact functor on the category of coherent sheaves.
Thus one can employ the above identity iteratingly to the tower
 of fibrations of $E_{(A;0)}$ by restrictive flag manifolds.
Corollary 3.5.5 and the fact that all the graded sheaves of modules
 involved are locally free together imply the following decomposition
 of $T_{\ast}E_{(A;0)}$ in the $K$-group of $E_{(A;0)}\,$:

 {\footnotesize
 \begin{eqnarray*}
  [T_{\ast}E_{(A;0)}]
    & =
    & \sum_{i=1}^I\,[T_{\ast}^{(vert,i)}E_{(A;0)}]\,,
       \hspace{2em}\mbox{where
        $\;T_{\ast}^{(vert,I)}E_{(A;0)}
         :=\;$ the pull-back of
        $\;T_{\ast}\Fl_{m_{I,\bullet}}({\footnotesizeBbb C}^n)\,$,}
                                                               \\[.6ex]
    & = &  \sum_{i=1}^I\,[\Homsheaf(F_{\bullet}\widehat{\cal S}_i,
                           F^{(i)}_{\bullet}{\cal E})]\,
           -\, \sum_{i=1}^I\,
                [\Homsheaf(F_{\bullet}\widehat{\cal S}_i,
                           F_{\bullet}\widehat{\cal S}_i)]     \\[.6ex]
    & = & \sum_{i=1}^I\;\sum_{1\le j^{\prime} \le j \le K_i}\,
              \left[\Homsheaf\left( \widehat{\cal S}_{i,j}
                                   /\widehat{\cal S}_{i,j-1}\,,\,
                          \widehat{\cal P}_{i+1,j^{\prime}}
                                 /\widehat{\cal P}_{i+1,j^{\prime}-1}
                        \right)\right]    \\[.6ex]
    &   & \hspace{7em}
          -\,\sum_{i=1}^I\;\sum_{1 \le j^{\prime} \le j \le K_i}\,
              \left[\Homsheaf\left( \widehat{\cal S}_{i,j}
                                   /\widehat{\cal S}_{i,j-1}\,,\,
                          \widehat{\cal S}_{i,j^{\prime}}
                                   /\widehat{\cal S}_{i,j^{\prime}-1}
                        \right)\right]    \\[.6ex]
    & = & \sum_{i=1}^I\;\sum_{1 \le j^{\prime} \le j \le K_i}\,
              \left[\Homsheaf\left( \widehat{\cal S}_{i,j}
                                   /\widehat{\cal S}_{i,j-1}\,,\,
                        \widehat{\cal S}_{i+1, I_A(i,j^{\prime})}
                          /\widehat{\cal S}_{i+1, I_A(i,j^{\prime}-1)}
                        \right)\right]    \\[.6ex]
    &   & \hspace{7em}
          -\,\sum_{i=1}^I\;\sum_{1 \le j^{\prime} \le j \le K_i}\,
              \left[\Homsheaf\left( \widehat{\cal S}_{i,j}
                                   /\widehat{\cal S}_{i,j-1}\,,\,
                          \widehat{\cal S}_{i,j^{\prime}}
                                   /\widehat{\cal S}_{i,j^{\prime}-1}
                        \right)\right]    \\[.6ex]
    & = & \sum_{i=1}^I\;\sum_{1 \le j^{\prime} \le j \le K_i}\;
            \sum_{I_A(i,j^{\prime}-1)+1\le k \le I_A(i,j^{\prime})}\,
              \left[\Homsheaf\left( \widehat{\cal S}_{i,j}
                                   /\widehat{\cal S}_{i,j-1}\,,\,
                        \widehat{\cal S}_{i+1,k}
                          /\widehat{\cal S}_{i+1,k-1}
                        \right)\right]    \\[.6ex]
    &   & \hspace{7em}
          -\,\sum_{i=1}^I\;\sum_{1 \le j^{\prime} \le j \le K_i}\,
              \left[\Homsheaf\left( \widehat{\cal S}_{i,j}
                                   /\widehat{\cal S}_{i,j-1}\,,\,
                          \widehat{\cal S}_{i,j^{\prime}}
                                   /\widehat{\cal S}_{i,j^{\prime}-1}
                        \right)\right]\,.
 \end{eqnarray*}
 } 

\bigskip

\noindent $(b)$
{\it The $[(T_{\ast}\HQuot_P({\cal E}^n))|_{E_{(A;B)}}]$-part}$\,$:
Recall the exact sequences
 $$
  0\;\longrightarrow\; \oplus_{i=1}^I\Homsheaf({\cal E}_i,{\cal E}_i)\;
     \longrightarrow\; {\cal G}\; \longrightarrow\; {\cal K}\;
     \longrightarrow\; 0\,,
 $$
 {\small
 $$
  0\;\longrightarrow\;
     \pi_{1\ast}
       \left(\oplus_{i=1}^I\Homsheaf({\cal E}_i,{\cal E}_i)\right)\;
     \longrightarrow\; \pi_{1\ast}{\cal G}\;
     \longrightarrow\; \pi_{1\ast}{\cal K}\;
     \longrightarrow\;
     R^1\pi_{1\ast}
      \left(\oplus_{i=1}^I\Homsheaf({\cal E}_i,{\cal E}_i)\right)\;
     \longrightarrow\; 0\,,
 $$
 {\normalsize and}} 
 the filtered sheaves $F_{\bullet\,}{\cal E}_i:= {\cal E}_{i,\bullet}$
 with
 $$
  {\cal E}_{i,j}/\mbox{\raisebox{-.4ex}{${\cal E}_{i,j-1}$}}\;
  =\; \left(\,\pi_1^{\ast}\left(
          \widehat{\cal S}_{i,j}/
           \mbox{\raisebox{-.4ex}{$\widehat{\cal S}_{i,j-1}$}}\right)\,
     \right)
     (-a_{i,j}\,z)\,,
 $$
 which is locally free of rank $m_{i,j}$.
Recall also that if $f:V\rightarrow W$ is a proper morphism and
 ${\cal F}$ is a coherent sheaf on $V$ then the map
 $f_!({\cal F}):= \sum_i (-1)^i\,[R^if_{\ast}{\cal F}]$
 extends to a morphism of $K$-groups $f_!: K(V)\rightarrow K(W)$.

In the $K$-group of $\HQuot_P({\cal E}^n)\times\CP^1$
{\footnotesize
\begin{eqnarray*}
 \lefteqn{
 [\Homsheaf({\cal E}_i,{\cal E}_i)]\;
   =\; \sum_{j, j^{\prime}=1}^{K_i}\, \left[\,
         \Homsheaf\left(
           {\cal E}_{i,j}/\mbox{\raisebox{-.4ex}{${\cal E}_{i,j-1}$}}\;,\;
           {\cal E}_{i,j^{\prime}}/
               \mbox{\raisebox{-.4ex}{${\cal E}_{i,j^{\prime}-1}$}}\;
                  \right) \right]  } \\[.6ex]
  && =\;  \sum_{j,j^{\prime}=1}^{K_i}\, \left[\,
          \Homsheaf\left(
              (\,\pi_1^{\ast}(
              \widehat{\cal S}_{i,j}/
                \mbox{\raisebox{-.4ex}{$\widehat{\cal S}_{i,j-1}$}})\,
                             )(-a_{i,j}\,z)\;,\;
                (\,\pi_1^{\ast}(
                 \widehat{\cal S}_{i,j^{\prime}}/
                 \mbox{\raisebox{-.4ex}{$\widehat{\cal S}_{i,
                                              j^{\prime}-1}$}})\,
                 )(-a_{i,j^{\prime}}\,z)\,
                    \right)       \right]\,.
\end{eqnarray*}
{\normalsize By}} 
the definition of ${\cal G}$, ${\cal G}$ admits a filtration
 $F_{\bullet}{\cal G}$ with the associated graded coherent sheaf
 $\oplus_{i=1}^I\Homsheaf({\cal E}_i,{\cal E}_{i+1})$,
 where ${\cal E}_{I+1}$ is set to be ${\cal E}$.
Together with the tautological filtration of ${\cal E}_i$, this
 gives rise to the identity

{\footnotesize
\begin{eqnarray*}
 [{\cal G}]
  & = & \sum_{i=1}^I[\Homsheaf({\cal E}_i,{\cal E}_{i+1})]\;
    =\; \sum_{i=1}^I\;\sum_{\mbox{\tiny
                            $\begin{array}{c}
                               1\le j \le K_i \\[.6ex]
                               1\le j^{\prime} \le K_{i+1}
                             \end{array}$
                                  }  }
         \left[\,
          \Homsheaf\left(
            {\cal E}_{i,j}/\mbox{\raisebox{-.4ex}{${\cal E}_{i,j-1}$}}\;,\;
            {\cal E}_{i+1,j^{\prime}}/
                \mbox{\raisebox{-.4ex}{${\cal E}_{i+1,j^{\prime}-1}$}}\;
                   \right) \right]     \\[.6ex]
  & = & \sum_{i=1}^I\;\sum_{\mbox{\tiny
                            $\begin{array}{c}
                               1\le j \le K_i \\[.6ex]
                               1\le j^{\prime} \le K_{i+1}
                             \end{array}$
                                  }  }
         \left[\,
          \Homsheaf\left(
              (\,\pi_1^{\ast}(
              \widehat{\cal S}_{i,j}/
                \mbox{\raisebox{-.4ex}{$\widehat{\cal S}_{i,j-1}$}})\,
                             )(-a_{i,j}\,z)\;,\;
                (\,\pi_1^{\ast}(
                 \widehat{\cal S}_{i+1,j^{\prime}}/
                 \mbox{\raisebox{-.4ex}{$\widehat{\cal S}_{i+1,
                                              j^{\prime}-1}$}})\,
                 )(-a_{i+1,j^{\prime}}\,z)\,
                    \right)       \right]\,.
\end{eqnarray*}
} 

\noindent
(By convention,
 ${\cal E}_{I+1}={\cal E}$ with trivial filtration and
 $[\Homsheaf({\cal E}_I),{\cal E}_{I+1}]
   =\sum_{j=1}^{K_I}\,
     [ {\cal E}_{I,j}/\mbox{\raisebox{-.4ex}{${\cal E}_{I,j-1}$}}\,,\,
      {\cal E} ]
   =\sum_{j=1}^{K_I}
     [ (\,\pi_1^{\ast}(
              \widehat{\cal S}_{I,j}/
                \mbox{\raisebox{-.4ex}{$\widehat{\cal S}_{I,j-1}$}})\,
                             )(-a_{I,j}\,z)\,,\ ,{\cal E} ]$
 in the summation.)

Recall Lemma 3.5.2 and [CF1: Appendix]. It follows that
 $$
  H^1(\HQuot_P({\cal E})\times\CP^1, {\cal G})\;
  =\; H^1(\HQuot_P({\cal E})\times\CP^1, {\cal K})\; =\;0\,.
 $$
Consequently,
\begin{eqnarray*}
 \lefteqn{
   [(T_{\ast}HQuot_P({\cal E}^n))|_{E_{(A;0)}}]\;
    =\; \pi_{1!}\,[{\cal K}] \;
    =\; \pi_{1!}[{\cal G}]
           -\pi_{1!}\,\sum_{i=1}^I
                          [\Homsheaf({\cal E}_i,{\cal E}_i)] } \\[.6ex]
   && =\; [\pi_{1\ast}{\cal G}]\,
           -\,\sum_{i=1}^I
                [\pi_{1\ast}\Homsheaf({\cal E}_i,{\cal E}_i)]\,
           +\,\sum_{i=1}^I
                [R^1\pi_{1\ast}\Homsheaf({\cal E}_i,{\cal E}_i)]\,
\end{eqnarray*}
in the $K$-group of $E_{(A;0)}$.

Putting all these together, expressing $\Homsheaf$ of locally free
 sheaves by tensors with duals and applying the projection formula
 (e.g.\ [Ha: III, Exercise 8.3])$\,$:
  $$
    R^i\pi_{1\ast}\left(\,\pi_1^{\ast}(\,\cdot\,)
         \otimes {\cal O}(m)\,\right)\;
     =\; (\,\cdot\,)\,
         \otimes\, R^i\pi_{1\ast}\left({\cal O}(m)\right)\,,
  $$
 one leads to the decomposition

{\footnotesize
\begin{eqnarray*}
 \lefteqn{[(T_{\ast}HQuot_P({\cal E}^n))|_{E_{(A;0)}}] } \\[.6ex]
  && =\;
       \sum_{i=1}^I\!\!\sum_{\mbox{
                            $\begin{array}{c}
                               1\le j \le K_i \\[.6ex]
                               1\le j^{\prime} \le K_{i+1}
                             \end{array}$
                                  }  }\!\!
         \left[\,
          \left(
             \widehat{\cal S}_{i,j}/
                \mbox{\raisebox{-.4ex}{$\widehat{\cal S}_{i,j-1}$}}\,
          \right)^{\vee}
          \otimes\,
          \left(
             \widehat{\cal S}_{i+1,j^{\prime}}/
                 \mbox{\raisebox{-.4ex}{$\widehat{\cal S}_{i+1,
                                          j^{\prime}-1}$}}\right)\,
             \otimes\, \pi_{1\ast}{\cal O}(
                 a_{i,j}-a_{i+1,j^{\prime}})
                           \right]          \\[.6ex]
  && \hspace{2em}
     -\, \sum_{i=1}^I \sum_{j,j^{\prime}=1}^{K_i}\, \left[\,
            \left(
              \widehat{\cal S}_{i,j}/
                \mbox{\raisebox{-.4ex}{$\widehat{\cal S}_{i,j-1}$}}
                             \right)^{\vee}\,
               \otimes\,
              \left(
                 \widehat{\cal S}_{i,j^{\prime}}/
                 \mbox{\raisebox{-.4ex}{$\widehat{\cal S}_{i,
                                              j^{\prime}-1}$}}\right)\,
               \otimes\,
               \pi_{1\ast}
               {\cal O}(a_{i,j}-a_{i,j^{\prime}})\,
                                     \right]   \\[.6ex]
  && \hspace{2em}
     +\, \sum_{i=1}^I \sum_{j,j^{\prime}=1}^{K_i}\, \left[\,
            \left(
              \widehat{\cal S}_{i,j}/
                \mbox{\raisebox{-.4ex}{$\widehat{\cal S}_{i,j-1}$}}
                             \right)^{\vee}\,
               \otimes\,
              \left(
                 \widehat{\cal S}_{i,j^{\prime}}/
                 \mbox{\raisebox{-.4ex}{$\widehat{\cal S}_{i,
                                              j^{\prime}-1}$}}\right)\,
               \otimes\,
               R^1\pi_{1\ast}
               {\cal O}(a_{i,j}-a_{i,j^{\prime}})\,
                                     \right]\,.
\end{eqnarray*}
} 

Since $\pi_1:E_{(A;0)}\times\CP^1\rightarrow E_{(A;0)}$
 is the projection map,
 $$
  \pi_{1\ast}{\cal O}(a)\;
   =\; {\cal O}_{E_{(A;0)}}\otimes H^0(\CP^1,{\cal O}(a))\;
   =\;\left\{  \begin{array}{cl}
                 {\cal O}_{E_{(A;0)}}\otimes{\Bbb C}^{a+1}
                     & \mbox{for $a\ge 0$} \\[.6ex]
                 0   & \mbox{else}
               \end{array}
      \right.
  $$
  and
  $$
   R^1\pi_{1\ast}{\cal O}(a)\;
    =\; {\cal O}_{E_{(A;0)}}\otimes H^1(\CP^1,{\cal O}(a))
    =\; \left\{  \begin{array}{cl}
                 {\cal O}_{E_{(A;0)}}\otimes{\Bbb C}^{-a-1}
                     & \mbox{for $a\le -2$} \\[.6ex]
                 0   & \mbox{else}
               \end{array}
      \right.\,.
 $$
Observe also that, for a fixed $i=1,\,\ldots,\,I$,
the set of indices in ${\Bbb N}\times{\Bbb N}$
$$
 \{\,(j,j^{\prime})\,|\, 1\le j \le K_i,\, 1\le j^{\prime}\le K_{i+1},\,
          a_{i,j}-a_{i+1,j^{\prime}}\ge 0  \,\}\,
$$
coincides with the set of indices
$$
 \{\,(j,j^{\prime})\,|\,
  1\le j\le K_i,\,
  I_A(i,j^{\prime\prime}-1)+1\le j^{\prime}\le I_A(i,j^{\prime\prime})
   \hspace{1ex}\mbox{
     with $j^{\prime\prime}$ running over $[1,j]$}\,\}\,,
$$
(cf.\ Figure 2-1-2).
Incorporating these, one has the final expression

{\footnotesize
\begin{eqnarray*}
 \lefteqn{[(T_{\ast}HQuot_P({\cal E}^n))|_{E_{(A;0)}}] } \\[.6ex]
  && =\;
       \sum_{i=1}^I\;
        \sum_{1\le j^{\prime\prime}\le j \le K_i}\;
         \sum_{I_A(i,j^{\prime\prime}-1)
                     \le j^{\prime} \le I_A(i,j^{\prime\prime})}\;
         \left[\,
          \left(
             \widehat{\cal S}_{i,j}/
                \mbox{\raisebox{-.4ex}{$\widehat{\cal S}_{i,j-1}$}}\,
          \right)^{\vee}
          \otimes\,
          \left(
             \widehat{\cal S}_{i+1,j^{\prime}}/
                 \mbox{\raisebox{-.4ex}{$\widehat{\cal S}_{i+1,
                                          j^{\prime}-1}$}}\right)\,
             \otimes\, \pi_{1\ast}{\cal O}(
                 a_{i,j}-a_{i+1,j^{\prime}})
                           \right]          \\[.6ex]
  && \hspace{2em}
     -\, \sum_{i=1}^I\;
          \sum_{1\le j^{\prime}\le j\le K_i}\;
          \left[\,
            \left(
              \widehat{\cal S}_{i,j}/
                \mbox{\raisebox{-.4ex}{$\widehat{\cal S}_{i,j-1}$}}
                             \right)^{\vee}\,
               \otimes\,
              \left(
                 \widehat{\cal S}_{i,j^{\prime}}/
                 \mbox{\raisebox{-.4ex}{$\widehat{\cal S}_{i,
                                              j^{\prime}-1}$}}\right)\,
               \otimes\,
               \pi_{1\ast}
               {\cal O}(a_{i,j}-a_{i,j^{\prime}})\,
                                     \right]   \\[.6ex]
  && \hspace{2em}
     +\, \sum_{i=1}^I\;
          \sum_{1\le j<j^{\prime}\le K_i}
          \left[\,
            \left(
              \widehat{\cal S}_{i,j}/
                \mbox{\raisebox{-.4ex}{$\widehat{\cal S}_{i,j-1}$}}
                             \right)^{\vee}\,
               \otimes\,
              \left(
                 \widehat{\cal S}_{i,j^{\prime}}/
                 \mbox{\raisebox{-.4ex}{$\widehat{\cal S}_{i,
                                              j^{\prime}-1}$}}\right)\,
               \otimes\,
               R^1\pi_{1\ast}
               {\cal O}(a_{i,j}-a_{i,j^{\prime}})\,
                                     \right]\,.
\end{eqnarray*}
} 

\bigskip

\noindent $(c)$
{\it The decomposition of $\nu(E_{(A;0)}/HQuot_P({\cal E}^n))\,$}:
Combining Part (a) and Part (b), one obtains

{\footnotesize
\begin{eqnarray*}
 \lefteqn{
  [\nu(E_{(A;0)}/HQuot_P({\cal E}^n) )]     } \\[.6ex]
  && =\;
       \sum_{i=1}^I\;
        \sum_{1\le j^{\prime\prime}\le j \le K_i}\;
         \sum_{I_A(i,j^{\prime\prime}-1)
                     \le j^{\prime} \le I_A(i,j^{\prime\prime})}\;
         \left[\,
          \left(
             \widehat{\cal S}_{i,j}/
                \mbox{\raisebox{-.4ex}{$\widehat{\cal S}_{i,j-1}$}}\,
          \right)^{\vee}
          \otimes\,
          \left(
             \widehat{\cal S}_{i+1,j^{\prime}}/
                 \mbox{\raisebox{-.4ex}{$\widehat{\cal S}_{i+1,
                                          j^{\prime}-1}$}}\right)\,
             \otimes\, \pi_{1\ast}{\cal O}(
                 a_{i,j}-a_{i+1,j^{\prime}})
                           \right]          \\[.6ex]
  && \hspace{2em}
     -\, \sum_{i=1}^I\;
          \sum_{1\le j^{\prime}\le j\le K_i}\;
          \left[\,
            \left(
              \widehat{\cal S}_{i,j}/
                \mbox{\raisebox{-.4ex}{$\widehat{\cal S}_{i,j-1}$}}
                             \right)^{\vee}\,
               \otimes\,
              \left(
                 \widehat{\cal S}_{i,j^{\prime}}/
                 \mbox{\raisebox{-.4ex}{$\widehat{\cal S}_{i,
                                              j^{\prime}-1}$}}\right)\,
               \otimes\,
               \pi_{1\ast}
               {\cal O}(a_{i,j}-a_{i,j^{\prime}})\,
                                     \right]   \\[.6ex]
  && \hspace{2em}
     +\, \sum_{i=1}^I\;
          \sum_{1\le j<j^{\prime}\le K_i}
          \left[\,
            \left(
              \widehat{\cal S}_{i,j}/
                \mbox{\raisebox{-.4ex}{$\widehat{\cal S}_{i,j-1}$}}
                             \right)^{\vee}\,
               \otimes\,
              \left(
                 \widehat{\cal S}_{i,j^{\prime}}/
                 \mbox{\raisebox{-.4ex}{$\widehat{\cal S}_{i,
                                              j^{\prime}-1}$}}\right)\,
               \otimes\,
               R^1\pi_{1\ast}
               {\cal O}(a_{i,j}-a_{i,j^{\prime}})\,
                                     \right]  \\[.6ex]
  && \hspace{2em}
      -\, \sum_{i=1}^I\;\sum_{1 \le j^{\prime} \le j \le K_i}\;
            \sum_{I_A(i,j^{\prime}-1)+1\le k \le I_A(i,j^{\prime})}\,
              \left[\left( \widehat{\cal S}_{i,j}
                      /\widehat{\cal S}_{i,j-1}\right)^{\vee}\,
                 \otimes\,
                   \left(\widehat{\cal S}_{i+1,k}
                          /\widehat{\cal S}_{i+1,k-1}
                        \right)\right]    \\[.6ex]
  && \hspace{2em}
     +\,\sum_{i=1}^I\;\sum_{1 \le j^{\prime} \le j \le K_i}\,
            \left[\left( \widehat{\cal S}_{i,j}
                        /\widehat{\cal S}_{i,j-1}\right)^{\vee}\,
                 \otimes\,
                    \left(\widehat{\cal S}_{i,j^{\prime}}
                                   /\widehat{\cal S}_{i,j^{\prime}-1}
                        \right)\right]\,.
\end{eqnarray*}
} 

\bigskip

\begin{flushleft}
{\bf An exact expression of the $S^1$-equivariant Euler class
     $e_{S^1}(E_{(A;0)}/\HQuot_P({\cal E}^n))$.}
\end{flushleft}
Since $\widehat{\cal S}_{i,j}/\widehat{\cal S}_{i,j-1}$ are bundles
 on $E_{(A;0)}$ rather than on $E_{(A;0)}\times\CP^1$, the $S^1$-action
 on $\CP^1$ induces only the trivial action on them.
Thus in terms of Chern roots and $S^1$-weights
 $$
  c_{S^1}\left(\widehat{\cal S}_{i,j}/\widehat{\cal S}_{i,j-1}\right)\;
  =\; c\left(\widehat{\cal S}_{i,j}/\widehat{\cal S}_{i,j-1}\right)\;
  =\;  \prod_{k=1}^{m_{i,j}}(1+y_{i,j; k})\,.
 $$

The $S^1$-action on
 ${\cal E}_{i,j}/\mbox{\raisebox{-.4ex}{${\cal E}_{i,j-1}$}}$
 induces an $S^1$-action on their dual, tensor products, and
 also direct image sheaves of any of these:
 $\pi_{1\ast}(\,\cdot\,)$ and $R^1\pi_{1\ast}(\,\cdot\,)$
 on $E_{(A;0)}$.
This induced $S^1$-action on $\pi_{1\ast}(\,\cdot\,)$ and
 $R^1\pi_{1\ast}(\,\cdot\,)$ coincides with the $S^1$-action
 induced from that on the related $H^0(\CP^1,{\cal O}(m))$ and
 $H^1(\CP^1,{\cal O}(m))$ respectively.
The $S^1$-weight system for the latter can be computed directly
 by the \v{C}ech representation of sheaf cohomologies,
 e.g.\ [Ha: III.5]:
 $$
  H^0(\CP^1,{\cal O}(m))\,,\,m\ge 0\;
  :\; \{\,0,\,1,\,\cdots,\,m\,\}\,,
 $$
 represented by $1,\, z,\,\ldots,\,z^m$ on $U_0$,
 $$
  H^1(\CP^1,{\cal O}(m))\,,\, m\le -2\;
  :\; \{\,m+1,\,m+2,\, \ldots,\,-2,\,-1\,\}\,,
 $$
 represented by $z^{m+1},\,z^{m+2},\,\ldots,\,z^{-2},\,z^{-1}$
 on $U_0\cap U_{\infty}$,
and the $S^1$-weight system of the sheaf cohomology groups is
 the empty set for any other choice of $m$.
Denote the irreducible representation of $S^1=U(1)$ with weight
 $w$ by $\gamma^{(w)} (\simeq {\Bbb C})$ and
define
 $I^{\prime}_A(i,j)$ to be the maximal $l$, $1\le l\le K_{i+1}$
  such that $\widehat{\cal S}_{i,j}\subset \widehat{\cal S}_{i+1, l}$
  with $a_{i,j}\le a_{i+1,l}-1$
 and $I^{\prime\prime}_A(i,j)$ to be the minimal $l$ such that
  $a_{i,j}\le a_{i,l}-2$.
Then, after cancellation of identical terms, one can express
 $[\nu(E_{(A;0)}/HQuot_P({\cal E}^n))]$ as

{\footnotesize
\begin{eqnarray*}
 \lefteqn{
  [\nu(E_{(A;0)}/HQuot_P({\cal E}^n) )]     } \\[.6ex]
  && \hspace{-3em}
     =\;
       \sum_{i=1}^I\;
        \sum_{1\le j^{\prime\prime}\le j \le K_i}\;
         \sum_{I_A(i,j^{\prime\prime}-1)
                     \le j^{\prime} \le I^{\prime}_A(i,j^{\prime\prime})}\;
         \left[\,
          \left(
             \widehat{\cal S}_{i,j}/
                \mbox{\raisebox{-.4ex}{$\widehat{\cal S}_{i,j-1}$}}\,
          \right)^{\vee}
          \otimes\,
          \left(
             \widehat{\cal S}_{i+1,j^{\prime}}/
                 \mbox{\raisebox{-.4ex}{$\widehat{\cal S}_{i+1,
                                          j^{\prime}-1}$}}\right)\,
             \otimes\,
              \left(\,\gamma^{(1)}\,\oplus\,\cdots\,\oplus\,
                   \gamma^{(a_{i,j}-a_{i+1,j^{\prime}})}\,\right)
                           \right]          \\[.6ex]
  && \hspace{2em}
     -\, \sum_{i=1}^I\;
          \sum_{1\le j^{\prime} < j\le K_i}\;
          \left[\,
            \left(
              \widehat{\cal S}_{i,j}/
                \mbox{\raisebox{-.4ex}{$\widehat{\cal S}_{i,j-1}$}}
                             \right)^{\vee}\,
               \otimes\,
              \left(
                 \widehat{\cal S}_{i,j^{\prime}}/
                 \mbox{\raisebox{-.4ex}{$\widehat{\cal S}_{i,
                                              j^{\prime}-1}$}}\right)\,
               \otimes\,
                \left(\,
                 \gamma^{(1)}\,\oplus\,\cdots\,
                      \oplus\,\gamma^{(a_{i,j}-a_{i,j^{\prime}})}\,
                \right)
          \right]   \\[.6ex]
  && \hspace{2em}
     +\, \sum_{i=1}^I\;
          \sum_{1\le j
                  <I^{\prime\prime}_A(i,j)\le j^{\prime}\le K_i}
          \left[\,
            \left(
              \widehat{\cal S}_{i,j}/
                \mbox{\raisebox{-.4ex}{$\widehat{\cal S}_{i,j-1}$}}
                             \right)^{\vee}\,
               \otimes\,
              \left(
                 \widehat{\cal S}_{i,j^{\prime}}/
                 \mbox{\raisebox{-.4ex}{$\widehat{\cal S}_{i,
                                              j^{\prime}-1}$}}\right)\,
               \otimes\,
                \left(\,
                  \gamma^{(a_{i,j}-a_{i,j^{\prime}}+1)}\,
                     \oplus\,\cdots\,\oplus\, \gamma^{(-1)}\,
                \right)\,    \right]\,.
\end{eqnarray*}
} 

Let $\alpha=c_1({\cal O}_{\CPscriptsize^{\infty}}(1))$.
Putting all these together, applying the rule for Chern roots
 under tensor products and Lemma 3.3.2 in [L-L-L-Y],
 one concludes that

{\footnotesize
\begin{eqnarray*}
 \lefteqn{
  e_{S^1}(\nu(E_{(A;0)}/HQuot_P({\cal E}^n) ))     } \\[.6ex]
  && =\; \left(
       \prod_{i=1}^I\;
       \prod_{1\le j^{\prime\prime}\le j \le K_i}\;
       \prod_{I_A(i,j^{\prime\prime}-1)
              \le j^{\prime} \le I^{\prime}_A(i,j^{\prime\prime})}\;
       \prod_{k=1}^{m_{i,j}}\;
       \prod_{k^{\prime}=1}^{m_{i+1,j^{\prime}}}\;
       \prod_{l=1}^{a_{i,j}-a_{i+1,j^{\prime}}}\;
         (\,-\,y_{i,j;k}\,+\, y_{i+1,j^{\prime};k^{\prime}}\,
            -\, l\alpha \,)  \right)             \\[.6ex]
  && \hspace{2em}
     \cdot\,\left(
        \prod_{i=1}^I\;
        \prod_{1\le j^{\prime} < j\le K_i}\;
        \prod_{k=1}^{m_{i,j}}\;
        \prod_{k^{\prime}=1}^{m_{i,j^{\prime}}}\;
        \prod_{l=1}^{a_{i,j}-a_{i,j^{\prime}}}\;
        (\, -\, y_{i,j;k}\, +\, y_{i,j^{\prime};k^{\prime}}\,
            -\, l\alpha \,)             \right)^{-1}\\[.6ex]
  && \hspace{2em}
     \cdot\,\left(
       \prod_{i=1}^I\;
       \prod_{1\le j< I^{\prime\prime}_A(i,j)\le j^{\prime}\le K_i}\;
       \prod_{k=1}^{m_{i,j}}\;
       \prod_{k^{\prime}=1}^{m_{i,j^{\prime}}}\;
       \prod_{l=a_{i,j}-a_{i,j^{\prime}}+1}^{-1}\;
        (\,-y_{i,j;k}\, +\, y_{i,j^{\prime};k^{\prime}}\, -\,l\alpha\,)
                                     \right)  \\[.6ex]
  && =\; \left(
       \prod_{i=1}^I\;
       \prod_{1\le j^{\prime\prime}\le j \le K_i}\;
       \prod_{I_A(i,j^{\prime\prime}-1)
              \le j^{\prime} \le I^{\prime}_A(i,j^{\prime\prime})}\;
       \prod_{k=1}^{m_{i,j}}\;
       \prod_{k^{\prime}=1}^{m_{i+1,j^{\prime}}}\;
       \prod_{l=1}^{a_{i,j}-a_{i+1,j^{\prime}}}\;
         (\,-\,y_{i,j;k}\,+\, y_{i+1,j^{\prime};k^{\prime}}\,
            -\, l\alpha \,)  \right)             \\[.6ex]
  && \hspace{2em}
     \cdot\,\left(
        \prod_{i=1}^I\;
        \prod_{1\le j<j^{\prime} \le K_i}\;
        \prod_{k=1}^{m_{i,j}}\;
        \prod_{k^{\prime}=1}^{m_{i,j^{\prime}}}\;
          (-1)^{m_{i,j}m_{i,j^{\prime}}(a_{i,j^{\prime}}-a_{i,j}-1) }\;
            \left(\,\rule[-.1ex]{0ex}{3ex}
                     -\, y_{i,j^{\prime};k^{\prime}}\, +\, y_{i,j;k}\,
                     -\, (a_{i,j^{\prime}}-a_{i,j} )\,\alpha \,
            \right)\, \right)^{-1}\,.
\end{eqnarray*}
}  

\bigskip

\noindent
{\it Remark 3.5.6
     $[\,e_{S^1}(\nu(E_{(A;0)}/HQuot_P({\cal E}^n)))$ invertible$\,]$.}
 Observe that in the K-group decomposition of
  $\nu(E_{(A;0)}/HQuot_P({\cal E}^n))$ all the direct summands with
  null $S^1$-weight are cancelled.
 Consequencely, $e_{S^1}(\nu(E_{(A;0)}/HQuot_P({\cal E}^n)))$
  is an invertible element in $A^{\ast}(E_{(A;0)})(\alpha)$,
  as it should be and is manifest from the final exact expression
  above.

\bigskip

\noindent
{\it Remark 3.5.7 $[$Grassmannian manifold$\,]$.}
 For $X=\Gr_r({\Bbb C}^n)$, let $(\alpha_{\bullet};0)=(A;0)$,
  $K=K_1$, $m_j=m_{1,j}$, $a_j=a_{1,j}$, and $y_{j;k}=y_{i,j;k}$.
 Then the above expression simplifies to
  \begin{eqnarray*}
   \lefteqn{
     e_{S^1}(\nu(E_{(\alpha_{\bullet};0)}/\Quot_P({\cal E}^n))) }\\[.6ex]
     &&
      =\;\frac{ \rule[-1.2ex]{0ex}{1.6ex}
           \prod_{1\le j \le K}\;
           \prod_{k=1}^{m_j}\;
           \prod_{l=1}^{a_j}\;
             (\,-\,y_{j;k}\,-\, l\alpha \,)^n
             }{
           \prod_{1\le j<j^{\prime} \le K}\;
           \prod_{k=1}^{m_j}\;
           \prod_{k^{\prime}=1}^{m_{j^{\prime}}}\;
             (-1)^{m_j m_{j^{\prime}}(a_{j^{\prime}}-a_j-1) }\;
                (\,\rule[-.1ex]{0ex}{3ex}
                     -\, y_{i,j^{\prime};k^{\prime}}\, +\, y_{i,j;k}\,
                     -\, (a_{j^{\prime}}-a_j )\,\alpha \,)
             }
  \end{eqnarray*}
  in [B-CF-K].

\bigskip

\bigskip

\subsection{An exact computation of
    $\int_{X}\tau^{\ast}e^{H\cdot t}\cap {\mathbf 1}_d$.}

Recall the tower of fibrations of $E_{(A;0)}$ obtained by forgetting
 one by one the subsheaves in an inclusion sequence of subsheaves.
It fits into the following commutative diagram:
 $$
  \begin{array}{ccccccl}
   E_{(A;0)}=E^{(1)}_{(A;0)}
    & \stackrel{f_1}{\longrightarrow}\; \cdots\;
      \stackrel{f_{i-1}}{\longrightarrow}
    & E_{(A;0)}^{(i)}
    & \stackrel{f_i}{\longrightarrow}\;\cdots\;
      \stackrel{f_{I-1}}{\longrightarrow}
    & E_{(A;0)}^{(I)}=\Fl_{m_{I,\bullet}}({\Bbb C}^n)
    & \stackrel{f_I}{\longrightarrow}\;\pt          \\[.6ex]
   \downarrow\;\mbox{\scriptsize $p$}
    & \downarrow  & \downarrow  & \downarrow
    & \downarrow    & \hspace{2em}\|                \\[.6ex]
   \Fl_{r_1,\,\ldots,\,r_I}({\Bbb C}^n)
    & \longrightarrow\; \cdots\;\longrightarrow
    & \Fl_{r_i,\,\ldots,\,r_I}({\Bbb C}^n)
    & \longrightarrow\;\cdots\;\longrightarrow
    & \Gr_{r_I}({\Bbb C}^n)
    & \longrightarrow\; pt &.
  \end{array}
 $$
Each $f_i$ is a bundle map with fiber the restrictive flag manifold
 $$
  \Fl_{r_{i,1},\,\ldots,\,r_{i,K_i}}({\Bbb C}^n,\Pi_{i+1,\bullet})\;
  =\; \Fl_{r_{i,1},\,\ldots.r_{i,K_i}}(\Pi_{i+1,K_i},\Pi_{i+1,\bullet})\,.
 $$
To integrate a cohomology class over $E_{(A;0)}$ is the same as
 to push forward that class from $E_{(A;0)}$ to a class on a point.
In this section, we shall give an exact expression of this integral via
 a sequence of push-forwards following the above tower of fibrations.

\bigskip

\begin{flushleft}
{\bf The associated roof of the tower of fibrations of $E_{(A;0)}$.}
\end{flushleft}
Let ${\cal S}_i$ be the tautological subbundle on
 $\Fl_{r_i,\,\ldots,\,r_I}({\Bbb C}^n)$.
Its pull-back to $E_{(A;0)}^{(i)}$ will be denoted the same.
Since ${\cal P}_{i+1,K_i}={\cal S}_{i+1}$, the restrictive flag
 manifold bundle $f_i:E_{(A;0)}^{(i)}\rightarrow E_{(A;0)}^{(i+1)}$
 over $E_{(A;0)}^{(i+1)}$ is the one associated to the data:
 (1) inclusion sequence of subbundles of ${\cal S}_{i+1}\,$:
     ${\cal P}_{i+1,\bullet}\,:\,
     {\cal P}_{i+1,1}\, \hookrightarrow\,\cdots\,,\hookrightarrow\,
     {\cal P}_{i+1,K_i}\,\hookrightarrow\, {\cal S}_{i+1}$, and
 (2) sequence of integers: $0< r_{i,1} < \,\cdots\, < r_{i,K_i}\,$.
In the notation of Sec.\ 3.2,
 $f_i:E_{(A;0)}^{(i+1)}\rightarrow E_{(A;0)}^{i+1}$ is
 simply the bundle map
 $\Fl_{r_{i,1,\,\ldots,\,r_{i,K_i}}}
               ({\cal S}_{i+1},{\cal P}_{i+1,\bullet})
   \rightarrow E_{(A;0)}^{(i+1)}$.
Let
 $E^{\,\prime\,(i)}_{(A;0)}
     := \Fl_{r_{i,1},\,\ldots,\,r_{i,K_i}}({\cal S}_{i+1})$.
Then from Sec.\ 3.2 one has the following commutative diagram:
$$
 \begin{array}{ccccccccc}
  &&& E^{\,\prime\,(i-1)}_{(A;0)}  &&  E^{\,\prime\,(i)}_{(A;0)}
                                   &&& \\[.6ex]
  &&& \mbox{\scriptsize $\iota_{i-1}$}\;\nearrow \hspace{2em}
        \searrow\;\mbox{\scriptsize $f^{\,\prime}_{i-1}$}
  &&  \mbox{\scriptsize $\iota_i$}\;\nearrow\hspace{2em}
        \searrow\;\mbox{\scriptsize $f^{\,\prime}_i$}
                                    &&& \\
  \cdots
   & \stackrel{f_{i-2}}{\longrightarrow}    & E_{(A;0)}^{i-1}
   & \stackrel{f_{i-1}}{\longrightarrow}    & E_{(A;0)}^{(i)}
   & \stackrel{f_i}{\longrightarrow}        & E_{(A;0)}^{(i+1)}
   & \stackrel{f_{i+1}}{\longrightarrow}    & \cdots\;, \\
 \end{array}
$$
where
 $\iota_i:E_{(A;0)}^{(i)}\rightarrow E^{\,\prime\,(i)}_{(A;0)}$
  is the canonical embedding and
 $f^{\prime}_i:E^{\,\prime\,(i)}_{(A;0)}\rightarrow E_{(A;0)}^{(i+1)}$
  is the natural flag manifold bundle map.
We shall call the above diagram the {\it associated roof} of the tower
 of fibrations of $E_{(A;0)}$.

\bigskip

\begin{flushleft}
{\bf The push-forward/integration formula for $f_i$.}
\end{flushleft}
We now discuss an explicit form for each
 $f_{i\ast}:A^{\ast}(E_{(A;0)}^{(i)})
            \rightarrow A^{\ast}(E_{(A;0)}^{(i+1)})$
 that follows from Lemma 3.2.2 and [Br: Proposition 2.1].

\bigskip

\noindent
{\bf Fact 3.6.1 [push-forward formula for $f^{\prime}_i\,$].}
 ([Br: Proposition 2.1]; see also [B-CF-K] and [H-B-J: Chapter 4].)
{\it
 Recall the Chern roots $\{y_{i,j;k}\}_{k=1}^{m_{i,j}}$ of
  ${\cal S}_{i,j}/{\cal S}_{i,j-1}$, $1\le j\le K_i$.
 Let $\{\,y_{i+1;k}\,\}_{k=1}^{r_{i+1}-r_i}$ be the Chern roots
  of ${\cal S}_{i+1}/{\cal S}_i$.
 For notation uniformality, let $y_{i,K_i+1;k}:= y_{i+1;k}$
  and $m_{i,K_i+1}:=r_{i+1}-r_i$.
 Let
  $$
   P\; \in\;
   A^{\ast}(E_{(A;0)}^{(i+1)})\,\left[\,
     y_{i,j;k}\,|\,
       1\le j\le K_i+1\,;\;\;
       \mbox{for each $j$}\,,\; 1\le k\le m_{i,j}\,\right]
  $$
  represent a class in $A^{\ast}(E_{(A;0)}^{\,\prime\,(i)})$.
 Then
  $$
   f^{\,\prime}_{i\,\ast}\,P\;
    =\; \sum_{\overline{\sigma}
              \in Sym_{(i+1)}/
                   Sym_{(i,1)}\times\cdots\times Sym_{(i,K_i+1)}}\,
          \overline{\sigma}\,\cdot\,
          \left(\,
           \frac{P
                }{
            \prod_{1\le j<j^{\prime}\le K_i+1}\,
            \prod_{k=1}^{m_{i,j}}\,
            \prod_{k^{\prime}=1}^{m_{i,j^{\prime}}}\,
              (y_{i,j^{\prime};k^{\prime}}- y_{i,j;k} )
                 }\,
          \right)
  $$
  in $A^{\ast}(E_{(A;0)}^{(i+1)})$,
  where
   $\Sym_{(i+1)}$ is the permutation group of $r_{i+1}$-many letters,
    acting on the set  $\{\,y_{i,j;k}\,\}_{j,k}$, and
   $\Sym_{(i,j)}$ is the permutation group for the set
    $\{\,y_{i,j;k}\,\}_k$.
} 

\bigskip

Note that both the numerator and the denominator of the above fraction
 are invariant under
 the $\Sym_{(i,1)}\times\,\cdots\,\times\Sym_{(i,K_i+1)}$-action;
 thus the
 $\Sym_{(i+1)}/
    \Sym_{(i,1)}\times\,\cdots\,\times\Sym_{(i,K_i+1)}$-action
 on the fraction is well-defined.

\bigskip

\noindent
{\bf Corollary 3.6.2 [push-forward formula for $f_i$].} {\it
 Let
  $P \in \iota_i^{\ast}A^{\ast}(E_{(A;0)}^{\,\prime\,(i)})
     \subset A^{\ast}(E_{(A;0)}^{(i)})$
  be expressed in terms of the Chern roots as in Fact 3.6.1.
 Then
  $$
   f_{i\,\ast}\,P\;
    =\; \sum_{\overline{\sigma}
              \in Sym_{(i+1)}/
                   Sym_{(i,1)}\times\cdots\times Sym_{(i,K_i+1)}}\,
          \overline{\sigma}\,\cdot\,
          \left(\,
           \frac{P\,\cdot\,\Omega({\cal P}_{i+1,\bullet})
                }{
            \prod_{1\le j<j^{\prime}\le K_i+1}\,
            \prod_{j^{\prime\prime}=1}^{m_{i,j}}\,
            \prod_{
             j^{\prime\prime\prime}=1}^{m_{i,j^{\prime}}}\,
            (y_{i,j^{\prime};j^{\prime\prime\prime}}
              - y_{i,j;j^{\prime\prime}} )
                 }\,
          \right)
  $$
  in $A^{\ast}(E_{(A;0)}^{(i+1)})$,
  where
  $$
   \Omega({\cal P}_{i+1,\bullet})\;
   =\; \prod_{j=1}^{K_i}\,
       \prod_{k^{\prime}=1}^{r_{i+1}-l_{i+1,j}}\,
       \prod_{k^{\prime\prime}=1}^{m_{i,j}}\,
         (q_{i+1,j;k^{\prime}}-y_{i,j;k^{\prime\prime}} )\,,
  $$
  with $\{\,q_{i+1,j;k^{\prime}}\,\}_{k^{\prime}}$ being the set
  of Chern roots of ${\cal S}_{i+1}/{\cal P}_{i+1,j}$,
  is the Poincar\'{e} dual of $[E_{(A;0)}^{(i)}]$
  in $A^{\ast}(E_{(A;0)}^{\,\prime\,(i)})$,
  described in Lemma 3.2.2.
} 

\bigskip

\begin{flushleft}
{\bf An exact expression of the integral.}
\end{flushleft}
Recall
 $X=\Fl_{r_1,\,\ldots,\,r_I}({\Bbb C}^n)$,
 the $I$-tuples of integers $d=(d_1,\,\ldots,\,d_I)$,
 the $I$-tuple of hyperplane classes $H=(H_1,\,\ldots,\,H_I)$
  from the embedding
  $X \hookrightarrow
     \CP^{{n\choose r_1}-1}\times\,\cdots\,\times\CP^{{n\choose r_I}-1}$
  that generate $H^2(X,{\Bbb Z})$, and
 the associated $I$-tuple of K\"{a}hler parameters
  $t=(t_1,\,\ldots,\,t_I)$.
Denote
 $\Sym_{(i+1)}/\Sym_{(i,1)}\times\cdots\times \Sym_{(i,K_i+1)}$
 by $\overline{\Sym}_{(i+1)}$.
Applying Corollary 3.6.2  
 to the sequence of fibrations $f_i$ as a subfibration of
 $f_i^{\,\prime}$, one concludes that

\begin{eqnarray*}
 \lefteqn{
  \int_X\,\tau^{\ast}\,e^{H\cdot t}\cap {\mathbf 1}_d\;
  =\; \sum_A\,\int_{E_{(A;0)}}\,
       \frac{ g^{\ast}\,\psi^{\ast}\,e^{\kappa\cdot\zeta}
           }{ e_{S^1}(E_{(A;0)}/\HQuot_P({\cal E}^n))     } } \\[.6ex]
  &&
  =\; \sum_A\,
      f_{I\ast}\circ\,\cdots\,\circ f_{1\ast}\,
       \left(\,
        \frac{ e^{-\sum_{i=1}^I\,\zeta\,\cdot\,c_1({\cal S}_{i,K_i})}
            }{ e_{S^1}(E_{(A;0)}/\HQuot_P({\cal E}^n))     }\,
       \right)   \\[.6ex]
  & &
  =\; \sum_A\;\;\;
      \sum_{\overline{\sigma}_{I+1}\,\in \overline{Sym}_{(I+1)} }\;
        \overline{\sigma}_{I+1}\,\cdot\;\;\;
        \cdots\;\;\;
      \sum_{\overline{\sigma}_2\,\in \overline{Sym}_{(2)} }\,
        \overline{\sigma}_2\,\cdot
               \\[.6ex]
  & & \hspace{2em}
   \left(
    \frac{  e^{-\sum_{i=1}^I\,\zeta\,\cdot\,c_1({\cal S}_{i,K_i})}\,
            \cdot\,
            \prod_{i=1}^I\,\Omega\left({\cal P}_{i+1,\bullet}\right)
         }{ e_{S^1}(E_{(A;0)})/HQuot_P({\cal E}^n)\,\cdot\,
            \prod_{i=1}^I\,
            \prod_{1\le j<j^{\prime}\le K_i+1}\,
            \prod_{k=1}^{m_{i,j}}\,
            \prod_{
             k^{\prime}=1}^{m_{i,j^{\prime}}}\,
            (y_{i,j^{\prime};k^{\prime}}
              - y_{i,j;k} )
         }
   \right) \\[1ex]
  & &
  =\; \sum_A\;\;\;
      \sum_{\overline{\sigma}_{I+1}\,\in \overline{Sym}_{(I+1)} }\;
        \overline{\sigma}_{I+1}\,\cdot\;\;\;
        \cdots\;\;\;
      \sum_{\overline{\sigma}_2\,\in \overline{Sym}_{(2)} }\,
        \overline{\sigma}_2\,\cdot
               \\[.6ex]
  && \hspace{2em}
     \left[\; e^{-\sum_{i=1}^I\,
                (\, y_{i,1;1} + \,\cdots\,+ y_{i,1; m_{i,1}}\,
                 +\, \cdots\,
                 +\, y_{i,K_i;1} + \,\cdots\,+ y_{i,K_i; m_{i,K_i}}\,)\,
                \zeta_i}
               \rule{0em}{2.2em}\right. \\[.6ex]
  && \hspace{2em}
     \cdot\,
     \left(
       \prod_{i=1}^I\;
       \prod_{1\le j^{\prime\prime}\le j \le K_i}\;
       \prod_{I_A(i,j^{\prime\prime}-1)
              \le j^{\prime} \le I^{\prime}_A(i,j^{\prime\prime})}\;
       \prod_{k=1}^{m_{i,j}}\;
       \prod_{k^{\prime}=1}^{m_{i+1,j^{\prime}}}\;
       \prod_{l=1}^{a_{i,j}-a_{i+1,j^{\prime}}}\;
         (\,-\,y_{i,j;k}\,+\, y_{i+1,j^{\prime};k^{\prime}}\,
            -\, l\alpha \,)
     \right)^{-1}   \\[.6ex]
  && \hspace{2em}
     \cdot\,
       \left(
        \prod_{i=1}^I\;
        \prod_{1\le j<j^{\prime} \le K_i}\;
        \prod_{k=1}^{m_{i,j}}\;
        \prod_{k^{\prime}=1}^{m_{i,j^{\prime}}}\;
          (-1)^{m_{i,j}m_{i,j^{\prime}}(a_{i,j^{\prime}}-a_{i,j}-1) }\;
          \left(\,\rule[-.1ex]{0ex}{3ex}
                  -\, y_{i,j^{\prime};k^{\prime}}\, +\, y_{i,j;k}\,
                  -\, (a_{i,j^{\prime}}-a_{i,j} )\,\alpha \,
          \right)\,
       \right) \\[.6ex]
  && \hspace{2em}
     \cdot\,
       \left(\,
         \prod_{i=1}^I\,
         \prod_{j=1}^{K_i}\,
         \prod_{k^{\prime}=1}^{r_{i+1}-l_{i+1,j}}\,
         \prod_{k^{\prime\prime}=1}^{m_{i,j}}\,
           (q_{i+1,j;k^{\prime}}-y_{i,j;k^{\prime\prime}} )
       \right)\\[.6ex]
  && \hspace{2em}
     \left.
     \cdot\,
       \left(\,
        \prod_{i=1}^I\,
        \prod_{1\le j<j^{\prime}\le K_i+1}\,
        \prod_{k=1}^{m_{i,j}}\,
        \prod_{k^{\prime}=1}^{m_{i,j^{\prime}}}\,
               (\, y_{i,j^{\prime};k^{\prime}}\, -\, y_{i,j;k}\,)\,
       \right)^{-1}\,
     \right]\,.
\end{eqnarray*}

\bigskip

\noindent
{\it Remark 3.6.3 $[\,$hypergeometric series$\,]$.}
 Note that a fixed ${\Bbb T}^n$-action on ${\Bbb C}^n$ induces
  a ${\Bbb T}^n$-action on $E_{(A;0)}$ and a compatible
  ${\Bbb T}^n$-action on the total space all the bundles on
  $E_{(A;0)}$ whose Chern roots are involved above.
 Thus, once a ${\Bbb T}^n$-action on ${\Bbb C}^n$ is fixed, all
  our discussion has a ${\Bbb T}^n$-equivariant extension.
 In particular, the class $A_d$ is the non-equivariant limit of
  a ${\Bbb T}^n$-equivariant class.
 Recall from [L-L-Y1, II: Lemma 2.5] the fact that
  {\it the zero class $\omega=0$ is the only class in
  $H^{\ast}_{{\scriptsizeBbb T}^n}(X)$ such that
  $\int_X e^{H\cdot\zeta}\cap \omega=0$ for all generic
  $\zeta\in{\Bbb C}$}.
 This implies that the integral
  $\int_X\tau^{\ast}e^{H\cdot t}\cap {\mathbf 1}_d$
  determines the class ${\mathbf 1}_d$ in $H_{S^1}^{\ast}(X)$
  uniquely and, hence, the fundamental hypergeometric series
  $\HG[{\mathbf 1}]^X(t)$.

\bigskip

\bigskip

\section{Remarks on the Hori-Vafa conjecture.}

We conclude this paper with some remarks on the Hori-Vafa conjecture.

There are three aspects of stringy dualities that have led to
 various miraculous conjectural relations among mathematical objects
 and quantities:
 the string world-sheet field theory aspect,
 the string target space-time field theory aspect, and
 the lower-dimensional effective field theory aspect after
  compactifications.
One important example is the phenomenon of mirror symmetry of
 Calabi-Yau $3$-folds:
  world-sheet aspects from nonlinear sigma models that give rise to
   equivalent $d=2$, $N=(2,2)$ conformal field theories via
   a $U(1)$-charge redefinition and
  effective field theory aspects from compactification of $d=10$
   superstring theories to equivalent $d=4$, $N=2$ supersymmetric field
   theories via a field redefinition.
(Cf.\ Key word search: ``{\it duality}", ``{\it mirror}" from
      {\tt www.arXiv.org/hep-th})

In [H-V], Hori and Vafa generalize the world-sheet aspects of mirror
 symmetry to being the equivalence of $d=2$, $N=(2,2)$ supersymmetric
 field theories
 (i.e.\ without imposing the conformal invariance on the theory).
This leads them to a much broader encompassing picture of
 mirror symmetry.
(See [HKKPTVVZ] for full explanations.)
Putting this in the frame work of abelian gauged linear sigma models
 (GLSM) ([Wi1]),
 studying the effective field theories expanded around points in
  various phases on the theory space of a GLSM, and
 taking the generalized mirror of these theories enable them to
  link many $d=2$ field theories together.
Generalization of this setting to nonabelian GLSM ([Wi1: Sec.\ 5.3])
 leads them to the following conjecture,
 when the physical path integrals are interpreted appropriately
 mathematically:

\bigskip

\noindent
{\bf Conjecture 4.1 [Hori-Vafa].} [H-V: Appendix A]. {\it
 The hypergeometric series for a given homogeneous space
  $($e.g.\ a Grassmannian manifold$)$ can be reproduced from
  the hypergeometric series of simpler homogeneous spaces
  $($e.g.\ product of projective spaces$)$.
 Similarly for the twisted hypergeometric series that are related
  to the submanifolds in homogeneous spaces.
} 

\bigskip
\noindent
(Cf.\ [H-V: Appendix A]; see also [B-CF-K].)
In other words, different homogeneous spaces (or some simple quotients
 of them) can give rise to generalized mirror pairs.

\bigskip

\begin{flushleft}
{\bf The Hori-Vafa formula for Grassmannian manifolds.}
\end{flushleft}
For $X=\Gr_r({\Bbb C}^n)$,
 $E_{(A;0)}=E_{(\alpha_{\bullet};0)}\subset \Quot_{(d)}$
 is naturally isomorphic to a flag manifold
 $\Fl_{m_1,m_1+m_2,\,\ldots,\,m_1+m_2+\,\cdots\,+m_k}({\Bbb C}^n)$,
 where
  ${\alpha_{\bullet}}$ is a partition of $d$ into nonnegative integers
   of length $r$,
  $m_i$ counts the multiplicity of the identical summands
   in the partition,
  $m_1+\,\,\cdots\,+m_k=r$, (cf.\ [L-L-L-Y]).
The tower of fibrations of $E_{(A;0)}$ in the beginning of Sec.\ 3.6
 is shortened to
 $$
  \begin{array}{cccl}
   E_{(A;0)}=E_{(\alpha_{\bullet};0)}
    & \stackrel{f}{\longrightarrow}      & \pt    \\[.6ex]
   \downarrow\;\mbox{\scriptsize $p$}    & & \|   \\[.6ex]
   \Gr_r({\Bbb C}^n)  & \longrightarrow  & \pt  &,
  \end{array}
 $$
 in which both $p$ and $f$ are flag manifold bundle maps.
Applying push-formula Fact 3.6.1 to $p$ instead,
 plugging the exact expression of
  $e_{S^1}(\nu(E_{(\alpha_{\bullet})/\Quot_P({\cal E}^n)}))$
  in Remark 3.5.7,
 and simplifying,
one concludes that ([B-CF-K: Theorem 1.5])
 \begin{eqnarray*}
  \lefteqn{\HG[{\mathbf 1}]^X(t)\;
   =\; e^{-\frac{H\cdot t}{\alpha}}\,
       \sum_{d \ge 0}\,e^{d\cdot t}\,
       \sum_{(\alpha_{\bullet})}\,
           p_{\ast}\left(\,\frac{1}{
                    e_{S^1}(\nu(E_{(\alpha_{\bullet};0)})/
                                     \Quot_P({\cal E}^n)) }\,\right)
                                                            } \\[.6ex]
  & &
  =\; e^{-\frac{H\cdot t}{\alpha}}\,
      \sum_{d \ge 0}\,e^{d\cdot t}\;
      \sum_{(\alpha_{\bullet})}\;
       \frac{ \rule[-1.4ex]{0em}{1em}
          (-1)^{(r-1)d}\,
          \prod_{1\le j\le j^{\prime} \le r}
            (-y_j+y_{j^{\prime}} -(\alpha_j-\alpha_{j^{\prime}})\,\alpha)
           }{
           \prod_{1\le j\le j^{\prime} \le r}
                   (-y_j+y_{j^{\prime}})
           \prod_{j^{\prime\prime}=1}^r\,
           \prod_{l=1}^{\alpha_{ j^{\prime\prime} }}\,
              (-y_{j^{\prime\prime}}-l\alpha )
           }\,,
  \end{eqnarray*}
 where $\{\,y_j\,\}_j$ are the Chern roots of the tautological
 subbundle on $\Gr_r({\Bbb C}^n)$.
Apply this to the simplest Grassmannian manifolds,
 i.e.\ projective spaces, one concludes that

\bigskip

\noindent {\bf Corollary/Fact 4.2 [Hori-Vafa formula].} {\it
 The Hori-Vafa conjecture holds for Grassmannian manifolds.
 Explicitly, let
  $X=\Gr_r({\Bbb C}^n)$ $($with hyperplane class $H$$)$ and
  ${\Bbb P}=\prod_{i=1}^r\CP^{n-1}$
   $($with hyperplane class $H_i$ from the $i$-th factor$)$.
 Define
  $$
   \Delta=\prod_{j<j^{\prime}}(y_{j^{\prime}}-y_j)
    \hspace{1em}\mbox{and}\hspace{1em}
   {\cal D}_{\Delta}\;
    :=\; \prod_{j<j^{\prime}}\,
           (\, \alpha\,\frac{\partial}{\partial t_{j^{\prime}}}
               -\,\alpha\,\frac{\partial}{\partial t_j} \,)\,.
  $$
 Then
  $$
   \HG[{\mathbf 1}]^X(t)\;
    =\; e^{H(r-1)\pi\sqrt{-1}/\alpha}\;
        \frac{1}{\Delta}\;
        \left.
         \left(\,
          {\cal D}_{\Delta}\,
            \HG[{\mathbf 1}]^{\scriptsizeBbb P}(t_1,\,\ldots,\,t_r)\,
         \right)
        \right|_{t_i=t+(r-1)\pi\sqrt{-1}}\,.
  $$
} 

The above formula was derived in [B-CF-K] by using the method and
key results in [L-L-L-Y].

\bigskip

\begin{flushleft}
{\bf General Hori-Vafa formula.}
\end{flushleft}
For general flag manifold $X=\Fl_{r_1,\,\ldots,\,r_I}({\Bbb C}^n)$,
 one way to obtain a push-forward formula for the the natural bundle
 map $p:E_{(A;0)}\rightarrow X$ is to modify the tower of fibrations
 of $E_{(A;0)}$ so far used into another tower of fibrations that
 is compatible with the morphisms on $X$:
 $$
  E_{(A;0)}=\widehat{E}^{(1)}_{(A;0)}\;
    \stackrel{g_1}{\longrightarrow}\;
     \cdots\; \stackrel{g_{I-1}}{\longrightarrow}\;
    \widehat{E}^{(I)}_{(A;0)}\;   \stackrel{g_I}{\longrightarrow}\;
    \Fl_{r_1,\ldots,r_I}({\Bbb C}^n)\,.
 $$
Indeed there is a god-given one obtained by forgetting all but the last
 element in a flag that characterizes $S^1$-invariant subsheaves
 when forgetting the $S^1$-invariant subsheaves of an inlcusion sequence
 one by one.
Unfortunately the $g_i$ thus obtained is no longer a bundle map:
 though all the fibers of $g_i$ are restrictive flag manifolds,
 the topology of these restrictive manifolds can change.
Thus techniques need to develop to take care of this if one follows this
 line.

A second line is to work out a clean inversion formula from
 [L-L-Y1,II: Lemma 2.5] for the case of flag manifolds.
$\HG[{\mathbf 1}]^X(t)$ can then be reconstructed from the integrals
 $\int_X\tau^{\ast} e^{H\cdot t}\cap {\mathbf 1}_d$.
Further simplification via Young tableaux combinatorics may hope
 to lead to a clean exact form for $\HG[{\mathbf 1}]^X(t)$.

On the other hand, since $\HG[{\mathbf 1}]^X(t)$ and
 $\{\,\int_X\tau^{\ast}e^{H\cdot t}\cap{\mathbf 1}_d\,\}
                                            _{d\succcurlyeq 0}$
 carry informations differ only by the integral over $X$ and are
 obtainable from each other, one may re-write our expression of
 $\int_X\tau^{\ast} e^{H\cdot t}\cap {\mathbf 1}_d$
 more suggestively (but only formally) as
 \begin{eqnarray*}
  \lefteqn{
   \int_X\,\tau^{\ast}\,e^{H\cdot t}\cap {\mathbf 1}_d\;
   =\; \sum_A\,\int_{E_{(A;0)}}\,
       \frac{ g^{\ast}\,\psi^{\ast}\,e^{\kappa\cdot\zeta}
           }{ e_{S^1}(E_{(A;0)}/\HQuot_P({\cal E}^n))     } } \\[.6ex]
  & &
  =\; \widetilde{\prod}_{i=1}^I\;
      \left(\;
      \sum_{\overline{\sigma}_{i+1}\,\in \overline{Sym}_{(i+1)} }\;
        \overline{\sigma}_{i+1}\,\cdot \right.
               \\[.6ex]
  && \hspace{2em}
     \left[\; e^{-(\, y_{i,1;1} + \,\cdots\,+ y_{i,1; m_{i,1}}\,
                 +\, \cdots\,
                 +\, y_{i,K_i;1} + \,\cdots\,+ y_{i,K_i; m_{i,K_i}}\,)\,
                \zeta_i}
               \rule{0em}{2.2em}\right. \\[.6ex]
  && \hspace{2em}
     \cdot\,
     \left(
       \prod_{1\le j^{\prime\prime}\le j \le K_i}\;
       \prod_{I_A(i,j^{\prime\prime}-1)
              \le j^{\prime} \le I^{\prime}_A(i,j^{\prime\prime})}\;
       \prod_{k=1}^{m_{i,j}}\;
       \prod_{k^{\prime}=1}^{m_{i+1,j^{\prime}}}\;
       \prod_{l=1}^{a_{i,j}-a_{i+1,j^{\prime}}}\;
         (\,-\,y_{i,j;k}\,+\, y_{i+1,j^{\prime};k^{\prime}}\,
            -\, l\alpha \,)
     \right)^{-1}   \\[.6ex]
  && \hspace{2em}
     \cdot\,
       \left(
        \prod_{1\le j<j^{\prime} \le K_i}\;
        \prod_{k=1}^{m_{i,j}}\;
        \prod_{k^{\prime}=1}^{m_{i,j^{\prime}}}\;
          (-1)^{m_{i,j}m_{i,j^{\prime}}(a_{i,j^{\prime}}-a_{i,j}-1) }\;
          \left(\,\rule[-.1ex]{0ex}{3ex}
                  -\, y_{i,j^{\prime};k^{\prime}}\, +\, y_{i,j;k}\,
                  -\, (a_{i,j^{\prime}}-a_{i,j} )\,\alpha \,
          \right)\,
       \right) \\[.6ex]
  && \hspace{2em}
     \cdot\,
       \left(\,
         \prod_{j=1}^{K_i}\,
         \prod_{k^{\prime}=1}^{r_{i+1}-l_{i+1,j}}\,
         \prod_{k^{\prime\prime}=1}^{m_{i,j}}\,
           (q_{i+1,j;k^{\prime}}-y_{i,j;k^{\prime\prime}} )
       \right)\\[.6ex]
  && \hspace{2em}
     \left. \left.
     \cdot\,
       \left(\,
        \prod_{1\le j<j^{\prime}\le K_i+1}\,
        \prod_{k=1}^{m_{i,j}}\,
        \prod_{k^{\prime}=1}^{m_{i,j^{\prime}}}\,
               (\, y_{i,j^{\prime};k^{\prime}}\, -\, y_{i,j;k}\,)\,
       \right)^{-1}\,
     \right]\;\right)\,.
 \end{eqnarray*}
 where
 $\widetilde{\prod}_{i=1}^I$
 is a constrained product
 (i.e.\ not all summands in each factor can be picked out for
  multiplication when expanding the product; rather they have to
  satisfty specific admissible conditions determined by $A$).
The level structure indexed by $i$ suggests
 a version of ``{\it family Hori-Vafa formula}" generalizing
 the case of Grassmannian manifolds
while the appearance of Thom classes in the expression suggests
 a version of ``{\it quantum submanifold formula}" generalizing
 the case of complete intersection submanifolds.

Finally, even if all these technicalities are settled and
 Hori-Vafa conjecture is checked, there is still a final question:
 {\it Why do they go this way?}
For that one has to turn to the most fundamental understanding of
 quantum field theories and path integrals.

\bigskip

With all these outlooks - and amazement and puzzles as well -,
 we conclude this paper temporarily here.

\newpage


{\footnotesize
\begin{thebibliography}{AAAAaa}
%




\bibitem[Br]{} M.\ Brion,
 {\it The push-forward and Todd class of flag bundles},
 in {\sl Parameter spaces
         - Enumerative geometry, algebra and combinatorics},
 P.\ Pragacz ed., Banach Center Publ.\ 36, pp.~45~-~50,
 Polish Acad.\ Sci.\ Warsaw, 1996.


\bibitem[B-CF-K]{} A.\ Bertram, I.\ Ciocan-Fontanine and B.\ Kim,
 {\it Two proofs of a conjecture of Hori and Vafa},
 math.AG/0304403.


\bibitem[B-H]{} A.\ Borel and F.\ Hirzebruch,
 {\it Characteristic classes and homogeneous spaces,}
 {\it I.} {\sl Amer.\ J.\ Math.}\ {\bf 80} (1958), pp.~458~-~538;
 {\it II.} {\sl Amer.\ J.\ Math.}\ {\bf 81} (1959), pp.~315~-~382;
 {\it III.} {\sl Amer.\ J.\ Math.}\ {\bf 82} (1960), pp.~491~-~504.

\bibitem[B-T]{} R.\ Bott and L.W.\ Tu,
 {\sl Differential forms in algebraic topology},
 GTM 82, Springer-Verlag, 1982.

\bibitem[CF1]{} I.\ Ciocan-Fontanine,
 {\it The quantum cohomology ring of flag varieties},
 {\sl Transactions Amer.\ Math.\ Soc.}\ {\bf 351} (1999),
 pp.~2695~-~2729.

\bibitem[CF2]{} --------,
 {\it On quantum cohomology rings of partial flag varieties},
 {\sl Duke Math.\ J.}\ {\bf 98} (1999), pp.~485~-~524.

\bibitem[Ch1]{} L.\ Chen,
 {\it Poincar\'{e} polynomials of hyperquot schemes},
 {\sl Math.\ Ann.}\ {\bf 321} (2001), pp.~235~-~251.

\bibitem[Ch2]{} --------,
 {\it Quantum cohomology of flag manifolds},
 {\sl Adv.\ Math.}\ {\bf 174} (2003), pp.~1~-~34.




\bibitem[C-V]{} S.\ Cecotti and C.\ Vafa,
 {\it On the classification of $N=2$ supersymmetric theories},
 {\sl Commun.\ Math.}\ {\bf 158} (1993), pp.~569~-~644.



\bibitem[E-H]{} D.\ Eisenbud and J. Harris,
 {\sl The geometry of schemes},
 GTM 197, Springer-Verlag, 2000.

\bibitem[E-H-X1]{} T.\ Eguchi, K.\ Hori, and C.-S.\ Xiong,
 {\it Gravitational quantum cohomology},
 {\sl Int.\ J.\ Mod.\ Phys.}\ {\bf A12} (1997), pp.~1743~-~1782.

\bibitem[E-H-X2]{} --------,
 {\it Quantum cohomology and Virasoro algebra},
 {\sl Phys.\ Lett.}\ {\bf B402} (1997), pp.~71~-~80.


\bibitem[Fu1]{} W.\ Fulton,
 {\sl Intersection theory},
 Ser.\ Mod.\ Surv.\ Math.\ 2, Springer-Verlag, 1984.

\bibitem[Fu2]{} --------,
 {\it Flags, Schubert polynomials, degeneracy loci, and determinantal
      formulas},
 {\sl Duke Math.\ J.}\ {\bf 65} (1991), pp.~381~-~420.

\bibitem[Fu3]{} --------,
 {\sl Young tableaux - with aplications to representation theory and
      geometry},
 London Math.\ Soc.\ Student Texts 35, Cambridge Univ.\ Press, 1997.

\bibitem[F-P]{} W.\ Fulton and R.\ Pandharipande,
 {\it Notes on stable maps and quantum cohomology},
 in {\sl Algebraic geometry - Stanta Cruz 1995},
 J.\ Koll\'{a}r, R.\ Lazarsfeld, and D.\ Morrison eds.,
 Proc.\ Symp.\ Pure Math.\ vol.\ 62, part 2, pp.\ 45 - 96,
 Amer.\ Math.\ Soc.\ 1997.

\bibitem[Ha]{} R.\ Hartshorne,
 {\sl Algebraic geometry},
 GTM 52, Springer-Verlag, 1977.

\bibitem[Hi]{} F.\ Hirzebruch,
 {\sl Topological methods in algebraic geometry},
 Grund.\ Math.\ Wiss.\ Ein.\ 131, Springer-Verlag, 1966.

\bibitem[H-B-J]{} F.\ Hirzebruch, T.\ Berger, and R.\ Jung,
 {\sl Manifolds and modular forms},
 translated by P.S.\ Landweber, 2nd ed.,
 Aspects Math.\ 20, Vieweg, 1994.

\bibitem[H-I-V]{} K.\ Hori, A.\ Iqbal, and C.\ Vafa,
 {\it D-branes and mirror symmetry},
 {\tt hep-th/0005247}.

\bibitem[HKKPTVVZ]{} K.\ Hori, S.\ Katz, A.\ Klemm, R.\ Pandharihande,
     R.\ Thomas, C.\ Vafa, R.\ Vakil, E.\ Zaslow,
 {\sl Mirror symmetry},
 Clay Math.\ Inst.\ Mono.\ 1, Amer.\ Math.\ Soc.\ 2003.

\bibitem[H-L]{} D.\ Huybrechts and M.\ Lehn,
 {\sl The geometry of moduli spaces of sheaves},
 Vieweg, 1997.

\bibitem[H-V]{} H.\ Hori and C.\ Vafa,
 {\it Mirror symmetry},
 {\tt hep-th/0002222}.

\bibitem[Il]{} L.\ Illusie,
 {\sl Complexe cotangent et d\'{e}formations, I},
 Lect.\ Notes Math.\ 239, Springer-Verlag, 1971.

\bibitem[Kim]{} B.\ Kim,
 {\sl Gromov-Witten invariants for flag manifolds},
 Ph.D.\ dissertation, University of California at Berkeley, 1996.

\bibitem[Ko]{} J.\ Koll\'{a}r,
 {\sl Rational curves on algebraic varieties},
 Ser.\ Mod.\ Surv.\ Math.\ 32, Springer-Verlag, 1996.

\bibitem[La]{} G.\ Laumon,
 {\it Un analogue global de c\^{o}ne nilpotent},
 Duke Math.\ J.\ {\bf 57} (1988), pp.~647~-~671.

\bibitem[Le]{} M.\ Lehn,
 {\it On the cotangent sheaf of Quot-schemes},
 {\sl International J.\ Math.}\ {\bf 9} (1998), pp.~513~-~522.



\bibitem[Liu]{} K.\ Liu,
 {\it Heat kernel and moduli space},
 {\sl Math.\ Res.\ Lett.}\ {\bf 3} (1996), pp.~743~-~762.

\bibitem[L-L-L-Y]{}
 B.\ Lian, C.-H.\ Liu, K.\ Liu, and S.-T. Yau,
 {\it The $S^1$ fixed points in Quot-schemes and mirror principle
      computations},
 in {\sl Vector bundles and representations theory},
 S.D.\ Cutkosky, D.\ Edidin, Z.\ Qin, and Q.\ Zhang eds.,
 Contemp.\ Math.\ 322, pp.~165~-~194, Amer.\ Math.\ Soc., 2003.

\bibitem[L-L-Y1]{}
 B.\ Lian, K.\ Liu, and S.-T.\ Yau,
 {\it Mirror principle, I},
 {\sl Asian J.\ Math.}\ {\bf 1} (1997), pp.~729~-~763;
 {\it II}, {\sl Asian J.\ Math.}\ {\bf 3} (1999), pp.~109~-~146;
 {\it III}, {\sl Asian J.\ Math.}\ {\bf 3} (1999), pp.~771~-~800;
 {\it IV}, {\tt math.AG/0007104}.

\bibitem[L-L-Y2]{} --------,
 {\it A survey of mirror principle},
 {\tt math.AG/0010064}.

\bibitem[L(CH)-L-Y]{} C.-H.\ Liu, K.\ Liu, and S.-T.\ Yau,
 {\it On A-twisted moduli stack for curves from Witten's gauged
      linear sigma models},
 {\tt math.AG/0212316}.

\bibitem[L-MB]{} G.\ Laumon and L.\ Moret-Bailly,
 {\sl Champs alg\'{e}briques},
 Springer-Verlag, 2000.

\bibitem[L-T1]{} J.\ Li and G.\ Tian,
 {\it Quantum cohomology of homogeneous varieties},
 {\sl J.\ Alg.\ Geom.}\ {\bf 6} (1997), pp.~267~-~305.

\bibitem[L-T2]{} --------,
 {\it Virtual moduli cycles and Gromov-Witten invariants of
      algebraic varieties},
 {\sl J.\ Amer.\ Math.\ Soc.}\ {\bf 11} (1998), pp.~119~-~174.

\bibitem[Ma]{} I.G.\ Macdonald,
 {\sl Symmetric functions and Hall polynomials},
 Oxford Univ.\ Press, 1979.

\bibitem[Mu]{} D.\ Mumford,
 {\sl Lectures on curves on an algebraic surface},
 Ann.\ Math.\ Studies 59, Princeton Univ.\ Press, 1966.

\bibitem[M-M]{} B.M.\ Mann and R.J.\ Milgram,
 {\it On the moduli space of $\SU(n)$ monopoles and holomorphic
      maps to flag manifolds},
 {\sl J.\ Diff.\ Geom.}\ {\bf 38} (1993), pp.~39~-~103.


\bibitem[So]{} F.\ Sottile,
 {\it Real rational curves in Grassmannians},
 {\sl J.\ Amer.\ Math.\ Soc.}\ {\bf 13} (2000), pp.~333~-~341.

\bibitem[Str]{} S.A.\ Str{\o}mme,
 {\it On parametrized rational curves in Grassmann varieties},
 in {\sl Space curves}, F.\ Ghione, C.\ Peskine, and E.\ Sernesi eds,
 pp.~251~-~272, Lect.\ Notes Math.\ 1266, Springer-Verlag, 1987.

\bibitem[S-S]{} F.\ Sottile and B.\ Sturmfels,
 {\it A Sagbi basis for the quantum Grassmannian},
 {\sl J.\ Pure Appl.\ Algebra.}\ {\bf 158} (2001), pp.~347~-~366.

\bibitem[Wi1]{} E.\ Witten,
 {\it Phases of $N=2$ theories in two dimensions},
 {\sl Nucl.\ Phys.}\ {\bf B403} (1993), pp.~159~-~222.

\bibitem[Wi2]{} --------,
 {\it The Verlinde algebra and the cohomology of the Grassmannian},
 in {\sl Geometry, topology, and physics for Raoul Bott},
 S.-T.\ Yau ed., pp.~357~-~422, International Press, 1994.

%
\end{thebibliography}
}
\end{document}